\documentclass{amsart}
\usepackage{fancyvrb}
\usepackage{amsmath}
\usepackage{amsfonts}
\usepackage{amssymb}
\usepackage{amsthm}
\usepackage{cases}
\usepackage{caption}
\usepackage{xcolor}
\usepackage{tikz}
\usepackage{centernot}
\usepackage{mathtools}
\usepackage[margin=1in]{geometry}
\usepackage[makeroom]{cancel}
\usepackage{chemarrow}
\usepackage{url}
\usetikzlibrary{decorations.pathreplacing, decorations.pathmorphing, math}

\newtheorem{lemma}{Lemma}[section]
\newtheorem{corollary}[lemma]{Corollary}
 
\newtheorem{proposition}[lemma]{Proposition} 
\newtheorem{theorem}[lemma]{Theorem}

\theoremstyle{definition}
\newtheorem{definition}[lemma]{Definition}
\newtheorem{example}[lemma]{Example}
\newtheorem{method}[lemma]{Method}
\allowdisplaybreaks
\theoremstyle{remark}

\newtheorem{remark}[lemma]{Remark}

\numberwithin{equation}{section}

\title{Triple quantum Schubert calculus}
\author{Matthew J. Samuel}

\DeclareMathOperator{\sch}{\mathfrak{S}}
\DeclareMathOperator{\qsch}{\widetilde{\mathfrak{S}}}

\newcommand{\xrightsquigarrow}[2][]{%
	\mathrel{\mathop{\rightsquigarrow}\limits^{\scriptstyle #2}%
	\ifx\hfuzz#1\hfuzz\else_{\scriptstyle #1}\fi}%
}

\newcommand{\tom}[1]{\ensuremath{\xrightsquigarrow{#1}}}
\newcommand{\fx}{\mathrm{fix}}

\newcommand{\code}{\mathfrak{c}}

\newcommand{\wte}{\mathbf{E}}
\newcommand{\qelem}{\widetilde{E}}
\newcommand{\pcomp}[1]{\widehat{#1}}

\begin{document}

\begin{abstract}
We prove a triple quantum Pieri formula: the product of a quantum double Schubert polynomial $\qsch_u(x;y)$ by a quantum factorial elementary polynomial $\qelem_{p,k}(x;z)$ in a second set of variables $z$ is a sum over chains in the quantum Bruhat graph, indexed by a quantum analogue of Sottile's relation, whose coefficients are monomials in the quantum parameters times polynomials in the differences $y_i-z_j$ with nonnegative integer coefficients, and which is Graham-positive at $z=y$. At $y=z=0$ it is the quantum Pieri rule of Ciocan-Fontanine and Postnikov; at $z=y$ it is an equivariant quantum Pieri rule for the complete flag variety, indexed by chains in the quantum Bruhat graph rather than, as in the rule of Huang and Li, by the data of a reduction to classical cohomology, and proved directly. A doubled Pieri formula, positive in the same senses, treats multiplication by the quantization of a product of two elementary polynomials in the same variables, which is not the product of the quantizations: it is two iterated Pieri steps with the incompatible pairs of chains omitted. More generally, the triple quantum Schubert structure constants $c_{u,v}^{w,d}(y;z)$, the coefficients of $\qsch_u(x;y)\qsch_v(x;z)$ in the basis $\{\qsch_w(x;y)\}$, are polynomials in the $y_i-z_j$ with nonnegative integer coefficients; the proof of this appears as an appendix. We extend the Peterson--Woodward comparison formula and the quantum-to-classical reduction of Leung--Li and Huang--Li to the triple constants, and we describe and prove correct the algorithms of \texttt{schubmult}, which compute products of arbitrary quantum double Schubert polynomials well beyond the reach of previous methods.
\end{abstract}
\maketitle



\section{Introduction}

Let $\sch_w(x;y)$, $w\in S_\infty$, denote the double Schubert polynomials of Lascoux and Sch\"utzenberger, which represent the torus-equivariant Schubert classes of the complete flag variety, and let $\qsch_w(x;y)$ denote the quantum double Schubert polynomials of Kirillov and Maeno \cite{kirillov2000quantum}, which represent the Schubert classes in torus-equivariant quantum cohomology \cite{lam2014quantum}. The structure constants of the equivariant quantum cohomology ring are the coefficients in
$$\qsch_u(x;y)\qsch_v(x;y)=\sum_{w,d}q^d\,c_{u,v}^{w,d}(y;y)\,\qsch_w(x;y),$$
and Mihalcea \cite{mihalcea2006positivity} proved that they are polynomials in the negative simple roots $y_{i+1}-y_i$ with nonnegative integer coefficients.

In this paper we study the coefficients obtained when the two factors are polynomials in different sets of variables, $y$ and $z$:
$$\qsch_u(x;y)\qsch_v(x;z)=\sum_{w,d}q^d\,c_{u,v}^{w,d}(y;z)\,\qsch_w(x;y).$$
We call the $c_{u,v}^{w,d}(y;z)$ the \emph{triple quantum Schubert structure constants}. At $z=y$ they recover the equivariant quantum structure constants, and at $q=0$ they recover the triple Schubert structure constants $c_{u,v}^w(y;z)$ of \cite{samuelmolev}, which were recently shown by Gao and Xiong \cite{gao2025grahampositivitytripleschubert} to be polynomials in the differences $y_i-z_j$ with nonnegative integer coefficients. The triple constants carry strictly more information than the equivariant ones: the equivariant quantum constant $c_{u,v}^{w,d}(y;y)$ is a specialization of $c_{u,v}^{w,d}(y;z)$, while conversely $c_{u,v}^{w,d}(y;z)$ is recovered from the constants $c_{u,v_2}^{w,d}(y;y)$, $v_2$ ranging over the right factors of the reduced factorizations $v=v_1\cdot v_2$, by the quantum Cauchy identity (Lemma \ref{lemma:quantumtoordinary}).

The main results of this paper are the following; the positivity theorem is due to the authors of Appendix \ref{appendix:positivity}, the rest are ours.
\begin{enumerate}
\item \emph{A triple quantum Pieri formula} (Theorem \ref{theorem:pieri}; general degree in Theorem \ref{theorem:pierigeneral}). Let $\qelem_{k,k}(x;z)=\qsch_{s_1s_2\cdots s_k}(x;z)$ be the quantum factorial elementary polynomial. We prove
$$\qsch_u(x;y)\qelem_{k,k}(x;z)=\sum_{u\tom{k} w}q_k(u,w)\prod_{a\in\fx_k(u,w)}(y_a-z_1)\,\qsch_w(x;y),$$
where $\tom{k}$ is a quantum analogue of the relation introduced by Sottile \cite{sottile} for the classical Pieri rule, $q_k(u,w)$ is a monomial in the $q_i$ recording the length-decreasing steps of a chain from $u$ to $w$ in the quantum Bruhat graph, and $\fx_k(u,w)$ is the set of values fixed by the chain in the first $k$ positions. The general-degree version expands $\qsch_u(x;y)\qelem_{p,k}(x;z)$ for every $p\leq k$ over the same chains, the coefficient of $\qsch_w(x;y)$ being $q_k(u,w)$ times a factorial elementary polynomial in the variables $y_a$, $a\in\fx_k(u,w)$, and $z$. In either form the coefficient of each $q$-monomial is a polynomial in the differences $y_i-z_j$ with nonnegative integer coefficients. This is the main result of the paper.

The word triple is essential. Mixed-variable formulas of this kind go back to the Littlewood--Richardson rule of Molev and Sagan for factorial Schur functions \cite{molev1999littlewood}, and the classical triple Pieri formula \cite[Theorem~7.1]{samuelmolev}, to which the present one reduces at $q=0$, has already shown what such a formula is good for: it is the engine of the multiplication algorithm of \texttt{schubmult} for double Schubert polynomials. The reason is that a formula with $z\neq y$ composes: an arbitrary $\qsch_v(x;z)$ expands into products of factorial elementary polynomials in the $z$-variables alone (Theorems \ref{theorem:quandouble} and \ref{theorem:quandoublemed}), and the triple Pieri formula then multiplies each factor against $\qsch_u(x;y)$ without any further change of variables. An equivariant formula, with $z=y$, cannot be iterated in this way. Conversely the triple formula determines the equivariant one by specialization, while recovering it from equivariant data requires the quantum Cauchy identity and a sum over factorizations (Lemma \ref{lemma:quantumtoordinary}).

The specializations locate the formula among known results. At $q=0$ it is the triple Pieri formula of \cite{samuelmolev}, and at $q=0$, $y=z=0$ it is Sottile's rule. At $y=z=0$ it is the quantum Pieri formula for complete flag manifolds, due to Ciocan-Fontanine \cite{ciocan1999partial} and reformulated by Postnikov \cite{postnikov1999quantum} in terms of what is now called the quantum Bruhat graph; the relation $\tom{k}$ is the factorial version of Postnikov's. At $p=1$ and $z=y$ it is the equivariant quantum Chevalley (Monk) formula of Mihalcea \cite{mihalcea2007equivariant}. At $z=y$ and general $p$ it is an equivariant quantum Pieri rule for the complete flag variety. Such a rule was first obtained by Huang and Li \cite[Theorem~3.10]{huang2015equivariant}, for all partial flag varieties of type A, by applying the quantum-to-classical reduction (a principle going back to Buch, Kresch and Tamvakis \cite{buch2003gwgrass}) repeatedly until no quantum parameters remain and then invoking Robinson's equivariant Pieri rule \cite{robinsonpieri}; its terms are indexed by the data of Ciocan-Fontanine's rule, namely a unimodal degree $d$, the permutations $\tau_d$ and $\phi_d$ it determines, and classical chains from $u\tau_d$ to $w\phi_d$, and its coefficients are restrictions of special Schubert classes of Grassmannians to torus-fixed points. Our formula differs from theirs in three respects. It has independent sets of variables $y$ and $z$, which is what makes it composable. Its terms are indexed intrinsically, by chains in the quantum Bruhat graph starting at $u$, with the quantum degree $q_k(u,w)$ and the weight $\fx_k(u,w)$ read off from the chain. This is closer in spirit to Postnikov's reformulation of Ciocan-Fontanine's rule than to the reduction-based formula of Huang and Li. And its proof is a direct analysis of the quantum Bruhat graph rather than a reduction.

It is worth separating the kinds of positivity in play. Positivity of the coefficients follows from Mihalcea's theorem at $z=y$ and from Theorem \ref{theorem:positivity} in general, without any formula; triple positivity does not by itself imply Graham positivity at $z=y$, but the Pieri coefficients are specializations of factorial elementary polynomials and are manifestly Graham-positive there (Corollary \ref{corollary:pierigraham}), exactly as in \cite{samuelmolev}. The quantum-to-classical reduction of Leung--Li \cite{leung2012classical} and Huang--Li \cite{huang2015equivariant} gives more: whenever one factor is Grassmannian, each structure constant equals an explicitly determined classical one, a statement we show holds verbatim for the triple constants (Proposition \ref{prop:grassmannianreduction}), and combining this with a positive classical rule gives a positive expression for each coefficient separately, which is how Huang and Li obtain their formula. A rule indexed by chains from $u$ is a different kind of statement: it describes the support of the product, the degrees $d$ and the permutations $w$ that occur, together with all coefficients, without reference to a target $(w,d)$, and it is this that the multiplication algorithm of Section \ref{section:algorithms} iterates. The proof occupies Section \ref{section:proofpieri}; it proceeds by induction on $k$ via the single-variable Monk formula, and rests on an intrinsic description of the chains (every nontrivial cycle of $u^{-1}w$ has exactly one element exceeding $k$, and a single length identity) which makes the term-by-term cancellations tractable.
\item \emph{A doubled quantum Pieri formula} (Theorem \ref{theorem:doubledpieri}). The quantization map of Fomin, Gelfand and Postnikov is not multiplicative on products of elementary polynomials in the same number of variables, so the quantum Schubert polynomial of a $k\times2$ rectangle is not the product of two quantum factorial elementary polynomials. We show that the product $\qsch_u(x;y)\qsch_{d_2[k]}(x;z)$, where $d_2[k]$ is the dominant permutation whose diagram is the $k\times2$ rectangle, is nevertheless given by two iterated Pieri steps, provided one omits the pairs of chains that are not \emph{compatible} in a precise cycle-theoretic sense (Definition \ref{def:compatible}). The doubled rule is equivalent to a termwise identity between pairs of chains (Theorem \ref{theorem:termwise}) which telescopes along explicit blocks of pairs and is then established by an injection built from a single local lemma, a generalization of the up/down pairing used in the proof of the Pieri formula (Theorem \ref{theorem:junction}). By quantum Giambelli the same product is an alternating sum of iterated Pieri steps; the doubled rule is the cancellation-free form of that sum, and it is the ``paired step'' of the multiplication algorithm.
\item \emph{Positivity} (Theorem \ref{theorem:positivity}). The constants $c_{u,v}^{w,d}(y;z)$ lie in $\mathbb Z_{\geq0}[y_i-z_j]$. This theorem and its proof are due to the authors of Appendix \ref{appendix:positivity}; the proof combines the geometric argument of Mihalcea with the Graham-positivity theorem of Gao and Xiong for $B^-$-stable effective cycles. We include it because it settles what would otherwise be a conjecture about the constants studied here, of which the Pieri formulas above are explicit combinatorial instances.
\item \emph{Parabolic and quantum-to-classical reductions} (Proposition \ref{prop:parabolic}, Theorem \ref{theorem:quantclassz}). The Peterson--Woodward comparison formula \cite{woodward2005d}, in the form proved by Huang and Li \cite{huang2015equivariant} for equivariant quantum cohomology, extends verbatim to the triple constants. Consequently the triple structure constants of the parabolic quantum double Schubert polynomials of Lam and Shimozono are equal to triple constants of the complete flag variety, and in particular they are also Graham-positive (Corollary \ref{corollary:parabolicpositivity}). The quantum-to-classical reduction of Leung--Li and Huang--Li extends to the triple constants as well. For Grassmannians, considerably more is known at $z=y$: Mihalcea \cite{mihalcea2006equivariant} gave an equivariant quantum Pieri rule and a recursive characterization of the equivariant quantum structure constants, and Buch's mutation algorithm for puzzles \cite{buch2015mutations}, combined with the equivariant ``quantum equals classical'' theorem of Buch and Mihalcea \cite{buch2011qkgrass}, provides a Graham-positive combinatorial Littlewood--Richardson rule for $QH_T^*(\mathrm{Gr}(k,n))$. No such rule is known to us for the triple constants with $z\neq y$, even in the Grassmannian case.
\item \emph{Algorithms} (Section \ref{section:algorithms}). Computing products of Schubert classes in $QH^*(Fl_n)$, even non-equivariantly, is difficult in practice: the general methods available, namely reduction of quantum Giambelli determinants modulo the quantum ideal and iteration of the quantum Monk or Chevalley formula with a recursive reduction, become infeasible well before $n$ reaches the range in which interesting phenomena appear. We describe the algorithms implemented in the package \texttt{schubmult} \cite{schubmult}, which compute products of arbitrary (quantum) double Schubert polynomials, with independent parameters $y$ and $z$, by combining the Pieri formulas on one factor with a signed expansion of the other factor into products of factorial elementary polynomials. In our experience this is by a wide margin the fastest method available, and it is the only one we know of that handles the triple constants directly. We prove the expansion, in its classical form (Theorem \ref{theorem:orddouble}) and in two quantum forms, over a strict dominant shape (Theorem \ref{theorem:quandouble}) and over the ``medium'' shape actually used by the implementation, in which a pair of equal columns is allowed when it is preceded by a drop of at least two and contributes a single doubled factor (Theorem \ref{theorem:quandoublemed}); the key step is a formula for the skew divided differences, in the sense of \cite{samuelleibniz,samuelmolev}, of a product $\prod_{j\leq k}(x_j-z_i)$ with respect to the $z$-variables (Theorem \ref{theorem:skewtop}), in which the terms of the iterated Monk rule that revisit a position cancel by a rotation argument (Lemma \ref{lemma:rotation}). The doubled Pieri rule is exactly what is needed to process the equal columns of the medium shape.
\end{enumerate}

The paper is organized as follows. Section \ref{section:definitions} defines the triple constants, states the positivity theorem and treats the parabolic case. Section \ref{section:pieri} states the Pieri formula, whose proof is given in Section \ref{section:proofpieri}. Section \ref{section:algorithms} states the expansion formulas for (quantum) double Schubert polynomials and describes the multiplication algorithms built on them; the expansion formulas are proved in Section \ref{section:proofvside}. Section \ref{section:quantclass} treats the quantum-to-classical reduction, Section \ref{section:additional} the general-degree and doubled Pieri formulas, and Section \ref{section:pw} proves the parabolic reduction via the Peterson--Woodward comparison formula. Appendix \ref{appendix:positivity} contains the proof of the positivity theorem.

\textbf{Acknowledgements.} We thank Rui Xiong and Changzheng Li for helpful discussions.

\section{Triple quantum Schubert structure constants} \label{section:definitions}

\subsection{Quantum double Schubert polynomials and the triple constants}
We fix $n\geq 1$ and write
$$x=(x_1,\ldots,x_n),\qquad y=(y_1,\ldots,y_n),\qquad q=(q_1,\ldots,q_{n-1}).$$
For $1\leq a<b\leq n$, set
$$q_{ab}:=q_aq_{a+1}\cdots q_{b-1},$$
and use the cohomological grading $\deg q_i=2$ (equivalently $\deg q_{ab}=2(b-a)$).

Let $s_i$ act on the $y$-variables by swapping $y_i,y_{i+1}$, and define the divided difference operators
$$\partial_i^{(y)}f:=\frac{f-s_i f}{y_i-y_{i+1}}\qquad(1\leq i\leq n-1).$$
For $w\in S_n$, if $w=s_{i_1}\cdots s_{i_r}$ is any reduced word, set
$$\partial_w^{(y)}:=\partial_{i_1}^{(y)}\cdots \partial_{i_r}^{(y)}.$$

Following Kirillov--Maeno \cite{kirillov2000quantum}, define the quantum elementary polynomials $e_j^{(q)}(x_1,\ldots,x_k)$ by
$$\Delta_k(t\mid x_1,\ldots,x_k):=\sum_{j=0}^k t^{k-j}e_j^{(q)}(x_1,\ldots,x_k),$$
where
$$\Delta_k(t\mid x_1,\ldots,x_k)=\det\!\begin{pmatrix}
x_1+t & q_1 & 0 & \cdots & 0\\
-1 & x_2+t & q_2 & \ddots & \vdots\\
0 & -1 & x_3+t & \ddots & 0\\
\vdots & \ddots & \ddots & \ddots & q_{k-1}\\
0 & \cdots & 0 & -1 & x_k+t
\end{pmatrix}.
$$

Let $w_0\in S_n$ be the longest permutation. The top quantum double Schubert polynomial is
$$\qsch_{w_0}(x;y):=\prod_{i=1}^{n-1}\Delta_i(-y_{n-i}\mid x_1,\ldots,x_i),$$
and for $w\in S_n$ we define
$$\qsch_w(x;y):=(-1)^{\ell(w_0) - \ell(w)}\partial_{ww_0}^{(y)}\qsch_{w_0}(x;y).$$
Its specialization
$$\qsch_w(x):=\qsch_w(x;0)$$
is the quantum Schubert polynomial. In particular, at $q=0$ these recover the ordinary double and single Schubert polynomials.

\subsubsection*{The stable setting}
The polynomials $\qsch_w(x;y)$ are stable \cite{kirillov2000quantum}: for $w\in S_n$, $\qsch_w(x;y)$ involves only $x_1,\ldots,x_{n-1}$, $y_1,\ldots,y_{n-1}$ and $q_1,\ldots,q_{n-2}$, and it does not change when $w$ is regarded as an element of $S_{n+1}$ under the standard inclusion $S_n\subset S_{n+1}$. (This is also immediate from the description $\qsch_w(x;y)=\mathcal Q(\sch_w(x;y))$ in the proof of Lemma \ref{lemma:quantumtoordinary} below, since the quantization map $\mathcal Q$ does not depend on $n$.) From now on we therefore let the number of variables be infinite,
$$x=(x_1,x_2,\ldots),\qquad y=(y_1,y_2,\ldots),\qquad z=(z_1,z_2,\ldots),\qquad q=(q_1,q_2,\ldots),$$
index the polynomials by $S_\infty=\bigcup_nS_n$, and work in the polynomial ring $\mathbb Z[q][x,y,z]$ in these (countably many) variables, where $\mathbb Z[q]=\mathbb Z[q_1,q_2,\ldots]$. A \emph{degree} is an element $d=(d_1,d_2,\ldots)$ of $\bigoplus_{i\geq1}\mathbb Z_{\geq0}$, i.e.\ a weak composition with finitely many nonzero parts, and $q^d=\prod_iq_i^{d_i}$; in Section \ref{section:definitions:parabolic} degrees will be regarded as effective elements of the coroot lattice. As a module over $\mathbb{Z}[q][y]$, the quantum double Schubert polynomials $\qsch_w(x;y)$, $w\in S_\infty$, form a basis of $\mathbb{Z}[q][x,y]$: the quantization map $\mathcal Q$ of Lemma \ref{lemma:quantumtoordinary} is a $\mathbb Z[q][y]$-linear bijection of $\mathbb Z[q][x,y]$ onto itself, since $\mathcal Q-\mathrm{id}$ lowers the total $x$-degree by at least two and is therefore locally nilpotent, so that $\sum_{r\geq0}(-1)^r(\mathcal Q-\mathrm{id})^r$ is a well-defined inverse on every polynomial. It carries the basis $\{\sch_w(x;y)\}$ to $\{\qsch_w(x;y)\}$. We may therefore write
$$\qsch_u(x;y)\qsch_v(x;z)=\sum_{w,d}q^dc_{u,v}^{w,d}(y;z)\qsch_w(x;y),$$
the sum being over $w\in S_\infty$ and all degrees $d$, with only finitely many nonzero terms. The coefficients $c_{u,v}^{w,d}(y;z)$ are the \emph{triple quantum Schubert structure constants}, which are polynomials in $y$ and $z$ with integer coefficients.

The finite-dimensional meaning of this expansion is the following. Given $u,v\in S_n$, choose $N$ so large that every $w$ occurring in the expansion lies in $S_N$ and every occurring $d$ has $d_i=0$ for $i\geq N$; this holds for all sufficiently large $N$, since there are finitely many terms. Then the expansion is an identity in $\mathbb Z[q_1,\ldots,q_{N-1}][x_1,\ldots,x_N,y,z]$, and reducing it modulo the ideal of Kim's presentation of $QH_T^*(\mathrm{Fl}_N)$ \cite{kim1996equivariant} (the ideal $J_{qy}^N$ recalled in Appendix \ref{appendix:positivity}) identifies it, by \cite{lam2014quantum}, with the expansion of the product of the two corresponding classes in the Schubert basis of $QH_T^*(\mathrm{Fl}_N)$, with $z=y$. Thus $c_{u,v}^{w,d}(y;y)$ is the equivariant quantum Littlewood--Richardson coefficient of $\mathrm{Fl}_N$ for every such $N$, and the degree vector $d$ is the degree $\sum_id_i\alpha_i^\vee$ of a rational curve in $\mathrm{Fl}_N$. Note that $w\notin S_n$ and $d_i\neq0$ for $i\geq n-1$ do occur: for instance $\qsch_{21}(x;y)^2=(y_2-y_1)\qsch_{21}(x;y)+\qsch_{231}(x;y)+q_1$, the last term coming from the quantum Monk formula in $\mathrm{Fl}_N$ for any $N\geq3$.

The product respects the grading where for each $i$, $\deg x_i=\deg y_i=\deg z_i=1$ and $\deg q_i=2$, in which $\qsch_w$ is homogeneous of degree $\ell(w)$. In particular, if $c_{u,v}^{w,d}(y;z)\neq 0$, then
$$\deg c_{u,v}^{w,d}(y;z)=\ell(u)+\ell(v)-\ell(w)-2\sum_i d_i;$$
for instance the term $q_1$ in $\qsch_{21}^2$ above has a coefficient of degree $0=1+1-0-2$. Unlike the usual triple Schubert structure constants, they do not exhibit triangularity in Bruhat order, and indeed we may have $\ell(w) < \ell(u) + \ell(v)$. There is an ordering in which the quantum double Schubert polynomials are triangular, but we will not need this in full generality.

Throughout, a factorization $w=v_1\cdot v_2$ written with a dot is understood to be \emph{reduced}, i.e.\ $\ell(w)=\ell(v_1)+\ell(v_2)$. The following identity will be used repeatedly to pass between the triple constants and the equivariant quantum structure constants $c_{u,v}^{w,d}(y;y)$. Setting $y=0$ in it recovers the expansion $\qsch_w(x;z)=\sum_{v_1\cdot v_2=w}\sch_{v_1}(0;z)\qsch_{v_2}(x)$ of Kirillov--Maeno \cite{kirillov2000quantum} in terms of quantum Schubert polynomials.

\begin{lemma}[Quantum Cauchy identity] \label{lemma:quantumtoordinary}
For all $w\in S_\infty$,
$$\qsch_w(x;z)=\sum_{v_1\cdot v_2=w}\sch_{v_1}(y;z)\,\qsch_{v_2}(x;y).$$
Consequently, for all $u,v,w$ and $d$,
$$c_{u,v}^{w,d}(y;z)=\sum_{v_1\cdot v_2=v}c_{u,v_2}^{w,d}(y;y)\,\sch_{v_1}(y;z).$$
\end{lemma}
\begin{proof}
Let $\mathcal Q$ denote the quantization map of Fomin, Gelfand and Postnikov \cite{fomin1997quantum}, used throughout the paper: the $\mathbb Z[q]$-linear map on polynomials in $x$ which sends each standard elementary monomial $\prod_j e_{i_j}(x_1,\ldots,x_j)$ to $\prod_j e^{(q)}_{i_j}(x_1,\ldots,x_j)$, extended linearly over any further coefficient variables. Since the divided differences in the $y$-variables commute with $\mathcal Q$, the definition of $\qsch_w(x;y)$ by $y$-divided differences gives $\qsch_w(x;y)=\mathcal Q(\sch_w(x;y))$ \cite[(10)]{lam2014quantum}. Now apply $\mathcal Q$, which is $\mathbb Z[y,z]$-linear, to the classical Cauchy-type expansion $\sch_w(x;z)=\sum_{v_1\cdot v_2=w}\sch_{v_1}(y;z)\sch_{v_2}(x;y)$ \cite[(6.3)]{notes}: the left side becomes $\qsch_w(x;z)$ and the right side becomes $\sum_{v_1\cdot v_2=w}\sch_{v_1}(y;z)\qsch_{v_2}(x;y)$. For the second assertion, multiply the first identity (with $w=v$) by $\qsch_u(x;y)$ and expand each $\qsch_u(x;y)\qsch_{v_2}(x;y)$ in the basis $\{\qsch_w(x;y)\}$ with coefficients $q^dc_{u,v_2}^{w,d}(y;y)$.
\end{proof}

\begin{remark}[$\tom{k}$ and the quantum $k$-Bruhat order] \label{remark:qkbruhat}
The \emph{quantum $k$-Bruhat order} was introduced by Benedetti, Bergeron, Colmenarejo, Saliola, and Sottile \cite{benedetti2025quantum} and studied further by Colmenarejo and Mayers \cite{colmenarejo2026quantum}. That order $\preceq_k^{\mathbf q}$ on $S_\infty[\mathbf q]$ is generated by cover relations
$$u'\;\lessdot_k^{\mathbf q}\;\mathbf q^{\gamma}u't_{ab},\qquad a\leq k<b,$$
each of exactly two kinds: a \emph{classical} cover, with $\gamma=0$ and $\ell(u't_{ab})=\ell(u')+1$, or a \emph{quantum} cover, with $\ell(u')+1=\ell(\mathbf q^{\gamma}u't_{ab})=\ell(u't_{ab})+2(b-a)$, that is $\ell(u't_{ab})=\ell(u')-2(b-a)+1$ and $\gamma=e_a+\cdots+e_{b-1}$. 
These are exactly the two alternatives of condition (4) of Definition \ref{definition:pieri}. Consequently every sequence $t_{a_1b_1},\ldots,t_{a_pb_p}$ witnessing $u\tom{k}w$ is a saturated chain from $u$ to $\mathbf q^{\alpha}w$ in the quantum $k$-Bruhat order, where $\mathbf q^{\alpha}$ is the monomial recorded by $q_k(u,w)$; in particular $u\tom{k}w$ implies $u\preceq_k^{\mathbf q}\mathbf q^{\alpha}w$.

The relation $\tom{k}$ is \emph{not} the full order $\preceq_k^{\mathbf q}$, however. Definition \ref{definition:pieri} imposes the additional requirements that the lower indices $a_1,\ldots,a_p$ be pairwise distinct, so that $\tom{k}$ admits only those saturated chains of a special form. Thus $\tom{k}$ is strictly larger than the covering relation $\lessdot_k^{\mathbf q}$ but strictly smaller than $\preceq_k^{\mathbf q}$: a general chain in the quantum $k$-Bruhat order may reuse a lower index $a_i$, and such chains do not contribute to the Pieri expansion. This is the quantum counterpart of the classical situation, where the Pieri relation of \cite{sottile,samuelmolev} is likewise the restriction of the $k$-Bruhat order to chains with distinct lower indices, and it is exactly this restriction that makes the Pieri coefficients well defined.

What we use below is only that each individual step of a $\tom{k}$ chain is a cover, and hence is governed by the following length-free description of the covers \cite[Proposition~3.1]{benedetti2025quantum} (see also \cite[Proposition~4.2]{colmenarejo2026quantum}): with $i\leq k<j$, the step $u'\lessdot_k^{\mathbf q}\mathbf q^{\gamma}u't_{ij}$ is a classical cover exactly when $u'(i)<u'(j)$ and no $l$ with $i<l<j$ satisfies $u'(i)<u'(l)<u'(j)$, and it is a quantum cover exactly when $u'(i)>u'(j)$ and every $l$ with $i<l<j$ satisfies $u'(j)<u'(l)<u'(i)$. This characterization is an immediate consequence of the standard formula for the length of $u't_{ij}$ (Lemma \ref{lemma:lengthchange} below), and it is the form in which the covers are used throughout Section \ref{section:proofpieri}.
\end{remark}

\subsection{Positivity}

For $d=(0,\ldots,0)$, this is the ordinary triple Schubert calculus structure constant, and it was proved in \cite{gao2025grahampositivitytripleschubert} that these have the positivity property that they are polynomials in the differences $y_i-z_j$ with nonnegative integer coefficients. We will refer to this as \emph{triple Schubert positivity}. The following theorem shows that it persists in the quantum setting.
\begin{theorem}[Triple quantum Schubert positivity] \label{theorem:positivity}
For all $u,v,w$ and $d$, the structure constants $c_{u,v}^{w,d}(y;z)$ are polynomials in the differences $y_i-z_j$ with nonnegative integer coefficients.
\end{theorem}
The proof, due to Tongren Xiao and Yihua Yue, is given in Appendix \ref{appendix:positivity}.

\subsection{The parabolic case} \label{section:definitions:parabolic}

Unlike ordinary cohomology or equivariant cohomology, quantum cohomology and equivariant quantum cohomology do not demonstrate simple functoriality, in the sense that a projection of the complete flag variety $G/B$ to a partial flag variety $G/P$ does not induce a homomorphism of rings $QH_T^*(G/P)\to QH_T^*(G/B)$, and indeed the quantum cohomology of $G/P$ is not isomorphic to the subring of the quantum cohomology of $G/B$ generated by the Schubert classes indexed by $W^P$. This is reflected in the fact that there are parabolic quantum double Schubert polynomials defined in a different manner that represent Schubert classes in $QH_T^*(G/P)$ \cite{lam2014quantum} and are not in general equal to the quantum double Schubert polynomials above. (The non-equivariant ring $QH^*(SL_n/P)$ itself is well understood: its presentation and its quantum Pieri and Giambelli formulas were proved by Ciocan-Fontanine \cite{ciocan1999partial} and, by elementary geometric arguments, by Buch \cite{buch2005partial}.)

\begin{definition}[Parabolic quantum double Schubert polynomials] \label{def:parabolic}
	Following Lam and Shimozono \cite{lam2014quantum}, we define parabolic quantum double Schubert polynomials as follows. Fix an integer $n>0$. A parabolic subgroup $W_P$ of $S_{n}$, which has Coxeter generators $s_1,\ldots,s_{n-1}$, is generated by a subset $P\subseteq \{s_1,\ldots,s_{n-1}\}$. In terms of the subgroup's properties, it is most natural to identify $P$ by the reflections that it \emph{excludes}. Suppose $P$ is generated by $S\setminus \{s_{i_1},\ldots,s_{i_k}\}$ with $1\leq i_1<\cdots<i_k\leq n-1$, and set $i_0=0$ and $i_{k+1}=n$. We define the \emph{block sizes} $n_1,\ldots,n_{k+1}>0$, with
	$$n_1+\cdots+n_{k+1}=n,$$
	by $n_p = i_{p} - i_{p-1}$ for $1\leq p\leq k+1$, and put
	$$N_j = n_1+\cdots+n_j=i_j\qquad(0\leq j\leq k+1),$$
	so that $N_0=0$ and $N_{k+1}=n$. Thus $W_P=S_{n_1}\times\cdots\times S_{n_{k+1}}$ and $SL_n/P$ is the variety of partial flags with subspaces of dimensions $N_1<\cdots<N_k$. The extreme cases are $k=0$, where $W_P=S_n$, $W^P=\{1\}$, $SL_n/P$ is a point and the only polynomial defined below is $\qsch^P_1=1$, and $P=\emptyset$ ($k=n-1$, all $n_p=1$), which recovers the quantum double Schubert polynomials of the complete flag variety.
	
	For such a setup, let $D^P$ be the $n\times n$ matrix with $x_1,\ldots,x_n$ on the diagonal, $-1$ on the superdiagonal, the entry $(-1)^{n_{j+1}+1}q_j$ in position $(N_{j+1},N_{j-1}+1)$ for $1\leq j\leq k$, and all other entries $0$, and let $D^P_j$ be its upper left $N_j\times N_j$ submatrix. Notice that there are precisely $k$ $q$ variables, $q_j$ corresponding to the excluded reflection $s_{i_j}$, and we set
	$$\deg q_j = n_j + n_{j+1}.$$  
	
	Define 
	$$\Delta^P_j(t;x)=\det(D^P_j+t\,\mathrm{Id})\qquad(1\leq j\leq k),$$
	so that for $P=\emptyset$ this is the polynomial $\Delta_j(t\mid x_1,\ldots,x_j)$ above (the matrices differ by a transposition). Then we define the parabolic quantum double Schubert polynomials $\qsch_w^P(x;y)$ for $w\in W^P$ by
	$$\qsch^P_w(x;y)=(-1)^{\ell(w(w_0^P)^{-1})}\partial_{w(w_0^P)^{-1}}^{(y)}\prod_{j=1}^k\ \prod_{i=n-N_{j+1}+1}^{n-N_j} \Delta^P_j(-y_i;x),$$
	where $w_0^P$ is the longest element of $W^P$. For $P=\emptyset$ the product is $\prod_{j=1}^{n-1}\Delta_j(-y_{n-j}\mid x_1,\ldots,x_j)=\qsch_{w_0}(x;y)$.
	
	Taken as defined, when we vary $n$ these polynomials are not stable under inclusion and do not form a basis of $\mathbb{Z}[q][x_1,\ldots,x_n;y]$ as a module over $\mathbb{Z}[q][y]$. To obtain stability and a basis, instead we fix $n$ and choose $N>n$ as large as we would like, maintaining $N_{k+1}=n$ and setting $n_{k+2} = N-n$. The added quantum parameter $q_{k+1}$ has degree $n_{k+1}+N-n$, which tends to infinity; hence it cannot occur in any polynomial or product of fixed degree once $N$ is sufficiently large. The remaining parameters are the original $q_1,\ldots,q_k$. From a cohomological perspective, this introduces ``extraneous'' polynomials that are indexed by elements not in $W^P$. The full indexing set is $W^{\pcomp{P}}\subset S_\infty$, where
	$$\pcomp{P} = P \cup \{s_{n+1},s_{n+2},\ldots\}$$
		The notation $W^{\pcomp P}$ specifies the stable indexing set only. Any use of longest elements below is made in a finite symmetric group, with the enlarged parabolic generated by $\pcomp P\cap\{s_1,\ldots,s_{N-1}\}$.
\end{definition}

\begin{example}
As equivariant cohomology rings of finite dimensional projective varieties over $\mathbb{C}$ are finitely generated modules over the coefficient ring, it is fully expected that products of double Schubert polynomials may include terms indexed by permutations that do not lie within the indexing set of the basis of Schubert classes. For example, consider the following non-quantum double Schubert polynomial product, with $P=\{s_1,s_3\}$ and basis elements in $S_4$:
\begin{align*}
\sch_{2413}(x;y)\sch_{3412}(x;y)
&= - (\alpha_1+\alpha_2)(\alpha_1+\alpha_2+\alpha_3)(\alpha_2+\alpha_3)\,\sch_{3412}(x;y) \\
&\quad + (\alpha_1+\alpha_2)(\alpha_1+2\alpha_2+2\alpha_3+\alpha_4)\,\sch_{35124}(x;y) \\
&\quad - (\alpha_1+\alpha_2)\,\sch_{361245}(x;y) \\
&\quad - (\alpha_1+2\alpha_2+2\alpha_3+\alpha_4)\,\sch_{45123}(x;y) \\
&\quad + \sch_{461235}(x;y).
\end{align*}
Not only does the product $\qsch_{2413}^{P}(x;y)\qsch_{3412}^{P}(x;y)$ include extra terms indexed by permutations in $W^P$, there are also terms with nonzero coefficients that are in $W^{\pcomp{P}}\setminus W^P$:
\begin{align*}
q_1(\alpha_1\alpha_2+\alpha_1\alpha_3+\alpha_2\alpha_3+\alpha_2^2)\,\qsch_{132}^{P}(x;y)
&+ q_1(\alpha_1+\alpha_2)\,\qsch_{231}^{P}(x;y) \\
&+ q_1(\alpha_2+\alpha_3)\,\qsch_{1423}^{P}(x;y) \\
&+ \color{red}{q_1(\alpha_1+\alpha_2)\,\qsch_{13254}^{\pcomp{P}}(x;y)} \\
&+ \color{red}{q_1\,\qsch_{14253}^{\pcomp{P}}(x;y)} \\
&+ q_1\,\qsch_{15234}^{P}(x;y) \\
&+ q_1\,\qsch_{2413}^{P}(x;y)
\end{align*}
\end{example}

\begin{definition}[Coroot lattices, $\langle d\rangle_i$, Peterson--Woodward lifts] \label{def:pwlift}
	Let $\alpha_1,\alpha_2,\ldots$ be the simple roots of $S_\infty$ and $\alpha_i^\vee$ the simple coroots, and let
	$$Q^\vee=\bigoplus_{i\geq1}\mathbb Z\alpha_i^\vee$$
	be the coroot lattice. We write its elements as integer vectors $d=(d_1,d_2,\ldots)$ with finitely many nonzero entries, $d=\sum_id_i\alpha_i^\vee$; the degrees of Section \ref{section:definitions} are the \emph{effective} elements, those with all $d_i\geq0$, and $q^d=\prod_iq_i^{d_i}$. (When working in $S_n$ we take $i\leq n-1$ and $d\in\mathbb Z^{n-1}$.) For $d\in Q^\vee$ and $i\geq1$ put
	$$\langle d\rangle_i = \langle\alpha_i,d\rangle=2d_i - d_{i-1} - d_{i + 1},$$
	with the convention $d_0=0$ (and $d_i=0$ for $i\geq n$ when working in $S_n$), and for $i<j$
	$$\langle d\rangle_{ij} = \sum_{l=i}^{j-1}\langle d\rangle_l=\langle\alpha_i+\cdots+\alpha_{j-1},d\rangle=\langle e_i-e_j,d\rangle.$$
	We write $\delta_i\in Q^\vee$ for the simple coroot $\alpha_i^\vee$, the vector with a single $1$ in position $i$, and $|d|=\sum_id_i$.

	For a parabolic subgroup $W_P$ generated by $P$, let
	$$Q^\vee_P=\bigoplus_{s_i\in P}\mathbb Z\alpha_i^\vee=\{d\in Q^\vee\mid d_i=0\text{ whenever }s_i\notin P\}$$
	be the coroot lattice of the Levi subgroup of $P$. The degrees of $QH_T^*(G/P)$ are indexed by the quotient $Q^\vee/Q^\vee_P$, which we identify with $\bigoplus_{s_i\notin P}\mathbb Z\alpha_i^\vee$ by sending the class of $d$ to $d_P:=(d_i)_{s_i\notin P}$, the vector of coordinates of $d$ at the excluded reflections; $d_P$ is called effective if these coordinates are nonnegative. In the notation of Definition \ref{def:parabolic}, the coordinate of $d_P$ at the excluded reflection $s_{i_j}$ is the exponent of $q_j$, so that $q^{d_P}=\prod_{j=1}^kq_j^{d_{i_j}}$. (For $P=\emptyset$, $Q^\vee_P=0$ and $d_P=d$.)

	Let $d_P\in Q^\vee/Q^\vee_P$. A \emph{Peterson--Woodward lift} of $d_P$ is a $d\in Q^\vee$ such that $d+ Q_P^\vee=d_P$ and $\langle d\rangle_{ij}\in \{0,-1\}$ for all $i<j$ such that $s_{i},\ldots,s_{j-1}\in P$. In the language of roots, the second condition says that $\langle\alpha,d\rangle\in\{0,-1\}$ for every positive root $\alpha$ of the Levi subgroup of $P$, these being the roots $\alpha_i+\cdots+\alpha_{j-1}$ with $s_i,\ldots,s_{j-1}\in P$. Every $d_P$ has exactly one Peterson--Woodward lift \cite{woodward2005d}; we denote it by $d$ or $\lambda_B$ when $d_P=\lambda_P$.
 \end{definition}

 \begin{method}[Checking the Peterson-Woodward lift property]
 While backing out a Peterson-Woodward lift is difficult, going in the other direction is not. Suppose $d\in Q^\vee$ and $P$ is parabolic. If there is any $i$ with $s_i\in P$ such that $\langle d\rangle_i\notin\{0,-1\}$, then $d$ is not a Peterson-Woodward lift. If this condition does hold, let $I_1,\ldots,I_m$ be the maximal intervals of consecutive indices $i$ with $s_i\in P$. Since $\langle d\rangle_{ij}=\sum_{l=i}^{j-1}\langle d\rangle_l$, the condition $\langle d\rangle_{ij}\in\{0,-1\}$ for all $i<j$ with $s_i,\ldots,s_{j-1}\in P$ holds if and only if each interval $I_r$ contains at most one index $l$ with $\langle d\rangle_l=-1$. If and only if this is the case, $d$ is the Peterson-Woodward lift of a unique $d_P\in Q^\vee/Q^\vee_P$, namely its image $d+Q^\vee_P$, which we denote by $\pi_P(d)$. In that case we let
$$z_P(d)=\{\,i\mid s_i\in P,\ \langle d\rangle_i=0\,\}.$$
When $P\subseteq\{s_1,\ldots,s_{N-1}\}$ is finite, write $w_{z_P(d)}$ for the longest element of the parabolic subgroup generated by $\{s_i\mid i\in z_P(d)\}$, and $w_P$ for the longest element of $W_P\subseteq S_N$.
 \end{method}

 \begin{example}
 Let $n=3$ and $P=\{s_1\}$, so that $W_P=S_2\times S_1$, $SL_3/P=\mathbb P^2$, $Q^\vee_P=\mathbb Z\alpha_1^\vee=\{(d_1,0)\}$, and $Q^\vee/Q^\vee_P\cong\mathbb Z$ via $d\mapsto d_2$. The lifts of $d_P=1$ are the vectors $(d_1,1)$, with $\langle d\rangle_1=2d_1-1$; the only one with $\langle d\rangle_1\in\{0,-1\}$ is $d=(0,1)$, with $\langle d\rangle_1=-1$ and $z_P(d)=\emptyset$, so that $w_Pw_{z_P(d)}=s_1$. For $d_P=2$ the lifts are $(d_1,2)$ with $\langle d\rangle_1=2d_1-2$, and the Peterson--Woodward lift is $d=(1,2)$, with $\langle d\rangle_1=0$, $z_P(d)=\{1\}$ and $w_Pw_{z_P(d)}=1$.
 \end{example}

 This algorithm is polynomial in the length of $d$, which in almost all currently feasible cases is quite small. In contrast, a naive search for $d$ with $d_P$ as the input has an exponential search space, though a more refined algorithm may exist.

 The parabolic triple quantum Schubert structure constants $c_{u,v}^{w,d_P}(y;z)$ are defined by
$$\qsch^P_u(x;y)\qsch^P_v(x;z)=\sum_{w\in W^{\pcomp P},\,d_P}q^{d_P}\,c_{u,v}^{w,d_P}(y;z)\,\qsch^P_w(x;y).$$
We use the same notation for the finite-rank structure constants of $QH_T^*(SL_N/P)$ when the ambient rank is specified. We prove the following finite-rank result in Section \ref{section:pw}, as a consequence of the Peterson--Woodward comparison formula of Huang and Li \cite{huang2015equivariant} and Lemma \ref{lemma:quantumtoordinary}, and then apply it to the stabilized polynomial coefficients.

 \begin{proposition} \label{prop:parabolic}
Fix $N$ and a parabolic subgroup $W_P\subseteq S_N$ generated by $P\subseteq\{s_1,\ldots,s_{N-1}\}$, and let $u,v,w\in W^{P}$. For the structure constants of $QH_T^*(SL_N/P)$ and $QH_T^*(\mathrm{Fl}_N)$, and any $d_P\in Q^\vee/Q_P^\vee$ with Peterson--Woodward lift $d$, we have
$$c_{u,v}^{w,d_P}(y;z) = c_{u,v}^{ww_Pw_{z_P(d)},d}(y;z).$$
\end{proposition}

For the stabilized polynomial coefficients, start with the block sizes of Definition \ref{def:parabolic} and apply Proposition \ref{prop:parabolic} after adjoining a final block of size $N-n$, with $N$ sufficiently large to contain all permutations in the product expansion. At fixed total degree the added parameter has degree too large to occur, so its degree coordinate is zero. The Peterson--Woodward lift for the original blocks, extended by zero in positions $i\geq n$, is also the lift for the enlarged parabolic: its root pairings on the added block are zero. Thus the longest element of that block occurs in both $w_P$ and $w_{z_P(d)}$ and cancels in their product. The twisted output permutation and the lift are consequently independent of further enlargement, and the finite-rank comparison gives the comparison for the stabilized coefficients. No longest element of $W_{\pcomp P}$ is required.

Thus, though we lack a simple functoriality property, we can reduce the computation of parabolic triple quantum Schubert structure constants completely to the complete flag variety case with an efficient algorithm. We obtain a nontrivial corollary from this.

\begin{corollary}[Triple quantum Schubert positivity for the parabolic case] \label{corollary:parabolicpositivity}
Let $W_P\subseteq S_N$ be a parabolic subgroup and let $u,v,w\in W^P$. Then for every $d_P$ the parabolic structure constant $c_{u,v}^{w,d_P}(y;z)$ exhibits triple Schubert positivity. The same holds for the stabilized polynomial coefficients. In particular this holds when $u,v,w$ are $k$-Grassmannian and $P$ is the corresponding maximal parabolic subgroup ($n_1=k$, $n_2=N-k$).
\end{corollary}
\begin{proof}
By Proposition \ref{prop:parabolic}, $c_{u,v}^{w,d_P}(y;z)$ equals a structure constant $c_{u,v}^{w',d}(y;z)$ of the complete flag variety, which lies in $\mathbb Z_{\geq0}[y_i-z_j]$ by Theorem \ref{theorem:positivity}. The stabilized assertion follows by the finite-rank interpretation just explained.
\end{proof}

\begin{remark}[The Grassmannian case] \label{remark:grassmannian}
When $P$ is a maximal parabolic subgroup, $G/P=\mathrm{Gr}(k,n)$, and at $z=y$ the parabolic constants $c_{u,v}^{w,d_P}(y;y)$ are the structure constants of $QH_T^*(\mathrm{Gr}(k,n))$, whose Graham positivity is Mihalcea's theorem \cite{mihalcea2006positivity}. In this case a positive \emph{combinatorial} rule is available. The equivariant ``quantum equals classical'' theorem of Buch and Mihalcea \cite{buch2011qkgrass} (the equivariant form of \cite{buch2003gwgrass}) identifies each $c_{u,v}^{w,d_P}(y;y)$ with an equivariant structure constant of a two-step flag variety $\mathrm{Fl}(k-d,k+d;n)$, and Buch's mutation algorithm for puzzles \cite{buch2015mutations} expresses the latter as a Graham-positive sum over equivariant puzzles; together they give a Littlewood--Richardson rule for $QH_T^*(\mathrm{Gr}(k,n))$ (see also the earlier equivariant quantum Pieri rule and recursion of Mihalcea \cite{mihalcea2006equivariant}, and the non-equivariant puzzle rule for two-step flag varieties of \cite{buch2016puzzle}). Proposition \ref{prop:parabolic} describes the same numbers as structure constants of the complete flag variety, with the degree $d_P$ replaced by its Peterson--Woodward lift and $w$ by $ww_Pw_{z_P(d)}$, which is how the package \cite{schubmult} computes them. For $z\neq y$ the triple Grassmannian constants are determined by the equivariant ones through the parabolic quantum Cauchy identity (proof of Corollary \ref{corollary:pwz}), $c_{u,v}^{w,d_P}(y;z)=\sum_{v_1\cdot v_2=v}c_{u,v_2}^{w,d_P}(y;y)\,\sch_{v_1}(y;z)$, but this expression is not manifestly positive in the differences $y_i-z_j$, and we do not know a positive combinatorial rule for them.
\end{remark}

\section{The quantum Pieri formula} \label{section:pieri}

\begin{definition}[$\tom{k}$, $\fx_{k}(u,w)$, $q_k(u,w)$, $E_{p,k}(y;z)$, $\wte(y_A;z)$] \label{definition:pieri}
We define a relation $\tom{k}$ on pairs of permutations, introduced in its non-quantum form by Sottile \cite{sottile}, which was used originally to indicate when the coefficient of a Schubert polynomial is nonzero in a suitable application of the Pieri formula. In the generalization to double Schubert polynomials in \cite{samuelmolev}, the same relation applies to all nonzero coefficients in the expansion of the generalized Pieri formula. The relation plays the same role in our case. 

Given $u,w\in S_\infty$ and a positive integer $k$, we declare that $u\tom{k} w$ if there exists a $p$ with $0\leq p\leq k$ such that there are transpositions $t_{a_1b_1},\ldots,t_{a_pb_p}$ satisfying 
\begin{enumerate}
\item $a_i\leq k<b_i$ for all $i$
\item $a_i\neq a_j$ whenever $i\neq j$
\item $b_i\geq b_j$ whenever $i<j$
\item We have either
$$\ell(ut_{a_1b_1}\cdots t_{a_ib_i})=\ell(ut_{a_1b_1}\cdots t_{a_{i-1}b_{i-1}})+1$$ 
or 
$$\ell(ut_{a_1b_1}\cdots t_{a_ib_i})=\ell(ut_{a_1b_1}\cdots t_{a_{i-1}b_{i-1}})-2(b_i-a_i)+1$$
for all $i$
\item $w=ut_{a_1b_1}\cdots t_{a_pb_p}$
\end{enumerate}
In contrast to the non-quantum case, $u\tom{k} w$ does not imply that $\ell(u) \leq \ell(w)$, as the final condition can cause the length of the permutation to drop. For intuition, the condition $\ell(ut_{ab}) = \ell(u) - 2(b - a) + 1$ means that if we write the reflection $t_{ab}$ as
$$s_{b-1}\cdots s_{a+1}s_as_{a+1}\cdots s_{b-1}$$
then each simple reflection in this word cancels a corresponding simple reflection in a reduced word for $u$, as opposite to the length-increase case where there is a word for which precisely one simple reflection is inserted, regardless of the length of the reflection as a permutation.

If $u\tom{k} w$, define 
$$\fx_{k}(u,w)=\{u(i)\mid i\leq k\mbox{ and }u(i)=w(i)\}$$
and
$$n_k(u,w)=k-\lvert\fx_k(u,w)\rvert.$$
Also, we define $q_k(u,w)$ recursively. If $u=w$, define $q_k(u,w)=1$. If we have defined $q_k(u,w')$ with $w'=ut_{a_1b_1}\cdots t_{a_{p-1}b_{p-1}}$, suppose $\ell(w't_{a_pb_p})=\ell(w')+1$. Then define
$$q_k(u,w)=q_k(u,w')$$
If instead $\ell(w't_{a_pb_p})=\ell(w')-2(b_p-a_p)+1$, define
$$q_k(u,w)=q_k(u,w')q_{a_pb_p}$$

Suppose $y$ and $z$ are sets of variables and $0\leq p\leq k$ are integers. We define $E_{p,k}(y;z)$ to be the factorial elementary symmetric polynomial in $y$ and $z$. For specificity, this is the double Schubert polynomial $\sch_{c[p,k]}(y;z)$ where
$$c[p,k] = s_{k-p+1}\cdots s_{k-1}s_k$$
If $p=0$, the polynomial is identically $1$, and we set $E_{p,k}=0$ for $p<0$. If $p=k$, $E_{k,k}(y;z)$ can be simply written as
$$E_{k,k}(y;z)=\prod_{i=1}^k(y_i-z_1)$$
An explicit formula can be written for general $p$ without too much trouble as well, using e.g. \cite{molev1999littlewood}:
$$E_{p,k}(y;z)=\sum_{1\leq i_1<\cdots<i_p\leq k}\prod_{j=1}^p(y_{i_j}-z_{i_j-j+1})$$

For a set $A=\{i_1,\ldots,i_{|A|}\}$, we define
$$\wte(y_A;z)=\prod_{i\in A}(y_i-z_1)=E_{|A|,|A|}(y_{i_1},\ldots,y_{i_{|A|}};z),$$
the specialization of the factorial elementary polynomial to which the coefficients of the Pieri formula reduce. It will be used to express those coefficients.
\end{definition}

\begin{definition}[Quantum factorial elementary polynomials $\qelem_{p,k}(x;z)$] \label{def:qelem}
For $0\leq p\leq k$ we define
$$\qelem_{p,k}(x;z)=\qsch_{c[p,k]}(x;z),$$
the quantum double Schubert polynomial of the Grassmannian permutation $c[p,k]=s_{k-p+1}\cdots s_k$, whose diagram is a single column of height $p$ with descent at $k$. At $q=0$ this is $E_{p,k}(x;z)$. We also set $\qelem_{0,k}=1$ and $\qelem_{p,k}=0$ for $p<0$.
\end{definition}

\begin{lemma} \label{lemma:qelemkk}
For $k\geq1$,
$$\qelem_{k,k}(x;z)=\Delta_k(-z_1\mid x_1,\ldots,x_k)=\sum_{j=0}^{k}(-z_1)^{k-j}e_j^{(q)}(x_1,\ldots,x_k).$$
\end{lemma}
\begin{proof}
The permutation $c[k,k]=s_1s_2\cdots s_k$ has a unique reduced word, so its reduced factorizations are $v_1=s_1\cdots s_j$, $v_2=s_{j+1}\cdots s_k$ for $0\leq j\leq k$. By Lemma \ref{lemma:quantumtoordinary} with $y=0$, $\qelem_{k,k}(x;z)=\sum_{j}\sch_{s_1\cdots s_j}(0;z)\,\qsch_{s_{j+1}\cdots s_k}(x)$. Here $s_1\cdots s_j$ is dominant with code $(1^j)$, so $\sch_{s_1\cdots s_j}(x;z)=\prod_{i\leq j}(x_i-z_1)$ and $\sch_{s_1\cdots s_j}(0;z)=(-z_1)^j$; and $s_{j+1}\cdots s_k$ is the Grassmannian permutation with descent at $k$ and a single column of height $k-j$, whose Schubert polynomial $e_{k-j}(x_1,\ldots,x_k)$ is a standard elementary monomial, so that its quantum Schubert polynomial is $e^{(q)}_{k-j}(x_1,\ldots,x_k)$ \cite{fomin1997quantum}. Summing gives the second expression, which equals $\Delta_k(-z_1\mid x_1,\ldots,x_k)$ by the definition of the quantum elementary polynomials.
\end{proof}

\begin{theorem}[The Pieri formula] \label{theorem:pieri}
Suppose $u\in S_\infty$ and $k\geq 1$. Then
$$\qsch_u(x;y)\qelem_{k,k}(x;z)=\sum_{u\tom{k} w}q_k(u,w)\wte(y_{\fx_{k}(u,w)};z)\qsch_w(x;y).$$
\end{theorem}

The proof is given in Section \ref{section:proofpieri}. Note that the right-hand side is manifestly triple Schubert positive, in accordance with Theorem \ref{theorem:positivity}, and that at $q=0$ it reduces to the case $p=k$ of the Pieri formula of \cite[Theorem~7.1]{samuelmolev}. The general-degree version, Theorem \ref{theorem:pierigeneral}, replaces $\wte(y_{\fx_k(u,w)};z)$ by $E_{p-n_k(u,w),\,k-n_k(u,w)}(y_{\fx_k(u,w)};z)$.

Triple Schubert positivity does not by itself imply Graham positivity at $z=y$: a difference $y_i-y_j$ with $i<j$ is not a nonnegative combination of the negative simple roots $y_{l+1}-y_l$. The Pieri coefficients nevertheless are Graham-positive after specialization, for the same reason as in \cite[Section~7]{samuelmolev}: they are factorial elementary polynomials evaluated at an increasing sequence of $y$-variables.

\begin{corollary}[Graham positivity of the Pieri coefficients] \label{corollary:pierigraham}
Let $u\tom{k}w$, $n=n_k(u,w)$ and $0\leq p\leq k$. The specialization $E_{p-n,\,k-n}(y_{\fx_k(u,w)};z)\big|_{z=y}$ is a polynomial in the differences $y_{l+1}-y_l$ with nonnegative integer coefficients. Consequently, at $z=y$, every coefficient in Theorems \ref{theorem:pieri} and \ref{theorem:pierigeneral} lies in $\mathbb Z_{\geq0}[y_{l+1}-y_l][q_1,q_2,\ldots]$.
\end{corollary}
\begin{proof}
Write $\fx_k(u,w)=\{a_1<a_2<\cdots<a_{k-n}\}$. Since the $a_i$ are distinct positive integers, $a_i\geq i$ for all $i$. By the formula of Definition \ref{definition:pieri},
$$E_{p-n,\,k-n}(y_{\fx_k(u,w)};z)=\sum_{1\leq i_1<\cdots<i_{p-n}\leq k-n}\ \prod_{j=1}^{p-n}\big(y_{a_{i_j}}-z_{i_j-j+1}\big),$$
and at $z=y$ each factor becomes $y_{a_{i_j}}-y_{i_j-j+1}$ with $a_{i_j}\geq i_j\geq i_j-j+1$, which is the sum $\sum_{l=i_j-j+1}^{a_{i_j}-1}(y_{l+1}-y_l)$ of negative simple roots (empty, hence $0$, when $a_{i_j}=i_j-j+1$). The product and the sum have nonnegative integer coefficients, and $q_k(u,w)$ is a monomial.
\end{proof}

\newcommand{\lom}{\widetilde{\lambda}}
\newcommand{\mom}{\widetilde{\mu}}
\newcommand{\thop}{\widetilde{\theta}}
\newcommand{\medperm}{\overline{\mu}}
\newcommand{\medtheta}{\overline{\theta}}
\newcommand{\Tom}[1]{\xRightarrow{#1}}
\newcommand{\dom}{\mu}

\section{Expansion formulas and the multiplication algorithms} \label{section:algorithms}

\subsection{Dominant shapes}

\begin{definition}[$\code$, (strict) dominant permutations, $\theta$, $\dom$, $\thop$, $\mom$, medium theta $\medtheta$, $\medperm$] \label{def:codetheta}
We say that a permutation $\mu$ is \emph{dominant} if $\code_i(\mu)\geq \code_{i+1}(\mu)$ for all $i$. A permutation $\mu$ is said to be \emph{strict dominant} if $\mu$ is dominant and whenever $\code_i(\mu)\neq 0$ we have $\code_i(\mu)>\code_{i+1}(\mu)$. For a permutation $v\in S_\infty$, we define a dominant permutation $\dom(v)$ recursively as follows. If $v$ is already dominant, define $\mu(v)=v$. Otherwise, let $i$ be the maximal index such that $\code_i(v)<\code_{i+1}(v)$, and define
$$\dom(v)=\dom(vs_i)$$
Define $\theta_i(v)$ to be $\code_i(\dom(v))$.

Similarly, for a permutation $v\in S_\infty$, we define a strict dominant permutation $\mom(v)$ recursively as follows. If $v$ is already strict dominant, define $\mom(v)=v$. Otherwise, let $i$ be the maximal index such that $\code_{i+1}(v)\neq 0$ and $\code_i(v)\leq \code_{i+1}(v)$, then define
$$\mom(v)=\mom(vs_i)$$
Define $\thop_i(v)$ to be $\code_i(\mom(v))$.

Between these two lies the \emph{medium theta} $\medtheta(v)$, the shape actually iterated by the quantum multiplication algorithm (the routine \texttt{medium\_theta} of \cite{schubmult}). Call an index $i\geq1$ \emph{defective} for $v$ if either
\begin{enumerate}
\item[(a)] $\code_i(v)<\code_{i+1}(v)$, or
\item[(b)] $\code_i(v)=\code_{i+1}(v)\neq0$, $i>1$, and $\code_{i-1}(v)\leq\code_i(v)+1$.
\end{enumerate}
If $v$ has no defective index, define $\medperm(v)=v$; otherwise let $i$ be the \emph{minimal} defective index and define
$$\medperm(v)=\medperm(vs_i).$$
Define $\medtheta_i(v)=\code_i(\medperm(v))$ and $|\medtheta(v)|=\sum_i\medtheta_i(v)$. In both cases (a) and (b) one has $\code_i(v)\leq\code_{i+1}(v)$, which forces $v(i)<v(i+1)$; hence $\ell(vs_i)=\ell(v)+1$, $\code(vs_i)$ is obtained from $\code(v)$ by replacing $(\code_i,\code_{i+1})$ with $(\code_{i+1}+1,\code_i)$ (in case (b): by incrementing $\code_i$), and the recursion terminates, since a defective index has $\code_{i+1}(v)\neq0$, so that $v\in S_n$ implies $vs_i\in S_n$, and the length increases at each step. By construction $\medtheta(v)$ is weakly decreasing, so $\medperm(v)$ is dominant; a pair of equal consecutive nonzero parts $\medtheta_i(v)=\medtheta_{i+1}(v)$ occurs only when $i=1$ or $\medtheta_{i-1}(v)\geq\medtheta_i(v)+2$, in particular no three consecutive nonzero parts are equal; and $\medperm(v)=v\cdot y$ is a reduced factorization with $\ell(\medperm(v))=|\medtheta(v)|$. Thus $\medtheta(v)$ permits some, but not all, repeated parts, and lies between $\theta(v)$ (arbitrary repeats) and $\thop(v)$ (no repeats). The number of nonzero parts of $\medtheta(v)$ equals the number of nonzero parts of $\code(v)$, since neither (a) nor (b) creates or destroys a zero. We say that a dominant permutation $\mu$ has \emph{medium shape} if its code satisfies the condition just stated: whenever two consecutive nonzero parts $\code_i(\mu)=\code_{i+1}(\mu)$ are equal, either $i=1$ or $\code_{i-1}(\mu)\geq\code_i(\mu)+2$.
\end{definition}

The three shapes are related by the right weak order, in which $u\leq_Rw$ means $w=u\cdot y$ reduced, equivalently $N(u)\subseteq N(w)$ for the sets of \emph{value inversions} $N(u)=\{(a,b):a<b,\ u^{-1}(b)<u^{-1}(a)\}$ \cite[Proposition~3.1.3]{combcox}. Recall that a permutation is dominant if and only if it avoids the pattern $132$; it is strict dominant if and only if it is decreasing and then increasing, i.e.\ if and only if it avoids both $132$ and $231$ (a strict dominant permutation has $\mu(i)>\mu(i+1)$ for all $i$ up to its last descent, and conversely a permutation without peaks has no $132$ and no $231$).

\begin{lemma} \label{lemma:forcedmove}
Let $v\in S_\infty$ and let $i\geq1$ be such that $v(i)<v(i+1)$ and $\code_{i+1}(v)\neq0$. Then every $w\geq_Rv$ avoiding $132$ and $231$ satisfies $w\geq_Rvs_i$. If moreover $\code_i(v)<\code_{i+1}(v)$, then every $w\geq_Rv$ avoiding $132$ satisfies $w\geq_Rvs_i$.
\end{lemma}
\begin{proof}
Put $a=v(i)<b=v(i+1)$; then $N(vs_i)=N(v)\cup\{(a,b)\}$, so we must show that $b$ precedes $a$ in $w$. Since $\code_{i+1}(v)\neq0$ there is a value $c<b$ at a position $j>i+1$ of $v$, and $(c,b)\in N(v)\subseteq N(w)$ says that $b$ precedes $c$ in $w$. Suppose that $a$ precedes $b$ in $w$. If $c<a$, then $(c,a)\in N(v)\subseteq N(w)$, so $w$ contains $a,b,c$ in this order, a $231$ pattern. If $a<c<b$, then $w$ contains $a,b,c$ in this order, a $132$ pattern. In the first assertion both cases are excluded. For the second, note that $\code_{i+1}(v)-\code_i(v)=\#\{j>i+1:a<v(j)<b\}$ when $v(i)<v(i+1)$, so $\code_i(v)<\code_{i+1}(v)$ allows us to choose $c$ with $a<c<b$, and only the $132$ case arises.
\end{proof}

\begin{proposition}[Weak order characterization] \label{prop:weakorder}
Let $v\in S_\infty$ and let $\mathcal P=\mathcal P(v)$ be the set of indices $i$ with $\medtheta_i(v)=\medtheta_{i+1}(v)\neq0$. Then:
\begin{enumerate}
\item[(i)] $\dom(v)$ is the minimum, for $\leq_R$, of the dominant permutations $w\geq_Rv$. (This is the upper projection of $v$ to the Tamari lattice, realized as the set of $132$-avoiding permutations in the weak order \cite{bjorner1997shellable}, i.e.\ to the quotient of the weak order by the sylvester congruence \cite{reading2004congruences,reading2006cambrian}.)
\item[(ii)] $\mom(v)$ is the minimum of the strict dominant permutations $w\geq_Rv$.
\item[(iii)] The simple reflections $s_i$, $i\in\mathcal P$, commute pairwise, and
$$\mom(v)=\medperm(v)\cdot\prod_{i\in\mathcal P}s_i$$
is a reduced factorization. Equivalently, $\thop(v)$ is obtained from $\medtheta(v)$ by increasing by one the first entry of each retained equal pair, and $|\thop(v)|=|\medtheta(v)|+|\mathcal P(v)|$.
\end{enumerate}
In (i) and (ii) the recursions of Definition \ref{def:codetheta} may take the defective indices in any order; $\medperm(v)$, by contrast, depends on the order (for $v=2341$ the permutations $4231=\medperm(v)$ and $3421$ are both minimal among the medium shape permutations above $v$), and it is normalized by the choice of the minimal defective index.
\end{proposition}
\begin{proof}
(i) Each step $v\mapsto vs_i$ of the recursion for $\dom$ has $\code_i(v)<\code_{i+1}(v)$, so by Lemma \ref{lemma:forcedmove} every dominant $w\geq_Rv$ satisfies $w\geq_Rvs_i$; by induction on $\ell(\dom(v))-\ell(v)$, $w\geq_R\dom(v)$, and $\dom(v)$ is dominant and $\geq_Rv$. (ii) is identical, each step of the recursion for $\mom$ having $\code_i(v)\leq\code_{i+1}(v)$, i.e.\ $v(i)<v(i+1)$, and $\code_{i+1}(v)\neq0$, so that the first assertion of Lemma \ref{lemma:forcedmove} applies. (iii) Put $\mu=\medperm(v)$ and $w=\mom(v)$. Each step of the recursion for $\medperm$ is of type (a) or (b) of Definition \ref{def:codetheta}, and in both cases $v(i)<v(i+1)$ and $\code_{i+1}(v)\neq0$; so Lemma \ref{lemma:forcedmove} and induction give $w\geq_R\mu$. The indices in $\mathcal P$ are pairwise at distance $\geq2$, since no three consecutive parts of $\medtheta(v)$ are equal, so the $s_i$ commute. For $i\in\mathcal P$ we have $\mu(i)<\mu(i+1)$ (as $\code_i(\mu)=\code_{i+1}(\mu)$) and $\code_{i+1}(\mu)\neq0$; applying Lemma \ref{lemma:forcedmove} successively to $\mu,\ \mu s_{i_1},\ \mu s_{i_1}s_{i_2},\ldots$ (the positions $i,i+1$ and the entry $\code_{i+1}$ being unaffected by the previous $s_{i'}$, $|i-i'|\geq2$) gives $w\geq_R\mu':=\mu\prod_{i\in\mathcal P}s_i$. Each $s_i$ increases the length, replacing $(\code_i,\code_{i+1})=(k,k)$ by $(k+1,k)$, so $\mu'=\mu\cdot\prod s_i$ is reduced and $\code(\mu')$ is strictly decreasing on its nonzero part (the part preceding a retained pair is $\geq k+2>k+1$): $\mu'$ is strict dominant and $\mu'\geq_R\mu\geq_Rv$, whence $w\leq_R\mu'$ by (ii). Thus $w=\mu'$. For the final remark, $N(2341)=\{(1,2),(1,3),(1,4)\}$ is contained in $N(4231)$ and in $N(3421)$, both of which are medium shape with codes $(3,1,1)$ and $(2,2,1)$, and neither is below the other.
\end{proof}

\subsection{Expansion of a quantum double Schubert polynomial}

\begin{definition}[$\xRightarrow{i}$, $\ell(u,v)$, $\sigma_i(u,v)$, $D_i(u,v)$, $z_D$] \label{def:itom}
Let $u,v\in S_\infty$ and $i\geq 1$ be an integer. We define $u\Tom{i} v$ if there exist distinct integers $b_1,b_2,\ldots,b_k$ such that $b_j\neq i$ for all $j$, $\ell(ut_{i,b_1}\cdots t_{i,b_j})=\ell(u)+j$ for all $j$, and 
$$v=ut_{i,b_1}\cdots t_{i,b_k}.$$
Thus $u\Tom{i}v$ means that $v$ is reached from $u$ by a chain of Bruhat covers $u\lessdot ut_{ib_1}\lessdot ut_{ib_1}t_{ib_2}\lessdot\cdots\lessdot v$, each of which exchanges the position $i$ with a position not used before. The sequence $b_1,\ldots,b_k$ is determined by $u$ and $v$: the product $u^{-1}v=t_{ib_1}\cdots t_{ib_k}$ is the cycle $i\mapsto b_k\mapsto b_{k-1}\mapsto\cdots\mapsto b_1\mapsto i$, from which the $b_j$ are read off. We write $\ell(u,v)=\ell(v)-\ell(u)=k$ and define
$$\sigma_i(u,v)=\#\{j\mid b_j<i\},\qquad D_i(u,v)=\{u(i)\}\cup\{u(j)\mid u(j)\neq v(j)\}=\{u(i),u(b_1),\ldots,u(b_k)\},$$
a set of $k+1$ values. Finally, for a finite set $D=\{d_1<\cdots<d_m\}$ of positive integers with $m=k-p+1$ we write $E_{p,k}(x;z_D)=E_{p,k}(x;z_{d_1},\ldots,z_{d_m})$ and $\qelem_{p,k}(x;z_D)=\qelem_{p,k}(x;z_{d_1},\ldots,z_{d_m})$; since these polynomials are symmetric in $z_1,\ldots,z_{k-p+1}$ (Lemma \ref{lemma:Ehe}), the order in which the elements of $D$ are listed is immaterial. As in Definition \ref{def:qelem}, $E_{p,k}=\qelem_{p,k}=0$ when $p<0$.
\end{definition}

The following formulas express a (quantum) double Schubert polynomial as a signed sum of products of (quantum) factorial elementary polynomials in nested initial alphabets; they are the $v$-side input of the multiplication algorithm, and they are proved in Section \ref{section:proofvside}. All of them start from a dominant permutation $\mu$ lying above $v^{-1}$ in the right weak order, that is $\mu=v^{-1}\cdot y$ with $\ell(\mu)=\ell(v)+\ell(y)$; the permutations $\dom(v^{-1})$, $\mom(v^{-1})$ and $\medperm(v^{-1})$ are of this kind, since each step of their construction multiplies on the right by a simple reflection $s_i$ with $\code_i\leq\code_{i+1}$, which forces $v(i)<v(i+1)$ and so increases the length. The implementation uses $\mu=\dom(v^{-1})$ in the classical case, and the parts of $\code(\mu)$ are then the $\theta_i(v^{-1})$ of Definition \ref{def:codetheta}. In the quantum case the dominant permutation cannot be arbitrary, because the quantum polynomial $\qsch_{\mu^{-1}}(x;z)$ which is being expanded factors into a product of quantum factorial elementary polynomials only when the column heights of the diagram of $\mu^{-1}$, which are the parts of $\code(\mu)$, are pairwise distinct (Lemma \ref{lemma:dominantfactor}); for repeated column heights it is \emph{not} a product (see Section \ref{section:additional}). Theorem \ref{theorem:quandouble} accordingly uses the strict shape $\mom(v^{-1})$. The implementation (the routine \texttt{schubmult\_q\_double\_fast} of \cite{schubmult}) uses the smaller intermediate shape $\medperm(v^{-1})$, in which an equal pair of columns is allowed when there is room for it; the quantization of such a pair is a single ``doubled'' factor, which Theorem \ref{theorem:quandoublemed} isolates and which the algorithm multiplies out by the doubled Pieri formula of Theorem \ref{theorem:doubledpieri} and Corollary \ref{corollary:lowerdegrees}.

\begin{theorem}[Formula for ordinary double Schubert polynomials] \label{theorem:orddouble}
Let $v\in S_\infty$ and let $\mu$ be a dominant permutation with $\mu=v^{-1}\cdot y$ reduced, for instance $\mu=\dom(v^{-1})$. Let $\theta=\code(\mu)$, with nonzero parts $\theta_1\geq\cdots\geq\theta_k>0$. Then
$$\sch_v(x;z)=\sum_{(v_0,\ldots,v_k)}\prod_{i=1}^k (-1)^{\sigma_i(v_{i-1},v_i)}E_{\theta_i-\ell(v_{i-1},v_i),\,\theta_i}\big(x;z_{D_i(v_{i-1},v_i)}\big),$$
the sum being over all sequences $1=v_0\Tom{1}v_1\Tom{2}v_2\Tom{3}\cdots\Tom{k}v_k=v\mu$.
\end{theorem}

\begin{theorem}[Formula for quantum double Schubert polynomials] \label{theorem:quandouble}
Let $v\in S_\infty$ and let $\mu$ be a strict dominant permutation with $\mu=v^{-1}\cdot y$ reduced, for instance $\mu=\mom(v^{-1})$. Let $\thop=\code(\mu)$, with nonzero parts $\thop_1>\cdots>\thop_k>0$. Then
$$\qsch_v(x;z)=\sum_{(v_0,\ldots,v_k)}\prod_{i=1}^k (-1)^{\sigma_i(v_{i-1},v_i)}\qelem_{\thop_i-\ell(v_{i-1},v_i),\,\thop_i}\big(x;z_{D_i(v_{i-1},v_i)}\big),$$
the sum being over all sequences $1=v_0\Tom{1}v_1\Tom{2}\cdots\Tom{k}v_k=v\mu$.
\end{theorem}

For the medium shape we need a name for the quantization of a product of two factorial elementary polynomials in the same $x$-variables. Let $\mathcal Q$ be the quantization map of Fomin, Gelfand and Postnikov (Lemma \ref{lemma:quantumtoordinary}), extended linearly over all coefficient variables, and for $0\leq p_1,p_2\leq k$ and two alphabets $z,z'$ put
\begin{equation}\label{eq:qelemdouble}
\qelem_{p_1,p_2;k}(x;z,z'):=\mathcal Q\big(E_{p_1,k}(x;z)\,E_{p_2,k}(x;z')\big),
\end{equation}
the \emph{doubled quantum factorial elementary polynomial}; for $p_1=p_2=k$ it is the quantum double Schubert polynomial $\qsch_{d_2[k]}$ of the $k\times2$ rectangle in the variables $z_1,z'_1$ (Section \ref{section:additional}). It is \emph{not} the product $\qelem_{p_1,k}(x;z)\qelem_{p_2,k}(x;z')$ in general.

\begin{theorem}[Formula for quantum double Schubert polynomials, medium shape] \label{theorem:quandoublemed}
Let $v\in S_\infty$ and let $\mu$ be a dominant permutation of medium shape (Definition \ref{def:codetheta}) with $\mu=v^{-1}\cdot y$ reduced, for instance $\mu=\medperm(v^{-1})$. Let $\medtheta=\code(\mu)$, with nonzero parts $\medtheta_1\geq\cdots\geq\medtheta_m>0$, let $\mathcal P=\{i:\medtheta_i=\medtheta_{i+1}\}$ be the set of \emph{paired} indices and $\mathcal S=\{1,\ldots,m\}\setminus(\mathcal P\cup(\mathcal P+1))$ the set of \emph{single} indices. Then
$$\qsch_v(x;z)=\sum_{(v_0,\ldots,v_m)}\ \prod_{i=1}^m(-1)^{\sigma_i(v_{i-1},v_i)}\ \prod_{i\in\mathcal S}\qelem_{\medtheta_i-\ell_i,\,\medtheta_i}\big(x;z_{D_i}\big)\ \prod_{i\in\mathcal P}\qelem_{\medtheta_i-\ell_i,\,\medtheta_i-\ell_{i+1};\,\medtheta_i}\big(x;z_{D_i},z_{D_{i+1}}\big),$$
where $\ell_i=\ell(v_{i-1},v_i)$, $D_i=D_i(v_{i-1},v_i)$, and the sum is over all sequences $1=v_0\Tom{1}v_1\Tom{2}\cdots\Tom{m}v_m=v\mu$. For a strict dominant $\mu$ this is Theorem \ref{theorem:quandouble}.
\end{theorem}

By Proposition \ref{prop:weakorder}(iii), $|\thop(v^{-1})|=|\medtheta(v^{-1})|+|\mathcal P|$: the medium shape is smaller than the strict shape by exactly the number of retained pairs (for instance $\code(v^{-1})=(1,1)$ gives $\medtheta=(1,1)$, $\thop=(2,1)$, and $\code(v^{-1})=(2,2,2)$ gives $\medtheta=(4,2,2)$, $\thop=(4,3,2)$). Since $\sum_i\ell_i=|\code(\mu)|-\ell(v)$ for every chain in these formulas, a smaller shape means shorter chains in the $v$-side sum and, in the multiplication algorithm, $u$-side Pieri steps of smaller degree; this is the reason the implementation prefers $\medperm(v^{-1})$ to $\mom(v^{-1})$. Proposition \ref{prop:weakorder}(iii) also explains how Theorem \ref{theorem:quandoublemed} arises from Theorem \ref{theorem:quandouble}: the dominant polynomial of the medium shape is obtained from that of the strict shape by one divided difference per retained pair, which turns the product of the two adjacent factors into the doubled factor (Corollary \ref{corollary:mediumfactor}).

\subsection{The classical algorithm}

Coupled with the Pieri formulas, these formulas give an algorithm for computing products of (quantum) double Schubert polynomials, which is the one implemented in \cite{schubmult}. (Algorithms for computing quantum Schubert polynomials and Gromov--Witten invariants of partial flag manifolds were given by Buch \cite{buch2005partial}.) We describe its classical form first. The $u$-side covers are generated by the classical Pieri rule, Theorem \ref{theorem:pierigeneral} at $q=0$ (the enumeration is Method \ref{method:qpieri} below with the quantum covers discarded); the $v$-side is the following signed path expansion.

The starting point of the $v$-side primitive is the structure of a \emph{dominant} double Schubert polynomial. If $\mu$ is dominant and the columns of its diagram have heights $\lambda_1,\ldots,\lambda_m$ (so that $(\lambda_1,\ldots,\lambda_m)=\code(\mu^{-1})$), then $\sch_\mu(x;z)$ factors as a product of factorial elementary symmetric polynomials of maximal degree, one per column of the shape,
$$\sch_\mu(x;z)=\prod_{i=1}^{m}E_{\lambda_i,\,\lambda_i}(x;z_i)$$
(Lemma \ref{lemma:dominantfactor}). For a general $v$ one writes $\sch_v(x;z)$ by applying divided differences in the $z$-variables to the dominant polynomial $\sch_{\mu^{-1}}(x;z)$ with $\mu=\dom(v^{-1})$, whose column heights are the parts of $\theta(v^{-1})$; because that polynomial factors by column, the divided-difference action decouples and may be carried out one column (one variable $z_i$) at a time. Applied to a single column this produces a \emph{signed} descending expansion. The following primitive carries it out.

\begin{definition}[$k$-strip descent beneath an edge, the down-Pieri problem] \label{def:kdownproblem}
Let $w,\mu\in S_n$ and $1\leq k<n$. Say $w'$ lies \emph{beneath the edge} $\mu$ if
$$\mathrm{inv}(w'\mu)=\mathrm{inv}(\mu)-\mathrm{inv}(w'),$$
that is, $\ell(w'\mu)=\ell(\mu)-\ell(w')$. A \emph{$k$-strip descent} of $w$ is a Bruhat descent $w\gtrdot wt_{ab}$ realized by a transposition $t_{ab}$ that moves the value in the pivot position $k$ either leftward (from a position $a<k$, with sign $-1$) or rightward (to a position $b>k$, with sign $+1$), touching only positions still holding their original value. Given a starting label $w$, the edge $\mu$, a degree $p$, and the column $k$, the \emph{down-Pieri problem} asks for the signed list of those length-$t$ descending $k$-strip chains from $w$ ($t\leq p$) whose endpoint lies beneath $\mu$.
\end{definition}

\begin{method}[$\textbf{\textsc{kdown}}(w,\mu,p,k)$] \label{method:kdown}
The algorithm returns a signed list $D$ of triples $(w',t,s)$, one for each length-$t$ descending $k$-strip chain from $w$ ($t\leq p$) whose endpoint $w'$ lies beneath $\mu$, with $s\in\{+1,-1\}$ its sign.
\begin{enumerate}
\item Initialize $D\gets[\,]$; if $w$ is already beneath $\mu$, append $(w,0,+1)$. Set a frontier $F\gets[(w,+1)]$.
\item For $t=1,\ldots,p$, form $F'\gets[\,]$ and, for each $(w'',s)\in F$:
  \begin{enumerate}
  \item (leftward) for each position $a<k$ with $w''(a)=w(a)$ such that $w''\gtrdot w''t_{a,k}$ is a Bruhat descent, set $w'''\gets w''t_{a,k}$ and append $(w''',-s)$ to $F'$; if $w'''$ is beneath $\mu$, append $(w''',t,-s)$ to $D$;
  \item (rightward) for each position $b>k$ with $w''(b)=w(b)$ such that $w''\gtrdot w''t_{k,b}$ is a Bruhat descent, set $w'''\gets w''t_{k,b}$ and append $(w''',+s)$ to $F'$; if $w'''$ is beneath $\mu$, append $(w''',t,+s)$ to $D$.
  \end{enumerate}
  Then set $F\gets F'$.
\item Return $D$.
\end{enumerate}
\end{method}

Applied without the edge $\mu$, the signed descending expansion produces many extraneous terms that ultimately cancel. The edge is exactly the bound that discards them: only endpoints beneath $\mu$ are retained in $D$. In the assembled $v$-side (Method~\ref{method:vpaths} below), only these retained labels seed the next column, so each stage passes forward only labels compatible with the remaining shape.

\begin{method}[$\textbf{\textsc{vpaths}}(v)$, the $v$-side path dictionaries] \label{method:vpaths}
Let $\theta:=\theta(v^{-1})$ with nonzero parts in columns $1,\ldots,L$, and let $\nu:=v\,\dom(v^{-1})$. Processing columns from $i=L$ down to $i=1$, maintain for each column a set of live labels (initialized at column $L$ to the single label $\nu$). At column $i$, with $k=i$ and edge $\mu_i:=\dom$ of the shape truncated below column $i$, call $\textbf{\textsc{kdown}}(w,\mu_i,\theta_i,k)$ for each live label $w$; record its signed output triples as the transitions out of $w$, and let their endpoints be the live labels seeding column $i-1$. The collected transitions form the path dictionaries used by $\textbf{\textsc{schubmult}}$.
\end{method}

\begin{method}[The algorithm \textbf{\textsc{schubmult}}] \label{method:schubmult}
Fix $n$ and let $u,v\in S_n$. The algorithm $\textbf{\textsc{schubmult}}(u,v)$ computes the expansion of $\sch_u(x;z)\,\sch_v(x;z)$ in the Schubert basis as follows. Set $\theta:=\theta(v^{-1})$ and let
$$L:=\#\{i\mid\theta_i\neq 0\}$$
be the number of nonzero parts. The algorithm proceeds in $L$ layers. It maintains a \emph{front}, a collection of nodes each consisting of an intermediate permutation on the $u$-side paired with a label on the $v$-side, initialized to the single node $(u,\,e)$. At layer $i$ it replaces each node by the nodes obtained by combining one $u$-side cover, generated by the classical Pieri rule, i.e.\ by $\textbf{\textsc{qpieri}}(\,\cdot\,,\theta_i,\theta_i)$ of Method~\ref{method:qpieri} below with the quantum covers discarded, with one $v$-side transition, generated by $\textbf{\textsc{vpaths}}(v)$ (Method~\ref{method:vpaths}), accumulating coefficients in the coefficient ring. After all $L$ layers the surviving coefficients at $v$-label $v\dom(v^{-1})$ give the structure constants.
\end{method}

\subsection{The quantum algorithm}
The same path organization gives the quantum multiplication algorithm, the routine \texttt{schubmult\_q\_double\_fast} of \cite{schubmult}. The $v$-side is unchanged except for the shape: the path dictionaries of Method~\ref{method:vpaths} are built over the columns of $\medtheta:=\medtheta(v^{-1})$ (Definition \ref{def:codetheta}), with $\nu:=v\,\medperm(v^{-1})$ in place of $v\dom(v^{-1})$, and by Theorem \ref{theorem:quandoublemed} they enumerate the chains $1=v_0\Tom1\cdots\Tom m v_m=\nu$ together with the data $\ell_i$, $\sigma_i$, $D_i$ of each step. On the $u$-side, the classical covers are replaced by the covers of the quantum Bruhat graph (Definition \ref{def:cover}, Remark~\ref{remark:qkbruhat}), including the length-decreasing moves, whose contributions carry the corresponding $q$-monomials. We first describe the two $u$-side primitives, which enumerate the chains of Theorem \ref{theorem:pierigeneral} and the compatible pairs of Corollary \ref{corollary:lowerdegrees}.

\begin{method}[$\textbf{\textsc{qpieri}}(u,p,k)$] \label{method:qpieri}
Given $u\in S_\infty$, $k\geq1$ and $0\leq p\leq k$, the algorithm returns a list $T$ of triples $(w,t,\mathbf q)$ with $w\in S_\infty$, $0\leq t\leq p$ and $\mathbf q$ a monomial in the $q_i$.
\begin{enumerate}
\item Initialize $T\gets[(u,0,1)]$ and a frontier $F\gets[(u,1,\infty)]$ of triples $(w,\mathbf q,b_{\mathrm{last}})$.
\item For $t=1,\ldots,p$, form $F'\gets[\,]$ and, for each $(w,\mathbf q,b_{\mathrm{last}})\in F$:
  \begin{enumerate}
  \item compute the set $A\gets\{a\leq k:\ w(a)=u(a)\}$ of positions $\leq k$ not yet moved;
  \item for each $b$ with $k<b\leq b_{\mathrm{last}}$, in decreasing order, and each $a\in A$ such that $t_{ab}$ is a cover of $w$ (Definition \ref{def:cover}: either $w(a)<w(b)$ and no $w(r)$, $a<r<b$, lies in $(w(a),w(b))$, or $w(a)>w(b)$ and every $w(r)$, $a<r<b$, lies in $(w(b),w(a))$), put $\mathbf q'=\mathbf q\,\kappa(w,t_{ab})$, i.e.\ $\mathbf q'=\mathbf q$ for a classical cover and $\mathbf q'=\mathbf q\,q_{ab}=\mathbf q\,q_a\cdots q_{b-1}$ for a quantum cover, and append $(wt_{ab},\mathbf q',b)$ to $F'$ and $(wt_{ab},t,\mathbf q')$ to $T$.
  \end{enumerate}
  Then set $F\gets F'$.
\item Return $T$.
\end{enumerate}
Only finitely many $b$ qualify in (b): if $w$ fixes every position $>N$ then no $t_{ab}$ with $b>\max(N,k)+1$ is a cover of $w$ (the values $w(r)=r$ at positions $\max(N,k)<r<b$ lie between $w(a)$ and $w(b)=b$), so $b$ ranges over $k<b\leq\min(b_{\mathrm{last}},\max(N,k)+1)$. The cover test in (b) is performed numerically: the implementation computes the signed length change
$$\ell(wt_{ab})-\ell(w)=\pm1\pm2\,\#\{r\mid a<r<b,\ w(r)\text{ between }w(a)\text{ and }w(b)\}$$
of Lemma \ref{lemma:lengthchange} (the routine \texttt{count\_bruhat}, with the sign $+$ if $w(a)<w(b)$ and $-$ otherwise) and accepts $t_{ab}$ exactly when this number is $1$ or $1-2(b-a)$, multiplying the current monomial by $q_a\cdots q_{b-1}$ in the second case. This is the routine \texttt{elem\_sym\_perms\_q} of \cite{schubmult}; at $q=0$, i.e.\ discarding the quantum covers, it enumerates the chains of the classical Pieri rule of \cite{sottile,samuelmolev} (Theorem \ref{theorem:pierigeneral} at $q=0$).
\end{method}

\begin{lemma}[Correctness of \textbf{\textsc{qpieri}}] \label{lemma:qpieri}
The list $T$ returned by $\textbf{\textsc{qpieri}}(u,p,k)$ contains exactly one triple $(w,t,\mathbf q)$ for each $w\in S_\infty$ with $u\tom kw$ and $n_k(u,w)\leq p$, and for it $t=n_k(u,w)$ and $\mathbf q=q_k(u,w)$. Consequently, for $0\leq p'\leq p$,
$$\qsch_u(x;y)\,\qelem_{p',k}(x;z)=\sum_{(w,t,\mathbf q)\in T,\ t\leq p'}\mathbf q\,E_{p'-t,\,k-t}\big(y_{\fx_k(u,w)};z\big)\,\qsch_w(x;y),$$
where $\fx_k(u,w)=\{u(a):a\leq k,\ u(a)=w(a)\}$ is read off from $u$ and $w$.
\end{lemma}
\begin{proof}
A triple $(w,t,\mathbf q)$ is appended to $T$ for each sequence of transpositions $t_{a_1b_1},\ldots,t_{a_tb_t}$ generated by the loop, and by construction these sequences are exactly those satisfying conditions (1)--(5) of Definition \ref{definition:pieri} with $p$ replaced by $t\leq p$: $a_i\leq k<b_i$, the $a_i$ are distinct (a position leaves $A$ once it is moved and, being moved only once, never returns, since $w(a)\neq u(a)$ after the move), the $b_i$ are weakly decreasing ($b\leq b_{\mathrm{last}}$), and each step is a cover (condition (4), by Definition \ref{def:cover}). The weight $\mathbf q$ is the product of the $\kappa$'s, which is $q_k(u,w)$ by the recursive definition in Definition \ref{definition:pieri}. It remains to see that each $w$ with $u\tom kw$ arises from exactly one such sequence. By Lemma \ref{lemma:piericycle}, a witnessing sequence is a pivot chain of $\pi=u^{-1}w$ in which the factors with a common $b$ are consecutive and the distinct values of $b$ decrease; the factors with a common $b$ form the pivot factorization \eqref{eq:pivot} of the cycle of $\pi$ with top $b$, whose order is determined by that cycle, and the order of the cycles is determined by their tops. Hence the sequence is determined by $w$, and $t=n_k(u,w)$ is its length. The displayed formula is Theorem \ref{theorem:pierigeneral} with $p'$ in place of $p$, the sum over $u\tom kw$ with $n_k(u,w)\leq p'$ being the sum over the triples with $t\leq p'$.
\end{proof}

\begin{method}[$\textbf{\textsc{doublepieri}}(u,p_1,p_2,k)$] \label{method:doublepieri}
Given $u\in S_\infty$, $k\geq1$ and $0\leq p_1,p_2\leq k$, the algorithm returns, for each triple $(v,t_1,\mathbf q_1)$ of $\textbf{\textsc{qpieri}}(u,p_1,k)$, a list $L(v,t_1,\mathbf q_1)$ of triples $(w,t_2,\mathbf q_2)$.
\begin{enumerate}
\item Compute $T_1\gets\textbf{\textsc{qpieri}}(u,p_1,k)$.
\item For each $(v,t_1,\mathbf q_1)\in T_1$: compute the cycle decomposition of $\pi_1=u^{-1}v$, and record for each top $B$ the set $\mathrm{supp}\,c_B\setminus\{B\}$ of non-top elements of the cycle $c_B$ of $\pi_1$ with top $B$ (empty if there is none); compute $T_2\gets\textbf{\textsc{qpieri}}(v,p_2,k)$; and for each $(w,t_2,\mathbf q_2)\in T_2$ compute the cycles of $\pi_2=v^{-1}w$ and append $(w,t_2,\mathbf q_2)$ to $L(v,t_1,\mathbf q_1)$ unless some cycle $d$ of $\pi_2$ has a non-top element in $\mathrm{supp}\,c_{\max d}\setminus\{\max d\}$.
\item Return the lists $L$.
\end{enumerate}
This is the routine \texttt{double\_elem\_sym\_q} of \cite{schubmult}.
\end{method}

\begin{lemma}[Correctness of \textbf{\textsc{doublepieri}}] \label{lemma:doublepieri}
The pairs $\big((v,t_1,\mathbf q_1),(w,t_2,\mathbf q_2)\big)$ with $(w,t_2,\mathbf q_2)\in L(v,t_1,\mathbf q_1)$ returned by $\textbf{\textsc{doublepieri}}(u,p_1,p_2,k)$ are in bijection with the compatible pairs $u\tom kv\tom kw$ with $n_k(u,v)\leq p_1$ and $n_k(v,w)\leq p_2$, with $t_1=n_k(u,v)$, $t_2=n_k(v,w)$, $\mathbf q_1\mathbf q_2=q_k(u,v)q_k(v,w)$. Consequently, for $0\leq p_1'\leq p_1$ and $0\leq p_2'\leq p_2$ and two alphabets $z,z'$,
$$\qsch_u(x;y)\,\qelem_{p_1',p_2';k}(x;z,z')=\sum_{\substack{(v,t_1,\mathbf q_1)\in T_1\\ (w,t_2,\mathbf q_2)\in L(v,t_1,\mathbf q_1)\\ t_1\leq p_1',\ t_2\leq p_2'}}\mathbf q_1\mathbf q_2\,E_{p_1'-t_1,\,k-t_1}\big(y_{\fx_k(u,v)};z\big)\,E_{p_2'-t_2,\,k-t_2}\big(y_{\fx_k(v,w)};z'\big)\,\qsch_w(x;y).$$
\end{lemma}
\begin{proof}
By Lemma \ref{lemma:qpieri} the pairs examined in step 2 are in bijection with the pairs of chains $u\tom kv\tom kw$ with $n_k(u,v)\leq p_1$, $n_k(v,w)\leq p_2$, with the stated $t$'s and $\mathbf q$'s. The rejection test is the negation of Definition \ref{def:compatible}: a cycle $d$ of $v^{-1}w$ whose top $\max d$ is the top of a cycle $c$ of $u^{-1}v$ must satisfy $c\cap d=\{\max d\}$, i.e.\ no non-top element of $d$ may lie in $c$ (it cannot equal $\max d=\max c$), and cycles of $v^{-1}w$ whose top is not a top of $u^{-1}v$ are unrestricted. The formula is Corollary \ref{corollary:lowerdegrees} with $(p_1',p_2')$ in place of $(p_1,p_2)$.
\end{proof}

The implementation calls $\textbf{\textsc{qpieri}}$ and $\textbf{\textsc{doublepieri}}$ once per front node with $p$, respectively $(p_1,p_2)$, equal to the largest degrees $k-\ell_i$, respectively $(k-\ell_i,k-\ell_{i+1})$, occurring among the $v$-side transitions of the layer, and then evaluates the coefficient of each transition by the displayed formulas with the actual degrees $p'\leq p$ (the routine \texttt{elem\_sym\_func\_q}, which reads $\fx_k(u,w)$ off the positions $a\leq k$ with $u(a)=w(a)$ and the alphabet $z_{D_i}$ off the transition); by the two lemmas this is the same as calling them separately for each transition.

\begin{method}[The algorithm \textbf{\textsc{schubmult}}$_q$] \label{method:schubmultq}
Fix $n$ and $u,v\in S_n$, and let $\medtheta=\medtheta(v^{-1})$ with nonzero parts $\medtheta_1\geq\cdots\geq\medtheta_m>0$, paired indices $\mathcal P$ and single indices $\mathcal S$ as in Theorem \ref{theorem:quandoublemed}. The algorithm computes the expansion of $\qsch_u(x;y)\,\qsch_v(x;z)$ in the basis $\{\qsch_w(x;y)\}$. It maintains a front of nodes $(w,\lambda)$, a $u$-side permutation paired with a $v$-side label, initialized to $(u,e)$ with coefficient $1$, and runs through the indices $i=1,\ldots,m$, skipping $i\in\mathcal P+1$:
\begin{enumerate}
\item If $i\in\mathcal S$, with $k=\medtheta_i$: for each node $(w,\lambda)$ and each $v$-side transition $\lambda\to\lambda'$ of column $i$, with data $(\ell_i,\sigma_i,D_i)$, and each $(w',t,\mathbf q)\in\textbf{\textsc{qpieri}}(w,k-\ell_i,k)$, create the node $(w',\lambda')$ with coefficient
$$(-1)^{\sigma_i}\,\mathbf q\,E_{k-\ell_i-t,\,k-t}\big(y_{\fx_k(w,w')};z_{D_i}\big)$$
times the coefficient of $(w,\lambda)$; by Lemma \ref{lemma:qpieri} this expands $\qsch_w(x;y)\,\qelem_{k-\ell_i,k}(x;z_{D_i})$.
\item If $i\in\mathcal P$, with $k=\medtheta_i=\medtheta_{i+1}$: for each node $(w,\lambda)$ and each pair of consecutive $v$-side transitions $\lambda\to\lambda'\to\lambda''$ of the columns $i$ and $i+1$, with data $(\ell_i,\sigma_i,D_i)$ and $(\ell_{i+1},\sigma_{i+1},D_{i+1})$, and each pair $\big((w',t_1,\mathbf q_1),(w'',t_2,\mathbf q_2)\big)$ returned by $\textbf{\textsc{doublepieri}}(w,k-\ell_i,k-\ell_{i+1},k)$, create the node $(w'',\lambda'')$ with coefficient
$$(-1)^{\sigma_i+\sigma_{i+1}}\,\mathbf q_1\mathbf q_2\,E_{k-\ell_i-t_1,\,k-t_1}\big(y_{\fx_k(w,w')};z_{D_i}\big)\,E_{k-\ell_{i+1}-t_2,\,k-t_2}\big(y_{\fx_k(w',w'')};z_{D_{i+1}}\big)$$
times the coefficient of $(w,\lambda)$; by Lemma \ref{lemma:doublepieri} this expands $\qsch_w(x;y)\,\qelem_{k-\ell_i,k-\ell_{i+1};k}(x;z_{D_i},z_{D_{i+1}})$ over the compatible pairs $w\tom kw'\tom kw''$.
\end{enumerate}
Coefficients of coinciding nodes are added. After the last index the coefficient of $(w,\nu)$ is $\sum_dq^dc_{u,v}^{w,d}(y;z)$.
\end{method}

\begin{proposition}[Correctness of \textbf{\textsc{schubmult}}$_q$] \label{prop:schubmultqcorrect}
Method \ref{method:schubmultq} computes $\sum_dq^dc_{u,v}^{w,d}(y;z)$ for all $w$.
\end{proposition}
\begin{proof}
By Theorem \ref{theorem:quandoublemed}, $\qsch_v(x;z)$ is the sum over the $v$-side chains $1=v_0\Tom1\cdots\Tom mv_m=\nu$ of $\prod_i(-1)^{\sigma_i}$ times the product of the single factors $\qelem_{k-\ell_i,k}(x;z_{D_i})$, $i\in\mathcal S$, and the doubled factors $\qelem_{k-\ell_i,k-\ell_{i+1};k}(x;z_{D_i},z_{D_{i+1}})$, $i\in\mathcal P$; the path dictionaries enumerate these chains step by step, the label $\lambda$ of a node recording the current $v_i$. Multiplying $\qsch_u(x;y)$ by these factors one at a time, in the order of the indices, by Lemma \ref{lemma:qpieri} for the single factors and Lemma \ref{lemma:doublepieri} for the doubled ones, the coefficient of the node $(w,\lambda)$ after the index $i$ is the coefficient of $\qsch_w(x;y)$ in $\qsch_u(x;y)$ times the product of the factors of the columns $\leq i$, summed over the partial chains ending at $\lambda$. After the last index, the nodes $(w,\nu)$ carry the coefficient of $\qsch_w(x;y)$ in $\qsch_u(x;y)\qsch_v(x;z)$. The two $v$-side transitions of a paired layer must be consumed together, since the doubled factor is not a product of two quantum factorial elementary polynomials. (The requirement that a retained pair be preceded by a drop of at least two is what makes the lifted shape $\nu$ of Corollary \ref{corollary:mediumfactor} strict dominant, so that the two columns of the pair are the only ones whose quantization interacts.)
\end{proof}

The classical algorithm is the specialization $q=0$ with $\medtheta$ replaced by $\theta$, in which the doubled factor is the product $E_{p_1,k}E_{p_2,k}$ and the paired step may be split into two ordinary ones. Thus the two algorithms share the same layer-by-layer structure, while their Pieri transitions reflect the classical or quantum multiplication rule, respectively.

Neither algorithm depends on the particular dominant permutation used for the $v$-side, beyond the hypotheses of the expansion theorems: the classical algorithm is correct for any dominant $\mu$ lying above $v^{-1}$ in the right weak order, that is, with $\mu=v^{-1}\cdot y$ reduced, equivalently with $\mu^{-1}$ above $v$ in the left weak order (Theorem \ref{theorem:orddouble}); the quantum algorithm with single steps only is correct for any strict dominant such $\mu$ (Theorem \ref{theorem:quandouble}); and the quantum algorithm with paired steps is correct for any such $\mu$ of medium shape (Theorem \ref{theorem:quandoublemed}), strict dominant permutations being a special case. The canonical choices $\dom(v^{-1})$, $\mom(v^{-1})$ and $\medperm(v^{-1})$ are simply the smallest admissible ones, and the length of $\mu$ governs the number of layers and the sizes of the path dictionaries, so they are the most efficient choices; the medium shape is preferred to the strict one for the same reason.

\begin{example}[The algorithm $\textbf{\textsc{schubmult}}_q$ on $\qsch_{21}(x;y)\,\qsch_{312}(x;z)$] \label{example:algorithm}
Let $u=21$ and $v=312$. The $v$-side: $v^{-1}=231$ has code $(1,1)$, which is of medium shape (the equal pair is at $i=1$), so $\medperm(v^{-1})=v^{-1}$, $\medtheta=(1,1)$, $\mathcal P=\{1\}$, $\mathcal S=\emptyset$, and $\nu=v\,\medperm(v^{-1})=1$. The only chain $1=v_0\Tom1v_1\Tom2v_2=\nu$ is the trivial one, with $\ell_1=\ell_2=0$, $\sigma_1=\sigma_2=0$, $D_1=\{v_0(1)\}=\{1\}$ and $D_2=\{v_1(2)\}=\{2\}$; so Theorem \ref{theorem:quandoublemed} reads $\qsch_{312}(x;z)=\qelem_{1,1;1}(x;z_1,z_2)=\mathcal Q\big((x_1-z_1)(x_1-z_2)\big)$, and there is a single paired layer with $k=1$ and $(p_1,p_2)=(1,1)$.

The $u$-side: $\textbf{\textsc{qpieri}}(21,1,1)$ examines the transpositions $t_{1b}$, $b\geq2$, of $w=21$: $t_{12}$ has length change $-1=1-2(2-1)$ (a quantum cover, $w(1)=2>w(2)=1$ with nothing between) and gives $(12,1,q_1)$; $t_{13}$ has $w(1)=2<w(3)=3$ with the intermediate value $1\notin(2,3)$, length change $+1$, and gives $(312,1,1)$; $t_{14}$ has the intermediate value $3\in(2,4)$ and is not a cover. Hence
$$T_1=\{(21,0,1),\ (12,1,q_1),\ (312,1,1)\}.$$
$\textbf{\textsc{doublepieri}}(21,1,1,1)$ then computes, for each $(v',t_1,\mathbf q_1)\in T_1$, the list $\textbf{\textsc{qpieri}}(v',1,1)$ and discards the incompatible pairs:
\begin{itemize}
\item $v'=21$ (trivial first step, no cycles): all of $T_2=T_1$ is retained.
\item $v'=12$ (first step $t_{12}$, cycle $(1,2)$ with top $2$, non-top element $1$): $T_2=\{(12,0,1),(21,1,1)\}$, the second step $t_{13}$ not being a cover of $12$. The pair $(12,1,q_1),(21,1,1)$ is \emph{rejected}: the cycle $(1,2)$ of the second step has the top $2$ of a first-step cycle and the non-top element $1$ in common with it. Retained: $(12,0,1)$.
\item $v'=312$ (first step $t_{13}$, cycle $(1,3)$): $T_2=\{(312,0,1),\ (132,1,q_1),\ (4123,1,1)\}$, from the quantum cover $t_{12}$ ($3>1$) and the classical cover $t_{14}$ ($3<4$, intermediate values $1,2\notin(3,4)$); $t_{13}$ is not a cover of $312$ ($3>2$ but the intermediate value $1\notin(2,3)$). The second-step cycles $(1,2)$ and $(1,4)$ have tops $2,4\neq3$, so both pairs are retained.
\end{itemize}
Each retained pair $\big((v',t_1,\mathbf q_1),(w,t_2,\mathbf q_2)\big)$ contributes $\mathbf q_1\mathbf q_2\,E_{1-t_1,1-t_1}(y_{\fx_1(u,v')};z_1)\,E_{1-t_2,1-t_2}(y_{\fx_1(v',w)};z_2)$ to the coefficient of $\qsch_w(x;y)$, where $E_{1,1}(y_a;z_i)=y_a-z_i$ and $E_{0,0}=1$:
\begin{center}
\begin{tabular}{llll}
$v'$ & $w$ & $\mathbf q_1\mathbf q_2$ & coefficient contributed\\\hline
$21$ & $21$ & $1$ & $(y_2-z_1)(y_2-z_2)$\\
$21$ & $12$ & $q_1$ & $q_1(y_2-z_1)$\\
$21$ & $312$ & $1$ & $y_2-z_1$\\
$12$ & $12$ & $q_1$ & $q_1(y_1-z_2)$\\
$312$ & $312$ & $1$ & $y_3-z_2$\\
$312$ & $132$ & $q_1$ & $q_1$\\
$312$ & $4123$ & $1$ & $1$
\end{tabular}
\end{center}
Adding the coefficients of coinciding $w$,
$$\qsch_{21}(x;y)\,\qsch_{312}(x;z)=(y_2-z_1)(y_2-z_2)\,\qsch_{21}+q_1(y_1+y_2-z_1-z_2)+(y_2+y_3-z_1-z_2)\,\qsch_{312}+q_1\,\qsch_{132}+\qsch_{4123},$$
all polynomials on the right being in $(x;y)$. The rejected pair would have contributed $q_1\qsch_{21}$, which is exactly the difference $q_1\qsch_u\qsch_{d_2[0]}$ between $\qsch_u\,\qelem_{1,1}(x;z_1)\qelem_{1,1}(x;z_2)$ and $\qsch_u\,\qsch_{312}(x;z)$ predicted by Lemma \ref{lemma:thetaPP} below; see Example \ref{example:doubled}. Every coefficient is a nonnegative combination of products of differences $y_a-z_b$, as Theorem \ref{theorem:positivity} requires, and has the degree $\ell(u)+\ell(v)-\ell(w)-2|d|$ of Section \ref{section:definitions} (for instance $3-0-2=1$ for the term $q_1(y_1+y_2-z_1-z_2)$).
\end{example}

\subsection{Comparison with straightening methods} \label{section:benchmark}

The general-purpose alternative to a Pieri-driven algorithm is straightening: expand both factors in a basis of standard monomials of $QH^*(Fl_n)$, multiply, and reduce modulo the quantum ideal. The most recent form of this is the quantum rim-hook rule of Gu and Kalashnikov \cite{gu2020rimhook}, which for the complete flag variety works in the basis of Fomin, Gelfand and Postnikov \cite{fomin1997quantum} of standard elementary monomials $\prod_m e_{j_m}(x_1,\ldots,x_m)$, $0\leq j_m\leq m$ (under flag duality, their monomials in the Chern classes of the tautological quotient bundles), and reduces a product vertex by vertex: $e_a(x_1,\ldots,x_m)e_b(x_1,\ldots,x_m)$ is expanded by the Pieri rule into two-column Schur polynomials, and each two-column shape is removed by a quantum rim-hook relation,
$$s_{(2^c,1^d)}[m]=e_{c+d}[m]\,e_c[m{+}1]-e_{c-1}[m]\,e_{c+d+1}[m{+}1]-q_m\,e_{c+d}[m]\,e_{c-2}[m{-}1]+q_m\,e_{c-1}[m]\,e_{c+d-1}[m{-}1]$$
(where $e_k[m]=e_k(x_1,\ldots,x_m)$ and $e_k[0]=e_k[n]=\delta_{k0}$), which pushes a factor to a neighboring vertex, after which the process repeats. Since the quantization map is multiplicative on standard elementary monomials, Schubert classes are converted into and out of this basis by the classical transition matrix, of size $n!\times n!$.

We implemented this algorithm, memoizing the straightening, and verified it against \texttt{schubmult} on all $576$ products in $S_4$ and $150$ random products in $S_5$. Table \ref{table:benchmark} compares running times on the same random pairs $u,v\in S_n$ (non-equivariant, $y=z=0$, the only case the rim-hook method covers), excluding the one-time construction of the transition matrix. For $n\leq6$ the two pure-Python implementations are comparable. At $n=7$ the straightening is already five to seven times slower than the pure-Python kernel of \texttt{schubmult} and more than $500$ times slower than its compiled kernel, and the cost of a single straightened product of two basis monomials grows by a factor of roughly $40$ to $50$ with each increment of $n$ ($18$\,ms, $0.69$\,s and $37$\,s on average at $n=6,7,8$); at $n=8$ a single Schubert product, which requires hundreds of such monomial products together with a $40320\times40320$ transition matrix, is out of reach, while \texttt{schubmult} takes $18$\,ms with the compiled kernel and $1.5$\,s in pure Python, and handles $n=9$ in about $0.1$\,s. The reason is structural. The algorithm of Section \ref{section:algorithms} is not cancellation-free, since the expansion of Theorem \ref{theorem:quandoublemed} is signed, but it expands only one factor into elementary polynomials, and the Pieri formula then multiplies each product of elementary polynomials against the other factor directly in the Schubert basis, so that the cancellation is confined to that one signed expansion. Straightening expands both factors into a basis in which Schubert classes are dense, and every reduction step creates further terms that later cancel. The triple constants, with $y\neq z$, are not accessible to the rim-hook method at all.

\begin{table}[ht]
\begin{center}
\begin{tabular}{c|c|c|c|c}
$n$ & products & \texttt{schubmult}, compiled & \texttt{schubmult}, pure Python & quantum rim-hook, pure Python\\
\hline
$5$ & $20$ & $0.1$\,ms & $2.5$\,ms & $1.4$\,ms\\
$6$ & $20$ & $0.7$\,ms & $38$\,ms & $46$\,ms\\
$7$ & $10$ & $5$\,ms & $0.53$\,s & $2.7$\,s\\
$8$ & $5$ & $18$\,ms & $1.5$\,s & --- \\
$9$ & $3$ & $0.12$\,s & $16$\,s & --- \\
\end{tabular}
\end{center}
\caption{Median wall-clock time per quantum product $\qsch_u\qsch_v$, $u,v\in S_n$ random, computed in the Schubert basis by \texttt{schubmult} (Pieri-based) and by the quantum rim-hook straightening of \cite{gu2020rimhook}; see Section \ref{section:benchmark}.} \label{table:benchmark}
\end{table}

\section{Quantum-to-classical reduction} \label{section:quantclass}

In this section we extend the quantum-to-classical reduction of Huang and Li \cite{huang2015equivariant} from the equivariant quantum structure constants $c_{u,v}^{w,d}(y;y)$ to the triple constants. The quantum-to-classical principle originates with Buch, Kresch and Tamvakis \cite{buch2003gwgrass}, who reduce the quantum product of Schubert classes of a Grassmannian to a classical product on a two-step flag variety; its equivariant (and $K$-theoretic) version is due to Buch and Mihalcea \cite{buch2011qkgrass}, and in combination with Buch's puzzle rule for two-step flag varieties \cite{buch2015mutations} it yields a positive combinatorial rule for $QH_T^*(\mathrm{Gr}(k,n))$ (Remark \ref{remark:grassmannian}). For the complete flag variety, Leung and Li \cite{leung2012classical} showed that Gromov--Witten invariants of different degrees coincide and lower the degree $d$ of a structure constant by one simple coroot at a time, and Huang and Li \cite{huang2015equivariant} extended this to equivariant quantum cohomology, which is the form we use. We write $\mathrm{sgn}_i(w)=1$ if $w(i)>w(i+1)$ and $\mathrm{sgn}_i(w)=0$ otherwise, and we use the notation $\delta_i=\alpha_i^\vee$, $|d|=\sum_id_i$ and $\langle d\rangle_i=2d_i-d_{i-1}-d_{i+1}$ (with $d_0=0$) of Definition \ref{def:pwlift}. We say that $u\in S_\infty$ is \emph{$k$-Grassmannian} if $k$ is its only descent.

\begin{theorem}[{\cite[Theorem~2.7]{huang2015equivariant}}] \label{theorem:quantclass}
	We have the following.
	\begin{enumerate}
	\item $c_{u,v}^{w,d}(y;y)=0$ unless $\mathrm{sgn}_i(w)+\langle d\rangle_i\leq \mathrm{sgn}_i(u)+\mathrm{sgn}_i(v)$ for all $i$.
	\item If $\mathrm{sgn}_k(w)+\langle d\rangle_k=\mathrm{sgn}_k(u)+\mathrm{sgn}_k(v)=2$ for some $k$, then
	$$c_{u,v}^{w,d}(y;y)=c_{us_k,vs_k}^{w,d-\delta_k}(y;y)=\begin{cases}
	c_{us_k,v}^{ws_k,d-\delta_k}(y;y) & \text{if $\mathrm{sgn}_k(w)=0$}\\
	c_{us_k,v}^{ws_k,d}(y;y) & \text{if $\mathrm{sgn}_k(w)=1$}
	\end{cases}$$
	\end{enumerate}
	This is symmetric in $u$ and $v$.
\end{theorem}

\begin{corollary}[{\cite[Corollary~2.8]{huang2015equivariant}}] \label{corollary:quantclass28}
	Let $u,v,w\in S_\infty$ and $d\in Q^\vee$. Fix an index $i$.
	\begin{enumerate}
	\item If $\langle d\rangle_i=\mathrm{sgn}_i(u)=0$ and $\mathrm{sgn}_i(w)=\mathrm{sgn}_i(v)=1$, then
	$$c_{u,v}^{w,d}(y;y)=c_{u,vs_i}^{ws_i,d}(y;y).$$
	\item If $\langle d\rangle_i=\mathrm{sgn}_i(u)=1$ and $\mathrm{sgn}_i(w)=\mathrm{sgn}_i(v)=0$, then
	$$c_{u,v}^{w,d}(y;y)=c_{us_i,v}^{ws_i,d-\delta_i}(y;y).$$
	\end{enumerate}
	This is symmetric in $u$ and $v$.
\end{corollary}

\begin{theorem}[Quantum reduction holds for triple quantum Schubert coefficients]\label{theorem:quantclassz}
We have the following.
	\begin{enumerate}
	\item $c_{u,v}^{w,d}(y;z)=0$ unless $\mathrm{sgn}_i(w)+\langle d\rangle_i\leq \mathrm{sgn}_i(u)+\mathrm{sgn}_i(v)$ for all $i$.
	\item If $\mathrm{sgn}_k(w)+\langle d\rangle_k=\mathrm{sgn}_k(u)+\mathrm{sgn}_k(v)=2$ for some $k$, then
	$$c_{u,v}^{w,d}(y;z)=c_{us_k,vs_k}^{w,d-\delta_k}(y;z)=\begin{cases}
	c_{us_k,v}^{ws_k,d-\delta_k}(y;z) & \text{if $\mathrm{sgn}_k(w)=0$}\\
	c_{us_k,v}^{ws_k,d}(y;z) & \text{if $\mathrm{sgn}_k(w)=1$}
	\end{cases}$$
	\end{enumerate}
	This is symmetric in $u$ and $v$.
\end{theorem}
\begin{proof}
For (1), we note that since $\mathrm{sgn}_i(v_2)\leq \mathrm{sgn}_i(v)$ for any factorization $v_1\cdot v_2=v$, we have that $\mathrm{sgn}_i(w)+\langle d\rangle_i\leq \mathrm{sgn}_i(u)+\mathrm{sgn}_i(v_2)$ is a necessary condition for $c_{u,v_2}^{w,d}(y;y)\neq 0$. The result then follows from Lemma \ref{lemma:quantumtoordinary}.

For (2), we note that if $\mathrm{sgn}_k(w)+\langle d\rangle_k=\mathrm{sgn}_k(u)+\mathrm{sgn}_k(v)=2$, then since $\mathrm{sgn}_k(v_2)=\mathrm{sgn}_k(v)$ for any factorization $v_1\cdot v_2=v$ such that $c_{u,v_2}^{w,d}(y;y)\neq 0$, Theorem \ref{theorem:quantclass}(2) holds for each term in the sum in Lemma \ref{lemma:quantumtoordinary}, and the result follows.
\end{proof}

Leung and Li \cite[Theorem~1.2]{leung2012classical} observed that when one of the two factors is Grassmannian, Theorem \ref{theorem:quantclass} can be iterated until no quantum parameter remains. Since their argument uses nothing but the two parts of that theorem, it applies verbatim to the triple constants.

\begin{proposition}[Grassmannian factor: reduction to classical triple constants] \label{prop:grassmannianreduction}
Let $v$ be $k$-Grassmannian, let $u,w\in S_\infty$ and let $d\neq0$ be effective. Then either $c_{u,v}^{w,d}(y;z)=0$, or there is an index $m$ with $\langle d\rangle_m>0$, $\mathrm{sgn}_m(u)=1$ and $\mathrm{sgn}_m(w)=0$ such that
$$c_{u,v}^{w,d}(y;z)=c_{us_m,v}^{ws_m,\,d-\delta_m}(y;z).$$
Consequently $c_{u,v}^{w,d}(y;z)$ is either $0$ or equal to a classical triple constant $c_{u',v}^{w',0}(y;z)$ with $u'=us_{m_1}\cdots s_{m_{|d|}}$ and $w'=ws_{m_1}\cdots s_{m_{|d|}}$, for a sequence of indices determined by $(u,w,d)$.
\end{proposition}
\begin{proof}
Let $I=\{m:\langle d\rangle_m>0\}$. If $I$ were empty, $d$ would be a nonpositive combination of the fundamental coweights, hence of the simple coroots, contradicting $d\neq0$ effective; so $I\neq\emptyset$. Suppose first that $|I|\geq2$ and choose $m\in I$ with $m\neq k$, so that $\mathrm{sgn}_m(v)=0$. By Theorem \ref{theorem:quantclassz}(1), $c_{u,v}^{w,d}(y;z)=0$ unless $\mathrm{sgn}_m(w)+\langle d\rangle_m\leq\mathrm{sgn}_m(u)\leq1$, which forces $\mathrm{sgn}_m(u)=1$, $\mathrm{sgn}_m(w)=0$ and $\langle d\rangle_m=1$. Then the quadruple $(u,vs_m,ws_m,d)$ satisfies $\mathrm{sgn}_m(ws_m)+\langle d\rangle_m=2=\mathrm{sgn}_m(u)+\mathrm{sgn}_m(vs_m)$, and Theorem \ref{theorem:quantclassz}(2) applied to it gives, on the one hand, $c_{u,vs_m}^{ws_m,d}=c_{us_m,v}^{ws_m,d-\delta_m}$, and on the other hand, since $\mathrm{sgn}_m(ws_m)=1$, $c_{u,vs_m}^{ws_m,d}=c_{u,v}^{w,d}$ by the second equality with the roles of $u$ and $v$ exchanged. Suppose next that $I=\{m\}$. Then $\langle d\rangle_m\geq2$ by \cite[Lemma~2.10]{leung2012classical}, and Theorem \ref{theorem:quantclassz}(1) forces $m=k$, $\mathrm{sgn}_k(u)=1$, $\mathrm{sgn}_k(w)=0$ and $\langle d\rangle_k=2$, whereupon Theorem \ref{theorem:quantclassz}(2) gives $c_{u,v}^{w,d}=c_{us_k,v}^{ws_k,d-\delta_k}$. In either case $|d|$ drops by one, and the last assertion follows by induction on $|d|$.
\end{proof}

\begin{remark}
Proposition \ref{prop:grassmannianreduction} with $v=c[p,k]$ expresses every coefficient of a Pieri product as a classical triple Pieri coefficient, which \cite[Theorem~7.1]{samuelmolev} describes positively. At $z=y$ this is exactly the route by which Huang and Li obtain their equivariant quantum Pieri rule \cite[Theorem~3.10]{huang2015equivariant}, with Robinson's rule \cite{robinsonpieri} in the role of the classical formula; the same route with the triple formula in that role yields a triple quantum Pieri rule of their form, in which the terms are indexed by a degree $d$, the twists $(u',w')$ of $(u,w)$ that $d$ determines, and classical chains from $u'$ to $w'$. Theorem \ref{theorem:pieri} is a statement of a different form: its terms are indexed by chains in the quantum Bruhat graph from $u$ itself, and the degree and the weight are read off from the chain. We have checked numerically, for $v=d_2[2]$ and all $u\in S_4$, that the recursion of the proposition carries each of the $125$ coefficients of positive degree in $q$ of $\qsch_u(x;y)\qsch_v(x;z)$ to the predicted classical triple constant.
\end{remark}

\section{Proof of the Pieri formula} \label{section:proofpieri}

\newcommand{\ntom}[1]{\centernot{\tom{#1}}}
\newcommand{\desc}{\mathrm{desc}}
\newcommand{\lup}{\triangleleft}
\newcommand{\supp}{\mathrm{supp}}
\newcommand{\Lw}{L}
\newcommand{\qw}{\mathbf{q}}

Throughout this section permutations act on positions: $(vt_{ab})(a)=v(b)$, $(vt_{ab})(b)=v(a)$, and $(vt_{ab})(r)=v(r)$ otherwise. A cycle is written $(a_p,\ldots,a_1,b)$ with the meaning $a_p\mapsto a_{p-1}\mapsto\cdots\mapsto a_1\mapsto b\mapsto a_p$; when $b=\max\{a_1,\ldots,a_p,b\}$ we call $b$ the \emph{top} of the cycle and refer to
\begin{equation}\label{eq:pivot}
(a_p,\ldots,a_1,b)=t_{a_1b}\,t_{a_2b}\cdots t_{a_pb}
\end{equation}
as its \emph{pivot factorization}. (Indeed, the right-hand side sends $a_p\mapsto b\mapsto a_{p-1}$, $b\mapsto a_p$ and $a_1\mapsto b$.)

\subsection{The single-variable Monk formula}

\begin{lemma} \label{lemma:singlevariable}
Let $u\in S_\infty$ and let $i,j\geq 1$. Then
$$(x_i-z_j)\qsch_u(x;y)=(y_{u(i)}-z_j)\qsch_u(x;y)+\sum_{\substack{p\neq i\\\ell(ut_{ip})=\ell(u)+1}}\mathrm{sign}(p-i)\qsch_{ut_{ip}}(x;y)$$
$$+\sum_{\substack{p\neq i\\\ell(ut_{ip})=\ell(u)-2|p-i|+1}}\mathrm{sign}(p-i)q_{ip}\qsch_{ut_{ip}}(x;y)$$
where $q_{ip}$ stands for $q_{\min(i,p)\max(i,p)}$.
\end{lemma}
\begin{proof}
The quantum double Schubert polynomials represent the Schubert classes in the torus-equivariant quantum cohomology of the flag manifold \cite{lam2014quantum}, so they satisfy the equivariant quantum Chevalley formula of Mihalcea \cite[Theorem~1]{mihalcea2007equivariant} (see \cite{buch2015curve} for a proof via curve neighborhoods; at $y=0$ it is the quantum Monk formula of \cite{fomin1997quantum}, for which Buch \cite{buch2003qmonk} gave a direct geometric proof). In type $A$, with $\qsch_{s_m}(x;y)=\sum_{a\leq m}(x_a-y_a)$, that formula reads
$$\qsch_{s_m}(x;y)\qsch_u(x;y)=\sum_{a\leq m}(y_{u(a)}-y_a)\,\qsch_u(x;y)+\sum_{\substack{a\leq m<b\\ \ell(ut_{ab})=\ell(u)+1}}\qsch_{ut_{ab}}(x;y)+\sum_{\substack{a\leq m<b\\ \ell(ut_{ab})=\ell(u)-2(b-a)+1}}q_{ab}\qsch_{ut_{ab}}(x;y),$$
the diagonal coefficient being the localization of $\qsch_{s_m}(x;y)$ at $u$, and the two sums ranging over the classical and quantum covers of $u$ in the quantum Bruhat graph (for $\alpha=e_a-e_b$ one has $\langle\omega_m,\alpha^\vee\rangle=1$ exactly when $a\leq m<b$, and $\langle 2\rho,\alpha^\vee\rangle=2(b-a)$). By stability we may work in $S_n$ for $n$ large. Now
$$x_i-z_j=\qsch_{s_i}(x;y)-\qsch_{s_{i-1}}(x;y)+y_i-z_j,$$
with the term $\qsch_{s_{i-1}}$ omitted when $i=1$. Subtracting the two Chevalley formulas, the diagonal coefficients combine to $(y_{u(i)}-y_i)+(y_i-z_j)=y_{u(i)}-z_j$, and the reflections $t_{ab}$ with $a<i<b$ occur in both sums with the same coefficient and cancel. What remains are the reflections $t_{ib}$ with $b>i$ (from $\qsch_{s_i}$, with sign $+1$) and $t_{ai}$ with $a<i$ (from $\qsch_{s_{i-1}}$, with sign $-1$), each carrying the factor $1$ or $q_{ab}$ according to whether the cover is classical or quantum. This is the asserted formula.
\end{proof}

\subsection{Covers and tight cycles}

The following classical count is the only fact about lengths that we shall use.

\begin{lemma} \label{lemma:lengthchange}
Let $v\in S_\infty$ and $a<b$. If $v(a)<v(b)$ then
$$\ell(vt_{ab})=\ell(v)+1+2\,\#\{r\mid a<r<b,\ v(a)<v(r)<v(b)\},$$
and if $v(a)>v(b)$ then
$$\ell(vt_{ab})=\ell(v)-1-2\,\#\{r\mid a<r<b,\ v(b)<v(r)<v(a)\}.$$
In particular $\ell(vt_{ab})\geq\ell(v)+1$ if $v(a)<v(b)$, and $\ell(vt_{ab})\geq\ell(v)-(2(b-a)-1)$ always.
\end{lemma}
\begin{proof}
Passing from $v$ to $vt_{ab}$ changes the inversion status of the pair of positions $(a,b)$, and of the pairs $(a,r)$, $(r,b)$ with $a<r<b$ exactly when $v(r)$ lies strictly between $v(a)$ and $v(b)$ (in which case both pairs change); all other pairs are unaffected.
\end{proof}

\begin{definition}[Covers, $\lup$, $\kappa$] \label{def:cover}
For $v\in S_\infty$ and $a<b$ we write $v\lup vt_{ab}$, and say that $t_{ab}$ is a \emph{cover} of $v$, if either $\ell(vt_{ab})=\ell(v)+1$ (a \emph{classical} cover) or $\ell(vt_{ab})=\ell(v)-2(b-a)+1$ (a \emph{quantum} cover). We set $\kappa(v,t_{ab})=1$ in the first case and $\kappa(v,t_{ab})=q_{ab}$ in the second. By Lemma \ref{lemma:lengthchange}, $t_{ab}$ is a classical cover of $v$ if and only if $v(a)<v(b)$ and no $r$ with $a<r<b$ satisfies $v(a)<v(r)<v(b)$, and a quantum cover if and only if $v(a)>v(b)$ and every $r$ with $a<r<b$ satisfies $v(b)<v(r)<v(a)$. These are the covers of the quantum Bruhat graph: the classical covers are the terms of Monk's formula and the quantum covers the terms carrying $q_{ab}$ in the quantum Chevalley formula \cite{fulton2004quantum,buch2003qmonk}; condition (4) of Definition \ref{definition:pieri} says precisely that each step of a witnessing sequence is a cover (Remark \ref{remark:qkbruhat}).
\end{definition}

It is convenient to separate the two values being exchanged from the values lying between them. For distinct numbers $\alpha\neq\beta$ and a finite set $J$ of numbers not containing $\alpha,\beta$, say that the pair $(\alpha,\beta)$ is \emph{$J$-admissible} if either $\alpha<\beta$ and $J\cap(\alpha,\beta)=\emptyset$, or $\alpha>\beta$ and $J\subseteq(\beta,\alpha)$. Thus $t_{ab}$ is a cover of $v$ if and only if $(v(a),v(b))$ is $J$-admissible for $J=\{v(r)\mid a<r<b\}$; the cover is classical when $\alpha<\beta$ and quantum when $\alpha>\beta$.

\begin{lemma} \label{lemma:admissible}
Let $\alpha,\beta,\gamma$ be distinct numbers and $I,J,J'$ finite sets of numbers avoiding them.
\begin{enumerate}
\item[(i)] $(\alpha,\beta)$ is $(J\cup J')$-admissible if and only if it is $J$-admissible and $J'$-admissible. In particular $J'\subseteq J$ and $J$-admissibility imply $J'$-admissibility.
\item[(ii)] (Composition.) If $(\alpha_1,\alpha_0),(\alpha_2,\alpha_1),\ldots,(\alpha_n,\alpha_{n-1})$ are $I$-admissible then so is $(\alpha_n,\alpha_0)$.
\item[(iii)] (Triangle.) If $(\gamma,\beta)$ is $I$-admissible and $(\alpha,\gamma)$ is $(J\cup\{\beta\}\cup I)$-admissible, then $(\alpha,\beta)$ is $I$-admissible and $(\alpha,\gamma)$ is $J$-admissible.
\end{enumerate}
\end{lemma}
\begin{proof}
(i) is immediate from the definition. For (ii), if $I=\emptyset$ there is nothing to prove. Otherwise $I$ cuts the complement of $I$ into intervals (``blocks''), among them $V_-=(-\infty,\min I)$ and $V_+=(\max I,\infty)$. An $I$-admissible pair $(\alpha',\alpha)$ with $\alpha'<\alpha$ has $\alpha,\alpha'$ in the same block, and one with $\alpha'>\alpha$ has $\alpha\in V_-$ and $\alpha'\in V_+$. Hence along the chain the values descend within a block, except for at most one jump from $V_-$ to $V_+$, after which they remain in $V_+$. Either $\alpha_n<\alpha_0$ lie in a common block, or $\alpha_0\in V_-$ and $\alpha_n\in V_+$; in both cases $(\alpha_n,\alpha_0)$ is $I$-admissible. For (iii), the second assertion is (i). For the first there are four cases. If $\gamma<\beta$ and $\alpha<\gamma$: then $\alpha<\beta$ and $I\cap(\alpha,\beta)\subseteq (I\cap(\alpha,\gamma))\cup(I\cap[\gamma,\beta))=\emptyset$. If $\gamma<\beta$ and $\alpha>\gamma$: then $\beta\in(\gamma,\alpha)$ gives $\beta<\alpha$, and $I\subseteq(\gamma,\alpha)$ together with $I\cap(\gamma,\beta)=\emptyset$ gives $I\subseteq(\beta,\alpha)$. If $\gamma>\beta$ and $\alpha<\gamma$: then $\beta\notin(\alpha,\gamma)$ forces $\beta<\alpha$, and $I\subseteq(\beta,\gamma)$ with $I\cap(\alpha,\gamma)=\emptyset$ gives $I\subseteq(\beta,\alpha)$. The case $\gamma>\beta$, $\alpha>\gamma$ cannot occur, since it would require $\beta\in(\gamma,\alpha)$.
\end{proof}

\begin{definition}[Descents of a cycle, $\desc_u(c)$]
Let $c$ be a cycle and suppose $u\in S_\infty$. Write $c$ as
$$(a_1,\ldots,a_n)$$
where $a_n>a_i$ for all $i<n$ (this is a unique way to express the cycle since the $a_j$ are distinct). Then for any $i$ with $1\leq i<n$ we say that $a_i$ is a \emph{descent} of $c$ with respect to $u$ if $u(a_i)>u(a_{i+1})$. We denote the set of all descents of $c$ with respect to $u$ by $\desc_u(c)$.
\end{definition}

In the notation $c=(a_p,\ldots,a_1,b)$ of \eqref{eq:pivot}, with $a_0:=b$, we have $a_j\in\desc_u(c)$ if and only if $u(a_j)>u(a_{j-1})$; the top $b$ is never a descent.

\begin{definition}[$\Lw_u(\sigma)$, $\qw_u(\sigma)$, $m$-shaped, $u$-tight] \label{def:tight}
Let $u,\sigma\in S_\infty$ and let $c_1,\ldots,c_n$ be the nontrivial cycles of $\sigma$, with $c_i=(a_{i,p_i},\ldots,a_{i,1},b_i)$ and $b_i=\max c_i$. Put
$$\Lw_u(\sigma)=\sum_{i=1}^n\Big(p_i-2\sum_{a\in\desc_u(c_i)}(b_i-a)\Big),\qquad \qw_u(\sigma)=\prod_{i=1}^n\prod_{a\in\desc_u(c_i)}q_{ab_i}.$$
Both depend only on the restriction of $u$ to $\supp(\sigma)$, and $\deg\qw_u(\sigma)=\sum_ip_i-\Lw_u(\sigma)$. We say that $\sigma$ is \emph{$m$-shaped} if every nontrivial cycle of $\sigma$ contains exactly one element greater than $m$, and that $\sigma$ is \emph{$u$-tight} if
$$\ell(u\sigma)-\ell(u)=\Lw_u(\sigma).$$
A \emph{pivot chain} of $\sigma$ is a sequence $t_1,\ldots,t_N$ ($N=\sum_ip_i$) of transpositions obtained by interleaving the pivot factorizations \eqref{eq:pivot} of $c_1,\ldots,c_n$ in any manner that preserves the order of the factors within each cycle; since the cycles are disjoint, $t_1\cdots t_N=\sigma$.
\end{definition}

\begin{lemma}[Tightness] \label{lemma:tight}
Let $u,\sigma\in S_\infty$, let $t_1,\ldots,t_N$ be a pivot chain of $\sigma$, and put $u_0=u$, $u_j=u_{j-1}t_j$.
\begin{enumerate}
\item[(i)] If $t_j=t_{a_{i,l}b_i}$ then $u_{j-1}(a_{i,l})=u(a_{i,l})$ and $u_{j-1}(b_i)=u(a_{i,l-1})$ (with $a_{i,0}=b_i$). Thus the $j$-th step exchanges the values $u(a_{i,l})$ and $u(a_{i,l-1})$, and $a_{i,l}\in\desc_u(c_i)$ exactly when the larger of the two sits at the smaller position $a_{i,l}$ before the step.
\item[(ii)] $\ell(u\sigma)-\ell(u)\geq\Lw_u(\sigma)$.
\item[(iii)] $\sigma$ is $u$-tight if and only if $u_{j-1}\lup u_j$ for every $j$. In that case the $j$-th step is a classical cover if $a_{i,l}\notin\desc_u(c_i)$ and a quantum cover if $a_{i,l}\in\desc_u(c_i)$, and $\prod_j\kappa(u_{j-1},t_j)=\qw_u(\sigma)$.
\item[(iv)] If $\sigma=\rho\tau$ where $\rho$ and $\tau$ are the products of two complementary sets of cycles of $\sigma$, then $\sigma$ is $u$-tight if and only if $\rho$ is $u$-tight and $\tau$ is $u\rho$-tight. Moreover $\Lw_u(\tau)=\Lw_{u\rho}(\tau)$ and $\qw_u(\tau)=\qw_{u\rho}(\tau)$.
\item[(v)] If $\sigma$ is $u$-tight then so is $t_1\cdots t_j$ for every $j$, and $u_j$ is obtained from $u$ by a sequence of covers.
\end{enumerate}
\end{lemma}
\begin{proof}
(i) The position $a_{i,l}$ is moved only by the factor $t_{a_{i,l}b_i}$, so before that factor is applied it carries $u(a_{i,l})$; the position $b_i$ is moved only by the factors of $c_i$, and after $t_{a_{i,l-1}b_i}$ it carries $u(a_{i,l-1})$ (and $u(b_i)$ before any factor of $c_i$). For (ii) and (iii), by (i) and Lemma \ref{lemma:lengthchange} the $j$-th step increases the length by at least $1$ if $a_{i,l}\notin\desc_u(c_i)$ and by at least $-(2(b_i-a_{i,l})-1)$ if $a_{i,l}\in\desc_u(c_i)$; in both cases equality holds exactly when $t_j$ is a cover of $u_{j-1}$, classical in the first case and quantum in the second (a step with $u_{j-1}(a_{i,l})<u_{j-1}(b_i)$ increases the length and so cannot be a quantum cover, and symmetrically). Summing over $j$ gives (ii), and shows that equality in (ii) holds exactly when every step is a cover; the product of the $\kappa$'s is then the product of $q_{a_{i,l}b_i}$ over the descents, which is $\qw_u(\sigma)$. For (iv), $\ell(u\rho\tau)-\ell(u)=(\ell(u\rho)-\ell(u))+(\ell(u\rho\tau)-\ell(u\rho))\geq \Lw_u(\rho)+\Lw_{u\rho}(\tau)=\Lw_u(\rho)+\Lw_u(\tau)=\Lw_u(\sigma)$, where $\Lw_{u\rho}(\tau)=\Lw_u(\tau)$ because $u\rho$ and $u$ agree on $\supp(\tau)$; equality holds if and only if both inequalities are equalities. (v) The partial product $t_1\cdots t_j$ is again a product of disjoint cycles of the form $(a_{i,l},\ldots,a_{i,1},b_i)$ (a cycle of $\sigma$ truncated after its first $l$ factors), and $t_1,\ldots,t_j$ is a pivot chain of it, all of whose steps are covers by (iii).
\end{proof}

With this language the cycle description of $\tom{k}$ becomes a statement in which the chain disappears entirely.

\begin{lemma} \label{lemma:piericycle}
Suppose $u,w\in S_\infty$ and let $k\geq0$ be an integer, and put $\pi=u^{-1}w$. Then $u\tom{k} w$ if and only if $\pi$ is $k$-shaped and $u$-tight. In that case the pivot chain of $\pi$ in which the cycles are taken in order of decreasing tops is a sequence $t_{a_1b_1},\ldots,t_{a_pb_p}$ as in Definition \ref{definition:pieri}, every pivot chain of $\pi$ is a chain of covers, and
$$q_k(u,w)=\qw_u(\pi).$$
In particular $q_k(u,w)$ does not depend on the witnessing sequence, and $\ell(u,w)+\deg q_k(u,w)$ equals the sum of the lengths of the cycles of $\pi$, each diminished by one.
\end{lemma}
\begin{proof}
Suppose first that $t_{a_1b_1},\ldots,t_{a_pb_p}$ satisfy conditions (1)--(5) of Definition \ref{definition:pieri}. The $b_i$ are weakly decreasing, so the factors sharing a common value $b$ of $b_i$ are consecutive, and by \eqref{eq:pivot} their product $t_{a_ib}\cdots t_{a_{i+r}b}$ is the cycle $(a_{i+r},\ldots,a_i,b)$, whose top is $b$ because $a_l\leq k<b$. Cycles belonging to distinct values of $b$ are disjoint, since the $a_i$ are pairwise distinct and no $a_i$ equals any $b_j$. Hence $\pi$ is the product of these cycles, each of which contains exactly one element greater than $k$, namely its top: $\pi$ is $k$-shaped. The given sequence is a pivot chain of $\pi$, and condition (4) says that each of its steps is a cover; so $\pi$ is $u$-tight by Lemma \ref{lemma:tight}(iii), and the recursive definition of $q_k(u,w)$ assigns to each step exactly the factor $\kappa$, whence $q_k(u,w)=\qw_u(\pi)$.

Conversely, suppose $\pi$ is $k$-shaped and $u$-tight. Order the nontrivial cycles of $\pi$ by decreasing tops and take the corresponding pivot chain. The top of each cycle is its unique element exceeding $k$, so condition (1) holds; (2) holds because the cycles are disjoint and the $a$'s within a cycle are distinct; (3) holds by the chosen order; (4) holds by Lemma \ref{lemma:tight}(iii); and (5) is $t_1\cdots t_N=\pi$. The final assertion is $\deg\qw_u(\pi)=\sum_ip_i-\Lw_u(\pi)$.
\end{proof}

\begin{corollary}[Cycle formula for the quantum weight] \label{corollary:cycleweight}
Suppose $u\tom{k}w$, and write the nontrivial cycles in the decomposition of
$u^{-1}w$ as $c_i=(a_{i,p_i},\ldots,a_{i,1},b_i)$, as in Lemma
\ref{lemma:piericycle}. Then
\[
q_k(u,w)=\prod_i\ \prod_{a\in\desc_u(c_i)}q_{a b_i}.
\]
Equivalently, the exponent of $q_r$ in $q_k(u,w)$ is
\[
\#\{(i,a):a\in\desc_u(c_i),\ a\leq r<b_i\}.
\]
Moreover, $q_k(u,w)$ is multiplicative over the cycles of $u^{-1}w$, and the factor attached to a cycle depends only on the values of $u$ on that cycle.
\end{corollary}
\begin{proof}
This is $q_k(u,w)=\qw_u(\pi)$ together with $q_{ab}=q_a\cdots q_{b-1}$.
\end{proof}

We next show that a tight cycle is very constrained: the values of $u$ along it, read cyclically, decrease except at exactly one place. For a cycle $c=(a_p,\ldots,a_1,b)$ with top $b$ and a permutation $u$ we write
$$\beta_0=u(b),\qquad \beta_j=u(a_j)\quad(1\leq j\leq p),$$
so that $a_j\in\desc_u(c)$ if and only if $\beta_j>\beta_{j-1}$, and we put $u_j=u\,t_{a_1b}\cdots t_{a_jb}$. By Lemma \ref{lemma:tight}(i), the $j$-th step of the pivot chain exchanges $\beta_j$ (at $a_j$) with $\beta_{j-1}$ (at $b$); moreover $u_m(a_l)=\beta_{l-1}$ for $l\leq m$ and $u_m(a_l)=\beta_l$ for $l>m$. If $c$ is $u$-tight, Definition \ref{def:cover} gives, for every $j$ and every position $r$ with $a_j<r<b$:
\begin{align}
&\text{if }\beta_j<\beta_{j-1}\text{ (classical step): } u_{j-1}(r)\notin(\beta_j,\beta_{j-1}); \label{eq:Cj}\\
&\text{if }\beta_j>\beta_{j-1}\text{ (quantum step): } u_{j-1}(r)\in(\beta_{j-1},\beta_j). \label{eq:Qj}
\end{align}

\begin{lemma} \label{lemma:onedescent}
Let $c=(a_p,\ldots,a_1,b)$ be a cycle with top $b$ which is $u$-tight. Then $\desc_u(c)$ has at most one element. In particular this applies to every cycle of $u^{-1}w$ whenever $u\tom{k}w$.
\end{lemma}
\begin{proof}
Suppose $i<j$ are two descents (quantum steps) with no descent strictly between them, so that $\beta_i>\beta_{i-1}$, $\beta_j>\beta_{j-1}$ and $\beta_i>\beta_{i+1}>\cdots>\beta_{j-1}$.

If $a_j>a_i$, then \eqref{eq:Qj} for the step $i$ and $r=a_j$ gives $\beta_j=u_{i-1}(a_j)\in(\beta_{i-1},\beta_i)$, so $\beta_{j-1}<\beta_j<\beta_i$. Hence $j\geq i+2$, and there is $m$ with $i<m<j$ and $\beta_m<\beta_j<\beta_{m-1}$. The step $m$ is classical. If $a_j>a_m$, then \eqref{eq:Cj} for the step $m$ and $r=a_j$ is violated by $u_{m-1}(a_j)=\beta_j$. So $a_m>a_j$, and then \eqref{eq:Qj} for the step $j$ and $r=a_m$ requires $u_{j-1}(a_m)=\beta_{m-1}\in(\beta_{j-1},\beta_j)$, contradicting $\beta_{m-1}>\beta_j$.

If $a_j<a_i$, then \eqref{eq:Qj} for the step $j$ and $r=a_i$ gives $u_{j-1}(a_i)=\beta_{i-1}\in(\beta_{j-1},\beta_j)$, so $\beta_{j-1}<\beta_{i-1}<\beta_i$. Hence again $j\geq i+2$ and there is $m$ with $i<m<j$ and $\beta_m<\beta_{i-1}<\beta_{m-1}$. If $a_i>a_m$, \eqref{eq:Cj} for the step $m$ and $r=a_i$ is violated by $u_{m-1}(a_i)=\beta_{i-1}$. So $a_m>a_i$, and \eqref{eq:Qj} for the step $i$ and $r=a_m$ requires $u_{i-1}(a_m)=\beta_m\in(\beta_{i-1},\beta_i)$, contradicting $\beta_m<\beta_{i-1}$.
\end{proof}

\begin{definition}[Redefinition of $\desc_u(c)$]
If $c$ is a tight cycle, we redefine $\desc_u(c)$ to be the unique element of the previously defined set $\desc_u(c)$ if it exists, or otherwise $0$.
\end{definition}

\begin{lemma} \label{lemma:geq}
Let $c=(a_p,\ldots,a_1,b)$ be a cycle with top $b$ which is $u$-tight, and suppose $\desc_u(c)=a_q$ with $q\geq1$. Then
$$\beta_q>\beta_{q+1}>\cdots>\beta_p>\beta_0>\beta_1>\cdots>\beta_{q-1}.$$
In particular $u(a_j)>u(b)$ for all $j\geq q$.
\end{lemma}
\begin{proof}
All steps other than the $q$-th are classical, which gives the two decreasing runs; it remains to prove $\beta_l>\beta_0$ for $l\geq q$.

First let $l=q$. For $l'<q$ with $a_{l'}<a_q$, \eqref{eq:Cj} for the step $l'$ and $r=a_q$ gives $\beta_q=u_{l'-1}(a_q)\notin(\beta_{l'},\beta_{l'-1})$. If some $l'<q$ has $a_{l'}>a_q$, let $l^\ast$ be the least such; then \eqref{eq:Qj} for the step $q$ and $r=a_{l^\ast}$ gives $u_{q-1}(a_{l^\ast})=\beta_{l^\ast-1}\in(\beta_{q-1},\beta_q)$, so $\beta_q>\beta_{l^\ast-1}$, and descending induction on $l'<l^\ast$ using $\beta_q\notin(\beta_{l'},\beta_{l'-1})$ yields $\beta_q>\beta_{l'-1}$ for all $l'\leq l^\ast$, in particular $\beta_q>\beta_0$. If no such $l'$ exists, the same induction starting from $\beta_q>\beta_{q-1}$ gives $\beta_q>\beta_0$.

Now let $l>q$ and assume $\beta_{l-1}>\beta_0$ (and $\beta_l<\beta_{l-1}$). Suppose, for a contradiction, that $\beta_l<\beta_0$. We claim that $a_{l'}<a_l$ and $\beta_l<\beta_{l'}$ for $1\leq l'\leq q-1$. Indeed, assume $\beta_l<\beta_{l'-1}$ (true for $l'=1$). If $a_{l'}>a_l$, then \eqref{eq:Cj} for the step $l$ and $r=a_{l'}$ requires $u_{l-1}(a_{l'})=\beta_{l'-1}\notin(\beta_l,\beta_{l-1})$, which is false as $\beta_l<\beta_{l'-1}\leq\beta_0<\beta_{l-1}$. So $a_{l'}<a_l$, and \eqref{eq:Cj} for the step $l'$ and $r=a_l$ gives $\beta_l\notin(\beta_{l'},\beta_{l'-1})$, whence $\beta_l<\beta_{l'}$. The claim yields $\beta_l<\beta_{q-1}$. Finally consider the quantum step $q$: if $a_q<a_l$, \eqref{eq:Qj} with $r=a_l$ requires $\beta_l\in(\beta_{q-1},\beta_q)$, contradicting $\beta_l<\beta_{q-1}$; if $a_q>a_l$, \eqref{eq:Cj} for the step $l$ with $r=a_q$ requires $u_{l-1}(a_q)=\beta_{q-1}\notin(\beta_l,\beta_{l-1})$, contradicting $\beta_l<\beta_{q-1}<\beta_0<\beta_{l-1}$. Hence $\beta_l>\beta_0$, and the induction is complete.
\end{proof}

\begin{corollary}[Cyclic monotonicity] \label{corollary:cyclic}
Let $c=(a_p,\ldots,a_1,b)$ be a $u$-tight cycle with top $b$. Then the cyclic sequence $(\beta_0,\beta_1,\ldots,\beta_p)$ has exactly one cyclic ascent: either $\beta_0>\beta_1>\cdots>\beta_p$ (and $\desc_u(c)=0$), or $\desc_u(c)=a_q$ with $q\geq1$ and $\beta_q>\cdots>\beta_p>\beta_0>\cdots>\beta_{q-1}$. Consequently
$$\Lw_u(c)=\begin{cases}p&\text{if the ascent is at the wrap }(\beta_p,\beta_0),\\ p-2(b-a_q)&\text{if the ascent is }(\beta_{q-1},\beta_q),\end{cases}\qquad \qw_u(c)=\begin{cases}1,\\ q_{a_qb}.\end{cases}$$
\end{corollary}

We shall repeatedly use the following elementary observation about such sequences. Call a cyclic sequence of distinct numbers \emph{unimodal} if it has exactly one cyclic ascent, i.e.\ if it is a cyclic rotation of a decreasing sequence. Every cyclically consecutive block of a unimodal sequence is unimodal, and so is any subsequence. If $\Gamma$ is unimodal and $\eta\notin\Gamma$, there is exactly one \emph{gap} of $\Gamma$ (a pair of cyclically consecutive entries) into which $\eta$ can be inserted so that the result is again unimodal, namely: the descent gap $(\gamma,\gamma')$ with $\gamma>\eta>\gamma'$ if $\min\Gamma<\eta<\max\Gamma$, and the ascent gap $(\min\Gamma,\max\Gamma)$ if $\eta<\min\Gamma$ or $\eta>\max\Gamma$. We call it the \emph{insertion gap} of $\eta$ in $\Gamma$.

\subsection{The inductive step: Monk neighbors at position $k$}

We organize the inductive step of Theorem \ref{theorem:pieri} around the position $k$. The following lemma records two facts about the interaction of $\tom{k}$ and $\tom{k-1}$ at position $k$.

\begin{lemma} \label{lemma:leveldichotomy}
Suppose $u,w\in S_\infty$ and let $k>1$ be an integer.
\begin{enumerate}
\item[(i)] If $u(k)=w(k)$, then $u\tom{k}w$ if and only if $u\tom{k-1}w$, and in that case $q_k(u,w)=q_{k-1}(u,w)$ and $\fx_k(u,w)=\fx_{k-1}(u,w)\cup\{u(k)\}$.
\item[(ii)] If $u(k)\neq w(k)$ and $u\tom{k}w$, then $u\ntom{k-1}w$.
\end{enumerate}
\end{lemma}
\begin{proof}
Let $\pi=u^{-1}w$. If $u(k)=w(k)$ then $k$ is a fixed point of $\pi$, so no nontrivial cycle of $\pi$ contains $k$, and such a cycle has exactly one element greater than $k$ if and only if it has exactly one element greater than $k-1$. Thus $\pi$ is $k$-shaped if and only if it is $(k-1)$-shaped, while tightness and $\qw_u(\pi)$ do not involve $k$ at all; (i) follows from Lemma \ref{lemma:piericycle}. For (ii), the cycle of $\pi$ containing $k$ is nontrivial and, $\pi$ being $k$-shaped, contains an element greater than $k$; so it contains two elements greater than $k-1$ and $\pi$ is not $(k-1)$-shaped.
\end{proof}

For the rest of this section we fix $u,w\in S_\infty$ and $k>1$, and write
$$\pi=u^{-1}w,\qquad C=\text{the cycle of }\pi\text{ containing }k\ (C=\{k\}\text{ if }\pi(k)=k),\qquad B=\max C.$$
For $m\geq0$ we set $c_m(u,w)=q_m(u,w)\wte(y_{\fx_m(u,w)};z)$ if $u\tom{m}w$ and $c_m(u,w)=0$ otherwise. The \emph{Monk summands} at $w$ are indexed by the reflections $t\in\{t_{ak}\mid a<k\}\cup\{t_{kb}\mid b>k\}$: writing $v=wt$ and $\sigma=u^{-1}v=\pi t$, the summand of $t$ is
$$\mathrm{sign}(t)\,\kappa(v,t)\,c_{k-1}(u,v),\qquad \mathrm{sign}(t_{ak})=-1,\ \mathrm{sign}(t_{kb})=+1,$$
where the summand is understood to be $0$ unless $v\lup w$. (These are precisely the terms of \eqref{eq:up} and \eqref{eq:down} below.) The summand of $t$ is therefore nonzero if and only if $v\lup w$ and $u\tom{k-1}v$, that is, $\sigma$ is $(k-1)$-shaped and $u$-tight.

\begin{lemma}[Classification of the Monk neighbors] \label{lemma:neighbors}
Suppose the Monk summand of $t$ is nonzero, and let $P=\fx_{k-1}(u,w)$. Then every nontrivial cycle of $\pi$ not meeting $t$ has exactly one element greater than $k$, and exactly one of the following holds.
\begin{enumerate}
\item[(U)] $t=t_{kB}$ with $B>k$; $C$ contains exactly two elements greater than $k-1$, namely $k$ and $B$; and $\fx_{k-1}(u,v)=P$. Writing $C=(a_p,\ldots,a_1,B)$ with $a_j=k$, the cycles of $\sigma$ are those of $\pi$ other than $C$ together with
$$C'=(a_{j-1},\ldots,a_1,B)\quad\text{and}\quad C''=(a_p,\ldots,a_{j+1},k).$$
\item[(S)] $t=t_{ak}$ with $a\in C$ and $B>k$; $C$ contains exactly two elements greater than $k-1$, namely $k$ and $B$; in the notation of (U), $a=a_m$ for some $m>j$; $\fx_{k-1}(u,v)=P$; and the cycles of $\sigma$ are those of $\pi$ other than $C$ together with
$$C'_m=(a_p,\ldots,a_{m+1},a_m,a_{j-1},\ldots,a_1,B)\quad\text{and}\quad C''_m=(a_{m-1},\ldots,a_{j+1},k).$$
\item[(P)] $t=t_{ak}$ with $a\in C$ and $B=k$; writing $C=(a_p,\ldots,a_1,k)$ we have $a=a_p=\pi(k)=u^{-1}(w(k))$; $\fx_{k-1}(u,v)=P\sqcup\{w(k)\}$; and the cycles of $\sigma$ are those of $\pi$ other than $C$ together with $C_-=(a_{p-1},\ldots,a_1,k)$ (trivial if $p=1$).
\item[(M)] $t=t_{ak}$ with $a\notin C$; then $B=k$, the cycle $D$ of $\pi$ containing $a$ (possibly $D=\{a\}$) satisfies $D\subseteq[1,k-1]$, and $\fx_{k-1}(u,v)=P\setminus\{u(a)\}$, which equals $P$ unless $D=\{a\}$. The cycles of $\sigma$ are those of $\pi$ other than $C$ and $D$ together with the single cycle $E$ with $\supp(E)=\supp(C)\cup\supp(D)$ obtained by setting $E(a)=\pi(k)$, $E(k)=\pi(a)$ and $E=\pi$ elsewhere.
\end{enumerate}
\end{lemma}
\begin{proof}
A cycle of $\pi$ disjoint from $\supp(t)$ is a cycle of $\sigma$, hence has exactly one element greater than $k-1$, and since it does not contain $k$ this element is greater than $k$. Recall that right multiplication of a permutation by a transposition $t_{xy}$ merges the cycles of $x$ and $y$ if they are distinct, and otherwise splits their common cycle $(x_1,\ldots,x_n)$, $x=x_s$, $y=x_t$, $s<t$, into $(x_t,x_{s+1},\ldots,x_{t-1})$ and $(x_s,x_{t+1},\ldots,x_n,x_1,\ldots,x_{s-1})$.

Let $t=t_{kb}$. If $b\notin C$, the cycle of $\sigma$ through $k$ contains $k$ and $b>k$, which is forbidden. So $b\in C$ and $t$ splits $C$ into a cycle through $k$ and a cycle through $b$; each must contain exactly one element greater than $k-1$, so $C\cap(k-1,\infty)=\{k,b\}$ and $b=B$. The explicit form of the two cycles is the splitting formula applied to $(x_1,\ldots,x_{p+1})=(a_p,\ldots,a_1,B)$ with $x_s=k$, $x_t=B$. Since $v$ agrees with $w$ at all positions less than $k$, $\fx_{k-1}(u,v)=P$. This is (U).

Let $t=t_{ak}$ with $a\in C$, so $t$ splits $C$. If $B>k$, the piece not containing $k$ must contain exactly one element greater than $k-1$, which forces $C\cap(k-1,\infty)=\{k,B\}$ and that piece to contain $B$. With $(x_1,\ldots,x_{p+1})=(a_p,\ldots,a_1,B)$, $k=a_j$ and $a=a_m$: if $m<j$ the piece through $a_m$ is $(a_m,a_{j-1},\ldots,a_{m+1})$, which has no element greater than $k-1$ unless it is trivial ($m=j-1$), in which case the other piece $(k,a_{m-1},\ldots,a_1,B,a_p,\ldots,a_{j+1})$ contains both $k$ and $B$; both are forbidden. So $m>j$, and the splitting formula gives $C'_m$ and $C''_m$. Here $v(a_m)=w(k)=u(a_{j-1})\neq u(a_m)$ and $w(a_m)\neq u(a_m)$, so $a_m$ is a fixed position neither for $v$ nor for $w$, and $\fx_{k-1}(u,v)=P$. This is (S). If $B=k$, the piece not containing $k$ has no element greater than $k-1$ and must be trivial, i.e.\ $\sigma(a)=a$, i.e.\ $\pi(k)=a$; in the notation $C=(a_p,\ldots,a_1,k)$ this is $a=a_p$, and $u(a_p)=u(\pi(k))=w(k)$. Then $\sigma$ has the cycle $C_-$, and $v(a_p)=w(k)=u(a_p)$ makes $a_p$ a new fixed position with value $w(k)\notin P$. This is (P).

Finally let $t=t_{ak}$ with $a\notin C$, so $t$ merges $C$ with the cycle $D$ of $a$ into a cycle $E$ containing $k$; $E$ may contain no further element greater than $k-1$, so $C\subseteq[1,k]$, i.e.\ $B=k$, and $D\subseteq[1,k-1]$. We have $v(a)=w(k)=u(\pi(k))\neq u(a)$ because $\pi(k)\in C\not\ni a$, so $a$ is not a fixed position of $v$; it is one of $w$ exactly when $D=\{a\}$. All other positions less than $k$ agree. This is (M).
\end{proof}

\subsubsection*{The case $B>k$: the up-neighbor and the split down-neighbors}

\begin{lemma}[Up/down pairing] \label{lemma:updown}
Let $\psi\in S_\infty$, let $k<B$, let $z_1,\ldots,z_s$ ($s\geq0$) be distinct positions less than $k$, let $X=(z_s,\ldots,z_1,k,B)$ and $w=\psi X$. Consider the following sequences of transpositions, each with product $X$:
\begin{align*}
F_{\mathrm{std}}&:\ t_{kB},\,t_{z_1B},\,\ldots,\,t_{z_sB} &&\text{(the pivot factorization of $X$)},\\
F_{\mathrm{up}}&:\ t_{z_1k},\,\ldots,\,t_{z_sk},\,t_{kB},\\
F_m&:\ t_{z_1k},\,\ldots,\,t_{z_{m-1}k},\,t_{z_mB},\,\ldots,\,t_{z_sB},\,t_{z_mk} &&(1\leq m\leq s).
\end{align*}
Call such a sequence a \emph{cover path} (from $\psi$) if each of its steps is a cover of the permutation reached before it, and define its \emph{weight} to be the product of the $\kappa$'s of its steps.
\begin{enumerate}
\item[(a)] If $F_{\mathrm{std}}$ is a cover path, then $F_{\mathrm{up}}$ is a cover path of the same weight, and no $F_m$ is a cover path.
\item[(b)] If $F_{\mathrm{up}}$ is a cover path but $F_{\mathrm{std}}$ is not, then exactly one $F_m$ is a cover path, and its weight equals that of $F_{\mathrm{up}}$.
\item[(c)] If some $F_m$ is a cover path, then so is $F_{\mathrm{up}}$.
\end{enumerate}
\end{lemma}
\begin{proof}
(The products are $X$: $F_{\mathrm{up}}$ is $(z_s,\ldots,z_1,k)\,t_{kB}$, and $F_m$ is $(z_{m-1},\ldots,z_1,k)\,(z_s,\ldots,z_m,B)\,t_{z_mk}$; both send $k\mapsto B\mapsto z_s\mapsto\cdots\mapsto z_1\mapsto k$.) If $s=0$ then $F_{\mathrm{std}}=F_{\mathrm{up}}$ and there is nothing to prove, so let $s\geq1$. Put
$$\beta=\psi(B),\qquad \zeta_0=\psi(k),\qquad \zeta_i=\psi(z_i)\ (1\leq i\leq s),\qquad I=\{\psi(r)\mid k<r<B\}.$$
The positions strictly between $k$ and $B$ are never moved, so $I$ is the set of values found there at every stage of every path. Let $\phi_m=\psi t_{z_1k}\cdots t_{z_{m-1}k}$ ($1\leq m\leq s+1$); thus $\phi_m(k)=\zeta_{m-1}$, $\phi_m(z_i)=\zeta_{i-1}$ for $i<m$, $\phi_m(z_i)=\zeta_i$ for $i\geq m$, $\phi_m(B)=\beta$, and $\phi_m$ agrees with $\psi$ elsewhere. Let $J_i=\{\phi_i(r)\mid z_i<r<k\}$. We introduce the predicates
\begin{align*}
U_i&:\ (\zeta_i,\zeta_{i-1})\text{ is }J_i\text{-admissible} &&(1\leq i\leq s),\\
T_i&:\ (\zeta_i,\beta)\text{ is }I\text{-admissible} &&(0\leq i\leq s),\\
A_i&:\ (\zeta_i,\zeta_{i-1})\text{ is }(J_i\cup\{\beta\}\cup I)\text{-admissible} &&(1\leq i\leq s).
\end{align*}
$U_i$ says that the $i$-th step of $F_{\mathrm{up}}$ (from $\phi_i$) is a cover, and $T_s$ that its last step is; so $F_{\mathrm{up}}$ is a cover path if and only if $U_1,\ldots,U_s,T_s$ hold. For $1\leq m\leq s+1$ let $\mathcal S_m$ be the statement that $t_{kB},t_{z_mB},\ldots,t_{z_sB}$ is a cover path from $\phi_m$; thus $F_{\mathrm{std}}$ is a cover path if and only if $\mathcal S_1$, and $\mathcal S_{s+1}=T_s$. The first step of $\mathcal S_m$ exchanges $\zeta_{m-1}$ at $k$ with $\beta$ at $B$ across $I$; the step $t_{z_iB}$ ($i\geq m$) exchanges $\zeta_i$ at $z_i$ with $\zeta_{i-1}$ at $B$ across the positions in $(z_i,k)$, which carry the same values as in $\phi_i$ (the positions $z_l$, $l<i$, carry $\zeta_{l-1}$ and the others their $\psi$-values), the position $k$, which carries $\beta$, and the positions in $(k,B)$. Hence
$$\mathcal S_m\iff T_{m-1}\wedge A_m\wedge\cdots\wedge A_s.$$
Let $\mathcal D_m$ ($1\leq m\leq s$) be the statement that $t_{z_mB},\ldots,t_{z_sB},t_{z_mk}$ is a cover path from $\phi_m$; then $F_m$ is a cover path if and only if $U_1,\ldots,U_{m-1}$ and $\mathcal D_m$ hold. The steps of $\mathcal D_m$ are: $D^{(m)}$, exchanging $\zeta_m$ at $z_m$ with $\beta$ at $B$ across $J_m\cup\{\zeta_{m-1}\}\cup I$ (the position $k$ now carries $\zeta_{m-1}$); for $i>m$, $D^{(i)}$, exchanging $\zeta_i$ with $\zeta_{i-1}$ across $J^\ast_i\cup I$, where $J^\ast_i$ consists of the values at the positions of $(z_i,k]$, so that $J^\ast_i=J_i\cup\{\beta\}$ if $z_m>z_i$ and $J^\ast_i=J_i\cup\{\zeta_{m-1}\}$ if $z_m<z_i$ (the position $z_m$ carries $\beta$ and $k$ carries $\zeta_{m-1}$); and $D^{\mathrm{last}}$, exchanging $\beta$ at $z_m$ with $\zeta_{m-1}$ at $k$ across $J^{\ast\ast}=N_m\cup\{\zeta_{l-1}\mid l>m,\ z_l>z_m\}$, where $N_m$ denotes the set of values carried by the positions of $(z_m,k)$ other than $z_{m+1},\ldots,z_s$ (these values are the same in $\phi_m$ and at the end of $\mathcal D_m$), so that $J_m=N_m\cup\{\zeta_l\mid l>m,\ z_l>z_m\}$. For $i>m$ the set across which $D^{(i)}$ is taken contains $J_i$, and the sets of $D^{(m)},\ldots,D^{(s)}$ all contain $\zeta_{m-1}$, the value at $k$.

\emph{Step 1: $\mathcal S_m\Rightarrow\mathcal S_{m+1}$, and $\mathcal S_m$ implies $U_m,\ldots,U_s,T_s$.} By Lemma \ref{lemma:admissible}(iii) applied with $(\alpha,\gamma,\beta)=(\zeta_m,\zeta_{m-1},\beta)$ and $J=J_m$, the conjunction $T_{m-1}\wedge A_m$ implies $T_m\wedge U_m$. The first claim follows, and the second by iterating it down to $\mathcal S_{s+1}=T_s$, using $A_i\Rightarrow U_i$ (Lemma \ref{lemma:admissible}(i)). In particular $\mathcal S_1$ implies that $F_{\mathrm{up}}$ is a cover path.

\emph{Step 2: $\mathcal D_m$ implies $U_m,\ldots,U_s,T_s$.} For $i>m$, $D^{(i)}$ gives $U_i$ by Lemma \ref{lemma:admissible}(i). By the same lemma, $D^{(m)}$ and $D^{(i)}$ ($i>m$) give that $(\zeta_m,\beta)$ and $(\zeta_i,\zeta_{i-1})$ are $I$-admissible, so $(\zeta_s,\beta)$ is $I$-admissible by Lemma \ref{lemma:admissible}(ii): this is $T_s$. Likewise the pairs $(\zeta_m,\beta),(\zeta_{m+1},\zeta_m),\ldots,(\zeta_s,\zeta_{s-1})$ are $\{\zeta_{m-1}\}$-admissible, which as in the proof of Lemma \ref{lemma:admissible}(ii) means:
\begin{equation}\label{eq:twoblocks}
\parbox{0.85\textwidth}{\centering the values $\beta,\zeta_m,\zeta_{m+1},\ldots,\zeta_s$ decrease on each side of $\zeta_{m-1}$,\\ with at most one jump from below $\zeta_{m-1}$ to above it.}
\end{equation}
For $U_m$, admissibility of $(\zeta_m,\beta)$ with respect to $\{\zeta_{m-1}\}$ leaves three possibilities: (1) $\zeta_{m-1}>\beta>\zeta_m$; (2) $\zeta_m>\zeta_{m-1}>\beta$; (3) $\beta>\zeta_m>\zeta_{m-1}$. We use $D^{(m)}$ restricted to $J_m$ and $D^{\mathrm{last}}$ restricted to $N_m$, and treat the remaining elements $\zeta_l$ ($l>m$, $z_l>z_m$) of $J_m$ by \eqref{eq:twoblocks}. In case (1), $J_m\cap(\zeta_m,\beta)=\emptyset$ and $N_m\cap(\beta,\zeta_{m-1})=\emptyset$, while by \eqref{eq:twoblocks} each $\zeta_l$ ($l>m$) is either less than $\zeta_m$ or greater than $\zeta_{m-1}$; hence $J_m\cap(\zeta_m,\zeta_{m-1})=\emptyset$. In case (2), $J_m\subseteq(\beta,\zeta_m)$ and $N_m\cap(\beta,\zeta_{m-1})=\emptyset$ give $N_m\subseteq(\zeta_{m-1},\zeta_m)$, while by \eqref{eq:twoblocks} the jump occurs at $(\beta,\zeta_m)$ and so $\zeta_{m-1}<\zeta_l<\zeta_m$ for $l>m$; hence $J_m\subseteq(\zeta_{m-1},\zeta_m)$. In case (3), $J_m\cap(\zeta_m,\beta)=\emptyset$ and $N_m\subseteq(\zeta_{m-1},\beta)$ give $N_m\subseteq(\zeta_{m-1},\zeta_m)$, while by \eqref{eq:twoblocks} all of $\beta,\zeta_m,\ldots,\zeta_s$ lie above $\zeta_{m-1}$ and decrease, so again $\zeta_{m-1}<\zeta_l<\zeta_m$ for $l>m$; hence $J_m\subseteq(\zeta_{m-1},\zeta_m)$. In each case this is $U_m$. Together with the prefix $U_1,\ldots,U_{m-1}$ this proves (c).

For the remainder assume that $F_{\mathrm{up}}$ is a cover path. Then the cycle $(z_s,\ldots,z_1,k)$ is $\psi$-tight, so by Corollary \ref{corollary:cyclic} the cyclic sequence $\Gamma=(\zeta_0,\zeta_1,\ldots,\zeta_s)$ is unimodal; let $\delta\in\{0,\ldots,s\}$ be the index of its ascent ($\zeta_\delta>\zeta_{\delta-1}$ if $\delta\geq1$; $\delta=0$ means the ascent is the wrap $(\zeta_s,\zeta_0)$), and let $g\in\{0,\ldots,s\}$ be the insertion gap of $\beta$ in $\Gamma$, where gap $g\geq1$ is the gap $(\zeta_{g-1},\zeta_g)$ and gap $0$ is the wrap. Explicitly: if $\min\Gamma<\beta<\max\Gamma$, then $g$ is determined by $\zeta_{g-1}>\beta>\zeta_g$ (indices modulo $s+1$); otherwise $g=\delta$.

\emph{Step 3: $\mathcal S_m$ implies $m>g$.} $\mathcal S_m$ says that the cycle $(z_s,\ldots,z_m,k,B)$ is $\phi_m$-tight, so by Corollary \ref{corollary:cyclic} the cyclic sequence $(\beta,\zeta_{m-1},\ldots,\zeta_s)$ is unimodal; that is, the wrap gap $(\zeta_s,\zeta_{m-1})$ of the block $\Gamma'=(\zeta_{m-1},\ldots,\zeta_s)$ is the insertion gap of $\beta$ in $\Gamma'$. Suppose $1\leq m\leq g$. Then $\zeta_{g-1},\zeta_g$ both lie in $\Gamma'$ and are consecutive there. If $\min\Gamma<\beta<\max\Gamma$ then $\zeta_{g-1}>\beta>\zeta_g$, so the insertion gap of $\beta$ in $\Gamma'$ is $(\zeta_{g-1},\zeta_g)$, which is not the wrap gap of $\Gamma'$ since $g-1\neq s$. If $\beta$ lies outside $[\min\Gamma,\max\Gamma]$ then $g=\delta$, so $\zeta_g=\max\Gamma$ and $\zeta_{g-1}=\min\Gamma$ lie in $\Gamma'$, whence $\beta$ lies outside the range of $\Gamma'$ and its insertion gap in $\Gamma'$ is the ascent gap $(\zeta_{g-1},\zeta_g)$ of $\Gamma'$, again not the wrap gap. In both cases $\mathcal S_m$ fails.

\emph{Step 4: $\mathcal S_{g+1}$ holds; hence $\mathcal S_m$ holds if and only if $m>g$.} We must verify $T_g$ and $A_i$ for $i>g$; by Lemma \ref{lemma:admissible}(i) we treat the sets $J_i$, $\{\beta\}$ and $I$ separately. Admissibility of $(\zeta_i,\zeta_{i-1})$ with respect to $J_i$ is $U_i$. With respect to $\{\beta\}$: if $i>g$ and $i\neq\delta$, the gap $(\zeta_{i-1},\zeta_i)$ is a descent gap and $\beta\notin(\zeta_i,\zeta_{i-1})$, for otherwise it would be the insertion gap of $\beta$; if $i=\delta>g$, then $\beta\in(\min\Gamma,\max\Gamma)=(\zeta_{\delta-1},\zeta_\delta)$, for otherwise $g=\delta$. With respect to $I$, and for $T_g$: the sequence $s_0=\beta,s_1=\zeta_g,\ldots,s_n=\zeta_s$ is a consecutive block of the unimodal sequence obtained by inserting $\beta$ into $\Gamma$, hence unimodal, and $T_s$ says that $(s_n,s_0)$ is $I$-admissible. If $s_0>s_1>\cdots>s_n$, then $T_s$ is classical, $I\cap(s_n,s_0)=\emptyset$, and every pair $(s_i,s_{i-1})$ is $I$-admissible. Otherwise the linear ascent is at some $j\geq1$, so $s_j>\cdots>s_n>s_0>\cdots>s_{j-1}$; then $T_s$ is quantum, $I\subseteq(s_0,s_n)$, and the pairs $(s_i,s_{i-1})$ with $i<j$ lie below $s_0$, the pair $(s_j,s_{j-1})$ has $I\subseteq(s_0,s_n)\subseteq(s_{j-1},s_j)$, and the pairs with $i>j$ lie above $s_n$; all are $I$-admissible. The pairs $(s_i,s_{i-1})$ are exactly $(\zeta_g,\beta)$ and $(\zeta_i,\zeta_{i-1})$, $i>g$. This proves $\mathcal S_{g+1}$, and the second assertion follows from Steps 1 and 3.

\emph{Step 5: $\mathcal D_m$ implies $m=g$.} Restricting $D^{(m)}$ to $\{\zeta_{m-1}\}$ leaves the three cases (1)--(3) of Step 2. In case (1) the gap $(\zeta_{m-1},\zeta_m)$ is a descent gap containing $\beta$, so $g=m$. In case (2) the gap is the ascent gap, $m=\delta$, and $\beta<\zeta_{m-1}=\min\Gamma$, so $g=\delta=m$. In case (3), $m=\delta$ and $\beta>\zeta_m=\max\Gamma$, so again $g=m$.

\emph{Step 6: if $g\geq1$ then $\mathcal D_g$ holds.} Write $a=z_g$. We are in one of the cases (1) $\zeta_{g-1}>\beta>\zeta_g$ ($g\neq\delta$), (2) $g=\delta$ and $\beta<\zeta_{g-1}=\min\Gamma$, (3) $g=\delta$ and $\beta>\zeta_g=\max\Gamma$. For $D^{(g)}$: admissibility of $(\zeta_g,\beta)$ with respect to $I$ is $T_g$ (Step 4); with respect to $\{\zeta_{g-1}\}$ it holds in each of the three cases by inspection; with respect to $J_g$ it follows from $U_g$: in case (1), $J_g\cap(\zeta_g,\zeta_{g-1})=\emptyset$ gives $J_g\cap(\zeta_g,\beta)=\emptyset$; in case (2), $J_g\subseteq(\zeta_{g-1},\zeta_g)\subseteq(\beta,\zeta_g)$; in case (3), $J_g\subseteq(\zeta_{g-1},\zeta_g)$ is disjoint from $(\zeta_g,\beta)$. For $D^{(i)}$, $i>g$: if $z_g>z_i$ the required statement is $A_i$ (Step 4); if $z_g<z_i$ then $z_g\in(z_i,k)$, so $\zeta_{g-1}=\phi_i(z_g)\in J_i$ already, $J^\ast_i=J_i$, and the statement follows from $U_i$ and the $I$-admissibility established in Step 4. For $D^{\mathrm{last}}$, the pair $(\beta,\zeta_{g-1})$ must be admissible with respect to $N_g\subseteq J_g$ and to the values $\zeta_{l-1}$, $l>g$, $z_l>z_g$: in case (1) it is classical, and $(\beta,\zeta_{g-1})\subseteq(\zeta_g,\zeta_{g-1})$ contains no element of $J_g$ (by $U_g$) and no $\zeta_l$ (as $\zeta_{g-1},\zeta_g$ are consecutive in the decreasing order of $\Gamma$); in case (2) it is classical, and $(\beta,\zeta_{g-1})$ lies below $\min\Gamma$ and below $J_g\subseteq(\zeta_{g-1},\zeta_g)$; in case (3) it is quantum, and $J_g\subseteq(\zeta_{g-1},\zeta_g)\subseteq(\zeta_{g-1},\beta)$ while every $\zeta_{l-1}$ with $l>g$ lies in $(\zeta_{g-1},\zeta_g]\subseteq(\zeta_{g-1},\beta)$.

\emph{Step 7: weights.} Suppose $F_{\mathrm{up}}$ is a cover path and $g\geq1$. The two paths $F_{\mathrm{up}}$ and $F_g$ share the prefix $t_{z_1k},\ldots,t_{z_{g-1}k}$; after it, the quantum steps of $F_{\mathrm{up}}$ are $t_{z_\delta k}$ (if $\delta\geq g$) and $t_{kB}$ (if $\zeta_s>\beta$), while those of $F_g$ are $t_{z_gB}$ (if $\zeta_g>\beta$), $t_{z_\delta B}$ (if $\delta>g$) and $t_{z_gk}$ (if $\beta>\zeta_{g-1}$). In case (1) with $\delta>g$ we have $\zeta_s>\zeta_0\geq\zeta_{g-1}>\beta$, and both weights are $q_{z_\delta k}q_{kB}=q_{z_\delta B}$; in case (1) with $\delta<g$ we have $\zeta_s<\zeta_g<\beta$ and both weights are $1$; in case (2) both are $q_{z_gk}q_{kB}=q_{z_gB}$; in case (3) both are $q_{z_gk}$. If instead $g=0$, so that $\mathcal S_1$ holds, the quantum steps of $F_{\mathrm{std}}$ are $t_{kB}$ (if $\zeta_0>\beta$) and $t_{z_\delta B}$ (if $\delta\geq1$); if $\zeta_s>\beta>\zeta_0$ (so $\delta\geq1$) both $F_{\mathrm{std}}$ and $F_{\mathrm{up}}$ have weight $q_{z_\delta B}$; if $\beta<\min\Gamma$ (so $\delta=0$) both have weight $q_{kB}$; and if $\beta>\max\Gamma$ both have weight $1$.

\emph{Conclusion.} (a): $\mathcal S_1$ gives $F_{\mathrm{up}}$ (Step 1) with the same weight (Step 7), and $g=0$ (Step 3), so no $\mathcal D_m$ holds (Step 5). (b): if $F_{\mathrm{up}}$ is a cover path and $\mathcal S_1$ fails, then $g\geq1$ (Step 4), $\mathcal D_g$ holds (Step 6) and no other $\mathcal D_m$ does (Step 5); as $U_1,\ldots,U_{g-1}$ hold, $F_g$ is the unique $F_m$ which is a cover path, and its weight is that of $F_{\mathrm{up}}$ (Step 7). (c) is Step 2.
\end{proof}

\begin{lemma} \label{lemma:pieriknotkn1}
Suppose $u,w\in S_\infty$ and let $k>1$ be an integer. If $u\tom{k} w$ but $u\ntom{k-1} w$, then there is a unique $b>k$ such that $wt_{kb}\lup w$ and $u\tom{k-1} wt_{kb}$, namely $b=B$. Furthermore, in that case there does not exist an $a<k$ such that $wt_{ak}\lup w$ and $u\tom{k-1} wt_{ak}$, and
$$\kappa(wt_{kB},t_{kB})\,q_{k-1}(u,wt_{kB})=q_k(u,w),\qquad \fx_{k-1}(u,wt_{kB})=\fx_k(u,w).$$
\end{lemma}
\begin{proof}
By Lemma \ref{lemma:piericycle}, $\pi$ is $k$-shaped and $u$-tight but not $(k-1)$-shaped, so $u(k)\neq w(k)$ and $B>k$ (Lemma \ref{lemma:leveldichotomy}). Write $C=(a_p,\ldots,a_1,B)$ with $a_j=k$, let $\rho$ be the product of the other nontrivial cycles of $\pi$, let $C'=(a_{j-1},\ldots,a_1,B)$, and put $\psi=u\rho C'$ and $X=(a_p,\ldots,a_{j+1},k,B)$, so that $w=\psi X$ and $X$ has the form of Lemma \ref{lemma:updown} with $s=p-j$ and $z_i=a_{j+i}$. Take the pivot chain of $\pi$ consisting of the pivot factorization of $\rho$, then that of $C$; by Lemma \ref{lemma:tight}(iii),(v) all its steps are covers, $\rho$ is $u$-tight, $C'$ is $u\rho$-tight (it is $C$ truncated after $j-1$ factors), and $F_{\mathrm{std}}$ is a cover path from $\psi$. By Lemma \ref{lemma:updown}(a), $F_{\mathrm{up}}$ is a cover path of the same weight. Its first $s$ steps form the pivot factorization of $C''=(a_p,\ldots,a_{j+1},k)$, so $\rho C'C''=u^{-1}wt_{kB}$ is $u$-tight by Lemma \ref{lemma:tight}(iii),(iv); it is $(k-1)$-shaped, its cycles being those of $\rho$ (each with exactly one element greater than $k$), $C'$ (top $B$, other elements less than $k$) and $C''$ (top $k$, other elements less than $k$). Hence $u\tom{k-1}wt_{kB}$, and the last step of $F_{\mathrm{up}}$ says $wt_{kB}\lup w$. Uniqueness of $b$ is Lemma \ref{lemma:neighbors}(U). The weights are $\qw_u(\rho)\qw_u(C')$ times the weights of $F_{\mathrm{std}}$, respectively $F_{\mathrm{up}}$, which agree; and $\fx_{k-1}(u,wt_{kB})=\fx_{k-1}(u,w)=\fx_k(u,w)$ since $u(k)\neq w(k)$.

Suppose finally that $a<k$ satisfies $wt_{ak}\lup w$ and $u\tom{k-1}wt_{ak}$. By Lemma \ref{lemma:neighbors} we are in case (S): $a=a_m$ with $m>j$, and $\sigma=\rho C'_mC''_m$ is $u$-tight. Taking the pivot chain of $\sigma$ in the order: factors of $\rho$, the first $j-1$ factors of $C'_m$ (which are those of $C'$), the factors of $C''_m$, the remaining factors $t_{a_mB},\ldots,t_{a_pB}$ of $C'_m$, and appending the cover $t_{a_mk}$, we obtain that $F_{m-j}$ is a cover path from $\psi$, contradicting Lemma \ref{lemma:updown}(a).
\end{proof}

\subsubsection*{The case $B=k$: the peel and the merges}

\begin{lemma}[Peel] \label{lemma:peel}
Let $\omega\in S_\infty$ and let $C=(a_p,\ldots,a_1,k)$ ($p\geq1$) be a cycle with top $k$; put $w=\omega C$, $C_-=(a_{p-1},\ldots,a_1,k)$ and $v=wt_{a_pk}=\omega C_-$. Then $C$ is $\omega$-tight if and only if $C_-$ is $\omega$-tight and $v\lup w$, and in that case $\kappa(v,t_{a_pk})\,\qw_\omega(C_-)=\qw_\omega(C)$.
\end{lemma}
\begin{proof}
The pivot factorization of $C$ is that of $C_-$ followed by $t_{a_pk}$, and $\omega C_-=v$. Apply Lemma \ref{lemma:tight}(iii).
\end{proof}

\begin{lemma}[Merge] \label{lemma:merge}
Let $\omega\in S_\infty$, $k>1$, $\eta=\omega(k)$, and let $D$ be a cycle with $\supp(D)\subseteq[1,k-1]$, allowing $D$ to be trivial, $D=\{a\}$ with $a\leq k-1$. Put $w=\omega D$, so $w(k)=\eta$. For $a\in\supp(D)$ let $\tau_a=Dt_{ak}$, the cycle with support $\supp(D)\cup\{k\}$ sending $a\mapsto k\mapsto D(a)$, and $v_a=wt_{ak}=\omega\tau_a$. Then:
\begin{enumerate}
\item[(a)] there is at most one $a\in\supp(D)$ such that $\tau_a$ is $\omega$-tight and $v_a\lup w$;
\item[(b)] such an $a$ exists if and only if $D$ is $\omega$-tight and $\max\supp(D)=k-1$ (for trivial $D=\{a\}$: if and only if $a=k-1$);
\item[(c)] for this $a$, $\kappa(v_a,t_{ak})\,\qw_\omega(\tau_a)=q_{k-1}\,\qw_\omega(D)$.
\end{enumerate}
\end{lemma}
\begin{proof}
\emph{Trivial $D=\{a\}$.} Then $w=\omega$, $\tau_a=t_{ak}$, and $\tau_a$ is $\omega$-tight if and only if $w\lup wt_{ak}$. If also $wt_{ak}\lup w$, then both $\ell(wt_{ak})-\ell(w)$ and its negative lie in $\{1,\,1-2(k-a)\}$, which forces $k-a=1$. Conversely for $a=k-1$ the transposition $t_{k-1,k}$ changes the length by $\pm1$, so both relations hold, one classical and one quantum with $\kappa=q_{k-1}$; thus $\kappa(v_a,t_{ak})\qw_\omega(\tau_a)=q_{k-1}=q_{k-1}\qw_\omega(D)$.

\emph{Nontrivial $D=(d_r,\ldots,d_1,d)$, $r\geq1$, $d=\max\supp(D)$.} Put $d_0=d$, $\gamma_0=\omega(d)$, $\gamma_i=\omega(d_i)$ and $\epsilon=k-d\geq1$. For $a=d_m$ ($0\leq m\leq r$) the cycle $\tau_m:=\tau_{d_m}$ is
$$\tau_m=(d_{m-1},\ldots,d_1,d,d_r,\ldots,d_m,k),$$
whose pivot factorization
$$t_{d_mk}\,t_{d_{m+1}k}\cdots t_{d_rk}\,t_{dk}\,t_{d_1k}\cdots t_{d_{m-1}k}$$
exchanges, in order, the values $\eta,\gamma_m,\gamma_{m+1},\ldots,\gamma_r,\gamma_0,\gamma_1,\ldots,\gamma_{m-1}$: the sequence of values of $\tau_m$ in the sense of Corollary \ref{corollary:cyclic} is the cyclic sequence $\Gamma=(\gamma_0,\ldots,\gamma_r)$ of $D$ with $\eta$ inserted into the gap $(\gamma_{m-1},\gamma_m)$ (the wrap gap $(\gamma_r,\gamma_0)$ when $m=0$). Note $w(d_m)=\gamma_{m-1}$ for $m\geq1$ and $w(d)=\gamma_r$; we write $w(a)$ for the value at $a=d_m$ in either case.

\emph{Necessity.} Suppose $\tau_m$ is $\omega$-tight and $v_m:=v_{d_m}\lup w$. By Corollary \ref{corollary:cyclic} the cyclic sequence of $\tau_m$ is unimodal; deleting $\eta$ shows that $\Gamma$ is unimodal, and then $m$ must be the insertion gap of $\eta$ in $\Gamma$. This already proves (a). Let $\delta$ be the ascent index of $\Gamma$ ($\gamma_\delta>\gamma_{\delta-1}$, or $\delta=0$ for the wrap), so that by Corollary \ref{corollary:cyclic} and Definition \ref{def:tight}
$$\Lw_\omega(D)=r-2(d-d_\delta),\qquad \qw_\omega(D)=q_{d_\delta d}$$
(read as $r$ and $1$ when $\delta=0$). Now compare lengths. Tightness gives $\ell(v_m)-\ell(\omega)=\Lw_\omega(\tau_m)$, Lemma \ref{lemma:tight}(ii) gives $\ell(w)-\ell(\omega)\geq\Lw_\omega(D)$ with equality if and only if $D$ is $\omega$-tight, and $v_m\lup w$ gives $\ell(v_m)-\ell(w)=-1$ if $w(a)>\eta$ and $\ell(v_m)-\ell(w)=2(k-a)-1$ if $w(a)<\eta$. Hence
\begin{equation}\label{eq:mergeineq}
\Lw_\omega(\tau_m)-\Lw_\omega(D)\ \geq\ \begin{cases}-1&\text{if }w(a)>\eta,\\ 2(k-a)-1&\text{if }w(a)<\eta,\end{cases}
\end{equation}
with equality if and only if $D$ is $\omega$-tight. We compute the left side from Corollary \ref{corollary:cyclic} in the three possible configurations of the insertion gap.
\begin{itemize}
\item[(I)] $\min\Gamma<\eta<\max\Gamma$, so $\gamma_{m-1}>\eta>\gamma_m$ and $m\neq\delta$. The ascent of the sequence of $\tau_m$ is the (interior) step onto $\gamma_\delta$, whose position is $d_\delta$, so $\Lw_\omega(\tau_m)=(r+1)-2(k-d_\delta)$ and $\Lw_\omega(\tau_m)-\Lw_\omega(D)=1-2\epsilon$. Here $w(a)=\gamma_{m-1}>\eta$, so \eqref{eq:mergeineq} reads $1-2\epsilon\geq-1$, forcing $\epsilon=1$ and equality.
\item[(II)] $\eta<\min\Gamma$, so $m=\delta$. The ascent of the sequence of $\tau_m$ is its first step $(\eta,\gamma_\delta)$, at position $d_\delta$; again $\Lw_\omega(\tau_m)-\Lw_\omega(D)=1-2\epsilon$, and $w(a)=\gamma_{\delta-1}>\eta$ forces $\epsilon=1$ and equality.
\item[(III)] $\eta>\max\Gamma$, so $m=\delta$. The sequence of $\tau_m$ is decreasing, its ascent being the wrap, so $\Lw_\omega(\tau_m)=r+1$ and $\Lw_\omega(\tau_m)-\Lw_\omega(D)=1+2(d-d_\delta)$. Here $w(a)=\gamma_{\delta-1}<\eta$, and \eqref{eq:mergeineq} reads $1+2(d-d_\delta)\geq2(k-d_\delta)-1$, i.e.\ $\epsilon\leq1$; so $\epsilon=1$ and equality.
\end{itemize}
In every case $d=k-1$ and $D$ is $\omega$-tight, proving the ``only if'' part of (b).

\emph{Sufficiency.} Suppose $D$ is $\omega$-tight and $d=k-1$; then $\Gamma$ is unimodal (Corollary \ref{corollary:cyclic}). Let $m$ be the insertion gap of $\eta$ and $a=d_m$. The computations above show that \eqref{eq:mergeineq} holds with equality, and since $\ell(w)-\ell(\omega)=\Lw_\omega(D)$, the tightness of $\tau_m$ is equivalent to $\ell(v_m)-\ell(w)$ taking the value $-1$ (if $w(a)>\eta$) or $2(k-a)-1$ (if $w(a)<\eta$), i.e.\ to $v_m\lup w$. It therefore suffices to show $wt_{ak}\lup w$, that is, that $(w(a),\eta)$ is $J$-admissible for $J=\{w(r)\mid a<r<k\}$. The positions in $(a,k)$ are: $d=k-1$ (if $a\neq d$), carrying $\gamma_r$; the $d_l>a$, carrying $\gamma_{l-1}$; and positions outside $\supp(D)$, carrying their $\omega$-values. If $a=d=k-1$ the set $J$ is empty and we are done; so assume $m\geq1$.
\begin{itemize}
\item[(I)] $\gamma_{m-1}>\eta>\gamma_m$: the pair $(w(a),\eta)=(\gamma_{m-1},\eta)$ is classical and we need $J\cap(\eta,\gamma_{m-1})=\emptyset$. No $\gamma_l$ lies in $(\gamma_m,\gamma_{m-1})$, these being consecutive in the decreasing order of $\Gamma$. The $m$-th step of the pivot chain of $D$ is classical, so by \eqref{eq:Cj} every position $r\in(d_m,d)$ outside $\supp(D)$ has $\omega(r)\notin(\gamma_m,\gamma_{m-1})$.
\item[(II)] $m=\delta$, $\eta<\gamma_{\delta-1}=\min\Gamma$: the pair $(\gamma_{\delta-1},\eta)$ is classical and we need $J\cap(\eta,\gamma_{\delta-1})=\emptyset$; all $\gamma_l\geq\gamma_{\delta-1}$, and by \eqref{eq:Qj} for the quantum step $\delta$ of $D$ the positions $r\in(d_\delta,d)$ outside $\supp(D)$ have $\omega(r)\in(\gamma_{\delta-1},\gamma_\delta)$.
\item[(III)] $m=\delta$, $\eta>\gamma_\delta=\max\Gamma$: the pair $(\gamma_{\delta-1},\eta)$ is quantum and we need $J\subseteq(\gamma_{\delta-1},\eta)$; the values $\gamma_l$ occurring in $J$ are different from $\gamma_{\delta-1}$ (which sits at $a$) and lie in $(\gamma_{\delta-1},\gamma_\delta]$, and the remaining positions have $\omega(r)\in(\gamma_{\delta-1},\gamma_\delta)$ by \eqref{eq:Qj}.
\end{itemize}
This proves the ``if'' part of (b).

\emph{Weights.} In cases (I) and (II), $\kappa(v_a,t_{ak})=1$ and $\qw_\omega(\tau_m)=q_{d_\delta k}=q_{d_\delta,k-1}q_{k-1}=\qw_\omega(D)q_{k-1}$ (read $q_{dk}=q_{k-1}$ and $\qw_\omega(D)=1$ when $\delta=0$). In case (III), $\kappa(v_a,t_{ak})=q_{ak}=q_{d_\delta k}$ and $\qw_\omega(\tau_m)=1$, giving the same product. This is (c).
\end{proof}

\subsubsection*{The cancellations}

We now assemble the Monk summands. Throughout, $\rho$ denotes the product of the nontrivial cycles of $\pi$ other than $C$ and other than the cycle $D$ of $\pi$ containing $k-1$ (we set $D=\{k-1\}$ if $k-1$ is fixed by $\pi$, and we do not single out $D$ when $k-1\in C$). By Lemma \ref{lemma:neighbors}, Lemma \ref{lemma:piericycle} and Lemma \ref{lemma:tight}(iv), each of the quantities $c_{k-1}(u,w)$, $c_{k-2}(u,w)$, $c_k(u,w)$ and each Monk summand vanishes unless $\rho$ is $u$-tight and each cycle of $\rho$ contains exactly one element greater than $k$; indeed each requires $\rho$ to be a tight product of cycles with exactly one element greater than $k$ or than $k-1$ or than $k-2$, and these conditions agree for cycles containing neither $k$ nor $k-1$. (For $c_{k-2}(u,w)$, note that $k-1$ lies in $C$ or in $D$, hence not in a cycle of $\rho$.) \emph{We assume this from now on}; otherwise all terms considered below are zero.

\begin{lemma}[Reflection--drop cancellation] \label{lemma:reflectioncancel}
Suppose $u,w\in S_\infty$, let $k>1$, and suppose $u(k)=w(k)$. Then
$$\sum_{\substack{a<k,\ wt_{ak}\lup w}}\!\!\!\mathrm{sign}(a-k)\,\kappa_{ak}\,c_{k-1}(u,wt_{ak})\ +\ \sum_{\substack{b>k,\ wt_{kb}\lup w}}\!\!\!\mathrm{sign}(b-k)\,\kappa_{kb}\,c_{k-1}(u,wt_{kb})\ +\ q_{k-1}\,c_{k-2}(u,w)=0,$$
where $\kappa_{rs}=\kappa(wt_{rs},t_{rs})$. In fact the two reflection sums contribute at most one nonzero summand in total: it occurs precisely when $u\tom{k-2}w$ (equivalently, when the quantum term is nonzero), it is indexed by some $a<k$, and it equals $-\,q_{k-1}c_{k-2}(u,w)$.
\end{lemma}
\begin{proof}
Here $C=\{k\}$, so by Lemma \ref{lemma:neighbors} only summands of type (M) can be nonzero: $t=t_{ak}$ with the $\pi$-cycle $D_a$ of $a$ contained in $[1,k-1]$. If $D_a$ is nontrivial then it is not a cycle of $\rho$ (those have an element greater than $k$), so $D_a=D$ and $\max D=k-1$; otherwise $a$ is a fixed point of $\pi$. Let $a$ be such a candidate and put $\omega_a=wD_a^{-1}$, so that $w=\omega_aD_a$, $\omega_a(k)=w(k)$, and $\omega_a=u\rho$ when $D_a=D$. Then $u^{-1}wt_{ak}=(\pi D_a^{-1})\,\tau_a$ with $\tau_a=D_at_{ak}$ as in Lemma \ref{lemma:merge}, where $\pi D_a^{-1}$ is the product of the cycles of $\pi$ other than $D_a$; by Lemma \ref{lemma:tight}(iv), $u^{-1}wt_{ak}$ is $u$-tight if and only if $\pi D_a^{-1}$ is $u$-tight and $\tau_a$ is $\omega_a$-tight. Now apply Lemma \ref{lemma:merge} to $(\omega_a,D_a)$. If $a$ is a fixed point, the summand can be nonzero only for $a=k-1$, which requires $k-1$ to be fixed, i.e.\ $D=\{k-1\}$; in that case $\pi D_a^{-1}=\rho$ and the summand is nonzero. If $D_a=D$ is nontrivial, at most one $a\in D$ works, and one does if and only if $D$ is $u\rho$-tight. Altogether at most one $a$ gives a nonzero summand, and one exists if and only if either $D=\{k-1\}$, or $D$ is nontrivial with $\max D=k-1$ and $u\rho$-tight. On the other hand $u\tom{k-2}w$ if and only if $\pi=\rho D$ is $(k-2)$-shaped and $u$-tight, i.e.\ (given our standing assumption on $\rho$) if and only if $D$ is trivial or has $k-1$ as its unique element greater than $k-2$, and $D$ is $u\rho$-tight: the same condition. When it holds, with $v=wt_{ak}$ the unique contributing neighbor and $\omega=\omega_a=u\rho$, Lemma \ref{lemma:merge}(c) and Lemma \ref{lemma:piericycle} give
$$\kappa_{ak}\,q_{k-1}(u,v)=\qw_u(\rho)\,\kappa_{ak}\qw_\omega(\tau_a)=\qw_u(\rho)\,q_{k-1}\qw_\omega(D)=q_{k-1}\,q_{k-2}(u,w),$$
and Lemma \ref{lemma:neighbors}(M) gives $\fx_{k-1}(u,v)=\fx_{k-1}(u,w)\setminus\{u(a)\}$, which equals $\fx_{k-2}(u,w)$: if $D=\{k-1\}$ then $a=k-1$ and we remove the fixed value $u(k-1)$; if $D$ is nontrivial then neither $a$ nor $k-1$ is a fixed position of $w$. Hence the summand is $-\kappa_{ak}c_{k-1}(u,v)=-q_{k-1}c_{k-2}(u,w)$, cancelling the quantum term. If the condition fails, both the reflection sums and the quantum term vanish.
\end{proof}

\begin{lemma}[Moved reflection cancellation] \label{lemma:movedcancel}
Let $u,w\in S_\infty$, let $k>1$, and suppose $u(k)\neq w(k)$, so that $C$ is nontrivial.
\begin{enumerate}
\item[(a)] If $B>k$, then $c_{k-1}(u,w)=c_{k-2}(u,w)=0$, and the sum of the Monk summands equals $c_k(u,w)$. More precisely: if $u\tom{k}w$ the only nonzero Monk summand is that of $t_{kB}$, and it equals $c_k(u,w)$; if $u\ntom{k}w$, the nonzero Monk summands are either none or exactly two, one of type (U) and one of type (S), with equal absolute value and opposite signs.
\item[(b)] If $B=k$, then $c_k(u,w)=0$, no summand of type (U) occurs, and
$$(y_{w(k)}-z_1)\,c_{k-1}(u,w)\;-\!\!\sum_{\substack{a<k,\ wt_{ak}\lup w}}\!\!\kappa_{ak}\,c_{k-1}(u,wt_{ak})\;+\;q_{k-1}c_{k-2}(u,w)=0.$$
In fact the diagonal term cancels against the summand of type (P), and the quantum term against the (at most one) summand of type (M).
\end{enumerate}
\end{lemma}
\begin{proof}
(a) $C$ contains $k$ and $B$, two elements greater than $k-1$, so $\pi$ is neither $(k-1)$- nor $(k-2)$-shaped and $c_{k-1}(u,w)=c_{k-2}(u,w)=0$. By Lemma \ref{lemma:neighbors}, a nonzero Monk summand is of type (U) or (S), and requires $C\cap(k-1,\infty)=\{k,B\}$; together with the standing assumption on $\rho$ this means that $\pi$ is $k$-shaped, which we may assume (otherwise all Monk summands vanish and $c_k(u,w)=0$). Write $C=(a_p,\ldots,a_1,B)$ with $a_j=k$, $C'=(a_{j-1},\ldots,a_1,B)$, $\psi=u\rho C'$ and $X=(a_p,\ldots,a_{j+1},k,B)$, so that $w=\psi X$ as in Lemma \ref{lemma:updown} with $s=p-j$, $z_i=a_{j+i}$. If $C'$ is not $u\rho$-tight then nothing is tight (Lemma \ref{lemma:tight}(iv),(v): $C'$ is a truncation of $C$ and of each $C'_m$, and a cycle of $u^{-1}wt_{kB}$), all Monk summands vanish, and $c_k(u,w)=0$. So assume $C'$ is $u\rho$-tight. Then, exactly as in the proof of Lemma \ref{lemma:pieriknotkn1}: $u\tom{k}w$ if and only if $F_{\mathrm{std}}$ is a cover path from $\psi$; the summand of type (U) is nonzero if and only if $F_{\mathrm{up}}$ is a cover path, and then equals $\qw_u(\rho)\qw_u(C')\,\mathrm{wt}(F_{\mathrm{up}})\,\wte(y_P;z)$ with $P=\fx_{k-1}(u,w)$; and the summand of type (S) indexed by $a_m$ ($m>j$) is nonzero if and only if $F_{m-j}$ is a cover path, and then equals $-\qw_u(\rho)\qw_u(C')\,\mathrm{wt}(F_{m-j})\,\wte(y_P;z)$ (Lemma \ref{lemma:neighbors}(S) for the fixed set). If $u\tom{k}w$, Lemma \ref{lemma:updown}(a) and Lemma \ref{lemma:pieriknotkn1} show that the only nonzero summand is that of $t_{kB}$ and that it equals $q_k(u,w)\wte(y_{\fx_k(u,w)};z)=c_k(u,w)$. If $u\ntom{k}w$, Lemma \ref{lemma:updown}(b),(c) show that either no summand is nonzero, or exactly two are, namely (U) and one (S), with equal weights; they cancel, and $c_k(u,w)=0$.

(b) $C\subseteq[1,k]$ has no element greater than $k$, so $\pi$ is not $k$-shaped and $c_k(u,w)=0$; type (U) is excluded by Lemma \ref{lemma:neighbors}. Write $C=(a_p,\ldots,a_1,k)$. For the diagonal and peel terms, let $\omega_0=u\rho D$ if $k-1\notin C$ and $\omega_0=u\rho$ if $k-1\in C$ (where $D$ is the cycle of $k-1$, trivial or not), so that $w=\omega_0C$ and $\omega_0$ agrees with $u$ on $\supp(C)$. The diagonal term is nonzero if and only if $\pi$ is $(k-1)$-shaped and tight, i.e.\ (given the standing assumption) $D$ is trivial or has an element greater than $k$, $D$ is $u\rho$-tight, and $C$ is $\omega_0$-tight (the conditions on $D$ being vacuous when $k-1\in C$). The summand of type (P) is nonzero if and only if $\rho D C_-$ is $(k-1)$-shaped and tight and $v=wt_{a_pk}\lup w$, i.e.\ the same conditions on $D$ together with: $C_-$ is $\omega_0$-tight and $v\lup w$. By Lemma \ref{lemma:peel} these two sets of conditions are equivalent, and when they hold, $\kappa_{a_pk}q_{k-1}(u,v)=\qw_u(\rho D)\kappa_{a_pk}\qw_{\omega_0}(C_-)=\qw_u(\rho D)\qw_{\omega_0}(C)=q_{k-1}(u,w)$, while $\fx_{k-1}(u,v)=\fx_{k-1}(u,w)\sqcup\{w(k)\}$ by Lemma \ref{lemma:neighbors}(P). Since $(y_{w(k)}-z_1)\wte(y_{\fx_{k-1}(u,w)};z)=\wte(y_{\fx_{k-1}(u,w)\cup\{w(k)\}};z)$, the diagonal term and the (negatively signed) summand of type (P) cancel.

For the merge and quantum terms, a summand of type (M) requires $a\notin C$ with the cycle $D_a$ of $a$ contained in $[1,k-1]$; as in the proof of Lemma \ref{lemma:reflectioncancel}, either $a$ is fixed by $\pi$ or $D_a=D$ with $\max D=k-1$ (in particular $k-1\notin C$). If $k-1\in C$, then $k-1$ is not fixed and no such $D$ exists; since a fixed-point candidate must be $a=k-1$ by the trivial case of Lemma \ref{lemma:merge} (applied as below), there are no nonzero summands of type (M), and $c_{k-2}(u,w)=0$ since $C$ contains the two elements $k-1,k$ greater than $k-2$. So suppose $k-1\notin C$. For a candidate $a$ put $\omega_a=wD_a^{-1}$, so that $w=\omega_aD_a$ and $\omega_a(k)=w(k)$; then $u^{-1}wt_{ak}=\rho_a\,E$, where $\rho_a$ is the product of the cycles of $\pi$ other than $C$ and $D_a$ and $E=C\tau_a$ is the cycle of Lemma \ref{lemma:neighbors}(M), $\tau_a=D_at_{ak}$. The pivot factorization of $E$ is that of $C$ followed by that of $\tau_a$, so by Lemma \ref{lemma:tight}(iii),(iv), $u^{-1}wt_{ak}$ is $u$-tight if and only if $\rho_a$ is $u$-tight, $C$ is $u\rho_a$-tight and $\tau_a$ is $\omega_a$-tight (note $u\rho_aC=\omega_a$), and then $\qw_u(E)=\qw_u(C)\qw_{\omega_a}(\tau_a)$. By Lemma \ref{lemma:merge} applied to $(\omega_a,D_a)$: a fixed-point candidate contributes only if $a=k-1$, i.e.\ $D=\{k-1\}$; and if $D_a=D$ is nontrivial, at most one $a\in D$ contributes, and one does if and only if $D$ is $\omega_a$-tight, where now $\omega_a=u\rho C=:\omega$. Altogether at most one summand of type (M) is nonzero, and one is if and only if $C$ is $u\rho$-tight and either $D=\{k-1\}$, or $D$ is nontrivial with $\max D=k-1$ and $\omega$-tight. This is exactly the condition for $u\tom{k-2}w$: $\pi=\rho CD$ is $(k-2)$-shaped (as $k-1\notin C$, $C$ has the unique element $k$ greater than $k-2$, and $D$ the unique element $k-1$) and it is $u$-tight if and only if $C$ is $u\rho$-tight and $D$ is $\omega$-tight. When it holds, Lemma \ref{lemma:merge}(c) gives
$$\kappa_{ak}q_{k-1}(u,wt_{ak})=\qw_u(\rho)\qw_u(C)\,\kappa_{ak}\qw_\omega(\tau_a)=\qw_u(\rho)\qw_u(C)\,q_{k-1}\qw_\omega(D)=q_{k-1}q_{k-2}(u,w),$$
and $\fx_{k-1}(u,wt_{ak})=\fx_{k-1}(u,w)\setminus\{u(a)\}=\fx_{k-2}(u,w)$ as in the proof of Lemma \ref{lemma:reflectioncancel}. Thus the summand of type (M) equals $-q_{k-1}c_{k-2}(u,w)$ and cancels the quantum term; when the condition fails, both vanish.
\end{proof}

We now record the recurrence satisfied by the quantum factorial elementary polynomials, which drives the induction.

\begin{lemma} \label{lemma:elemrecursion}
For all integers $k\geq2$ we have
$$\qelem_{k,k}(x;z)=(x_k-z_1)\qelem_{k-1,k-1}(x;z)+q_{k-1}\qelem_{k-2,k-2}(x;z),$$
with the conventions $\qelem_{0,0}=1$ and $\qelem_{k,k}=0$ if $k<0$. Moreover, the ordinary polynomials $\wte(y_A;z)=\prod_{i\in A}(y_i-z_1)$ satisfy, for any set $A$ and any $c\notin A$, the single-variable identity
$$\wte(y_{A\cup\{c\}};z)=(y_c-z_1)\,\wte(y_A;z).$$
\end{lemma}
\begin{proof}
By Lemma \ref{lemma:qelemkk}, $\qelem_{k,k}(x;z)=\Delta_k(-z_1\mid x_1,\ldots,x_k)$, and expanding the tridiagonal determinant $\Delta_k(t\mid x_1,\ldots,x_k)$ along its last row gives $\Delta_k=(x_k+t)\Delta_{k-1}+q_{k-1}\Delta_{k-2}$. The second identity is clear.
\end{proof}

The first recurrence is the quantization of the classical relation $E_{k,k}=(x_k-z_1)E_{k-1,k-1}$; the extra term $q_{k-1}\qelem_{k-2,k-2}$ is the only feature absent from the non-quantum setting. Discarding all quantum terms recovers the classical Pieri formula \cite[Theorem~7.1]{samuelmolev}, whose term-by-term cancellation is carried out in \cite[Proof of Proposition~4.4]{samuelmolev} by means of the reflection-matching Lemmas 7.6--7.9 of that paper; the lemmas above are the quantum counterparts, the essential new points being the length-decreasing covers permitted by Lemma \ref{lemma:piericycle}, the attendant bookkeeping of the powers of $q$, and the extra term of the recurrence, which is absorbed by the merge summands of Lemma \ref{lemma:merge}.

\begin{proof}[Proof of Theorem \ref{theorem:pieri}]
For a permutation $w$, write $[\qsch_w]F$ for the coefficient of $\qsch_w(x;y)$ in an expansion $F$ of a quantum double Schubert polynomial in the basis $\{\qsch_w(x;y)\}$. With this notation the theorem is the assertion that, for all $u$ and all $k\geq1$,
$$[\qsch_w]\big(\qsch_u(x;y)\qelem_{k,k}(x;z)\big)=\begin{cases}q_k(u,w)\,\wte(y_{\fx_{k}(u,w)};z) & \text{if }u\tom{k} w,\\[2pt] 0 & \text{otherwise.}\end{cases}$$
For brevity, when $u\tom{m} w$ we set
$$c_m(u,w)=q_m(u,w)\,\wte(y_{\fx_{m}(u,w)};z),$$
and $c_m(u,w)=0$ if $u\ntom{m}w$; note that $u\tom{0}w$ only for $w=u$, with $c_0(u,u)=1$, in accordance with $\qelem_{0,0}=1$. The claim is that $[\qsch_w](\qsch_u\qelem_{k,k})=c_k(u,w)$.

We argue by strong induction on $k$.

\emph{Base case $k=1$.} We have $\qelem_{1,1}(x;z)=x_1-z_1$, and Lemma \ref{lemma:singlevariable} with $i=j=1$ gives
$$(x_1-z_1)\qsch_u(x;y)=(y_{u(1)}-z_1)\qsch_u(x;y)+\sum_{\substack{b>1\\\ell(ut_{1b})=\ell(u)+1}}\qsch_{ut_{1b}}(x;y)+\sum_{\substack{b>1\\\ell(ut_{1b})=\ell(u)-2(b-1)+1}}q_{1b}\qsch_{ut_{1b}}(x;y).$$
The permutations $w$ with $u\tom{1}w$ are exactly $w=u$ and $w=ut_{1b}$ for $b>1$ with $u\lup ut_{1b}$, since Definition \ref{definition:pieri} permits at most one reflection when $k=1$. For $w=u$ we have $q_1=1$, $\fx_{1}(u,u)=\{u(1)\}$, and $\wte(y_{\{u(1)\}};z)=y_{u(1)}-z_1$, matching the diagonal term. For $w=ut_{1b}$ we have $\fx_{1}=\emptyset$ and $\wte(y_\emptyset;z)=1$; if the step is length-increasing then $q_1(u,w)=1$, matching the second sum, and if it is length-decreasing then $q_1(u,w)=q_{1b}$, matching the third sum. All other coefficients vanish on both sides. This proves the base case.

\emph{Inductive step.} Fix $k>1$ and assume the claim at all levels $m<k$. By Lemma \ref{lemma:elemrecursion},
\begin{equation}\label{eq:pieriexpand}
\qsch_u\qelem_{k,k}=(x_k-z_1)\big(\qsch_u\qelem_{k-1,k-1}\big)+q_{k-1}\,\qsch_u\qelem_{k-2,k-2}.
\end{equation}
By the inductive hypothesis, $\qsch_u\qelem_{k-1,k-1}=\sum_v c_{k-1}(u,v)\,\qsch_v$. Applying Lemma \ref{lemma:singlevariable} with $i=k$, $j=1$ to each $\qsch_v$ and extracting the coefficient of $\qsch_w$ from \eqref{eq:pieriexpand}, we obtain
\begin{align}
[\qsch_w]\big(\qsch_u\qelem_{k,k}\big)
&=(y_{w(k)}-z_1)\,c_{k-1}(u,w) \label{eq:diag}\\
&\quad+\sum_{\substack{b>k,\ w t_{kb}\lup w\\ \ell(w)=\ell(wt_{kb})+1}} c_{k-1}(u,wt_{kb})
+\sum_{\substack{b>k,\ w t_{kb}\lup w\\ \ell(w)=\ell(wt_{kb})-2(b-k)+1}} q_{kb}\,c_{k-1}(u,wt_{kb}) \label{eq:up}\\
&\quad-\sum_{\substack{a<k,\ w t_{ak}\lup w\\ \ell(w)=\ell(wt_{ak})+1}} c_{k-1}(u,wt_{ak})
-\sum_{\substack{a<k,\ w t_{ak}\lup w\\ \ell(w)=\ell(wt_{ak})-2(k-a)+1}} q_{ak}\,c_{k-1}(u,wt_{ak}) \label{eq:down}\\
&\quad+q_{k-1}\,c_{k-2}(u,w), \label{eq:rec}
\end{align}
where we used that $\mathrm{sign}(b-k)=1$ for $b>k$ and $\mathrm{sign}(a-k)=-1$ for $a<k$, that the term $\qsch_w$ arises from $\qsch_v$ with $v=wt$ exactly when $v\lup w$, and that a summand in \eqref{eq:up}--\eqref{eq:down} is present precisely when $u\tom{k-1}wt_{kb}$, respectively $u\tom{k-1}wt_{ak}$ (otherwise $c_{k-1}=0$). In the terminology introduced before Lemma \ref{lemma:neighbors}, \eqref{eq:up} and \eqref{eq:down} are the Monk summands at $w$, with $\kappa_{kb}$, respectively $\kappa_{ak}$, equal to $1$ in the first sum and to $q_{kb}$, respectively $q_{ak}$, in the second. We must show that \eqref{eq:diag}--\eqref{eq:rec} equals $c_k(u,w)$. We distinguish three cases according to whether $u(k)=w(k)$ and, when $u(k)\neq w(k)$, whether the cycle $C$ of $u^{-1}w$ containing $k$ has its top $B$ greater than $k$ or equal to $k$.

\emph{Case 1: $u(k)=w(k)$.} By Lemma \ref{lemma:leveldichotomy}(i) we have $u\tom{k}w$ if and only if $u\tom{k-1}w$; suppose first that both hold. Then, by the same lemma,
$$\fx_{k}(u,w)=\fx_{{k-1}}(u,w)\cup\{u(k)\},\qquad q_k(u,w)=q_{k-1}(u,w),$$
so the diagonal contribution \eqref{eq:diag} gives
$$(y_{u(k)}-z_1)\,c_{k-1}(u,w)=q_{k-1}(u,w)\,(y_{u(k)}-z_1)\,\wte(y_{\fx_{{k-1}}(u,w)};z)=q_k(u,w)\,\wte(y_{\fx_{k}(u,w)};z)=c_k(u,w),$$
by Lemma \ref{lemma:elemrecursion}. The remaining contributions are the reflection sums \eqref{eq:up}, \eqref{eq:down} and the quantum term $q_{k-1}c_{k-2}(u,w)$ of \eqref{eq:rec}, which sum to $0$ by Lemma \ref{lemma:reflectioncancel}. Hence the total is $c_k(u,w)$, as required. If instead $u(k)=w(k)$ but $u\ntom{k}w$ (equivalently $u\ntom{k-1}w$), then $c_k(u,w)=0$ and the diagonal term vanishes, while the reflection sums and quantum term again cancel by Lemma \ref{lemma:reflectioncancel}; the total is $0=c_k(u,w)$.

\emph{Case 2: $u(k)\neq w(k)$ and $B>k$.} By Lemma \ref{lemma:movedcancel}(a), the diagonal term \eqref{eq:diag} and the quantum term \eqref{eq:rec} vanish, and the Monk summands \eqref{eq:up}--\eqref{eq:down} sum to $c_k(u,w)$. (When $u\tom{k}w$, Lemma \ref{lemma:pieriknotkn1} identifies the single surviving summand: it is indexed by $t_{kB}$, with $v=wt_{kB}$, $u\tom{k-1}v$, $\fx_k(u,w)=\fx_{k-1}(u,v)$ and $\kappa_{kB}\,q_{k-1}(u,v)=q_k(u,w)$, the factor $\kappa_{kB}$ being $1$ or $q_{kB}$ according as $t_{kB}$ is a classical or a quantum cover of $v$; this is the step at which the quantum drop terms of Lemma \ref{lemma:singlevariable} and the quantum weight of Definition \ref{definition:pieri} interlock.)

\emph{Case 3: $u(k)\neq w(k)$ and $B=k$.} Here $u\ntom{k}w$, since $C$ contains no element greater than $k$, so $c_k(u,w)=0$; and \eqref{eq:diag}--\eqref{eq:rec} vanishes by Lemma \ref{lemma:movedcancel}(b).

By Lemma \ref{lemma:leveldichotomy} these cases are exhaustive. In all cases $[\qsch_w](\qsch_u\qelem_{k,k})=c_k(u,w)$, completing the induction and the proof.
\end{proof}

\section{Proofs of the expansion formulas} \label{section:proofvside}

In this section we prove Theorems \ref{theorem:orddouble}, \ref{theorem:quandouble} and \ref{theorem:quandoublemed}. Throughout, divided difference operators act on the $z$-variables, and in order to free the subscript for the skew operators we write
$$\partial^a=\partial^{(z)}_a,\qquad \partial^w=\partial^{(z)}_w=\partial^{a_1}\cdots\partial^{a_m}\quad(w=s_{a_1}\cdots s_{a_m}\text{ reduced}),$$
the rightmost operator being applied first. The simple reflection $s_a$ acts on polynomials by exchanging $z_a$ and $z_{a+1}$, and a permutation $w$ acts by $w(z_i)=z_{w(i)}$, so that $w=s_{a_1}\circ\cdots\circ s_{a_m}$ as operators. All operators of this section are linear over $\mathbb Z[q][x]$, which we treat as the ring of scalars. Recall that the quantization map $\mathcal Q$ of Fomin, Gelfand and Postnikov (Lemma \ref{lemma:quantumtoordinary}) is linear over $\mathbb Z[q][z]$, sends every standard elementary monomial $\prod_je_{i_j}(x_1,\ldots,x_j)$ to $\prod_je^{(q)}_{i_j}(x_1,\ldots,x_j)$, and satisfies $\mathcal Q(\sch_w(x;z))=\qsch_w(x;z)$.

\subsection{Skew divided difference operators}

We use the skew divided difference operators of \cite{samuelleibniz,samuelmolev}. They are the operators $\partial^w_u$, $u,w\in S_\infty$, determined by the Leibniz formula
\begin{equation}\label{eq:leibniz2}
\partial^w(PQ)=\sum_{u\in S_\infty}\partial^u(P)\,\partial^w_u(Q),
\end{equation}
which differs by a normalization from the convention of Macdonald \cite{notes}. We shall use the following facts.

\begin{lemma} \label{lemma:skewfacts}
Let $u,w\in S_\infty$.
\begin{enumerate}
\item[(i)] (Subword formula.) If $(a_1,\ldots,a_m)$ is a reduced word for $w$, then $\partial^w_u$ is the sum, over all reduced subwords $(b_1,\ldots,b_m)$ for $u$ (that is, $b_j\in\{1,s_{a_j}\}$ for each $j$, and the nonunit entries form a reduced word for $u$), of the compositions $\partial^{a_1}_{b_1}\cdots\partial^{a_m}_{b_m}$, where $\partial^{a}_{1}=\partial^a$ and $\partial^{a}_{s_a}=s_a$.
\item[(ii)] $\partial^w_u=0$ unless $u\leq w$ in Bruhat order; $\partial^w_1=\partial^w$ and $\partial^w_w=w$; and $\partial^w_u$ is homogeneous of degree $\ell(u)-\ell(w)$ in the $z$-variables. In particular, for a scalar $c$, $\partial^w_u(c)=c$ if $u=w$ and $\partial^w_u(c)=0$ otherwise.
\item[(iii)] (Ascending Leibniz formula.) For polynomials $f_1,\ldots,f_k$,
$$\partial^w_u(f_1f_2\cdots f_k)=\sum_{u=u_0\leq u_1\leq\cdots\leq u_k=w}\ \prod_{j=1}^k\partial^{u_j}_{u_{j-1}}(f_j).$$
\end{enumerate}
\end{lemma}
\begin{proof}
(i) is proved in \cite{samuelmolev} (it goes back to Macdonald \cite{notes}). (ii) follows from (i): a reduced subword for $u$ of a reduced word for $w$ exists only if $u\leq w$; the only subword for $u=1$ is the empty one and the only one for $u=w$ is the full word, whose composition $s_{a_1}\circ\cdots\circ s_{a_m}$ is the action of $w$; and each composition in (i) consists of $\ell(u)$ reflections, which preserve degree, and $\ell(w)-\ell(u)$ divided differences, which lower it by one. (iii) is a rearrangement of the generalized Leibniz rule of \cite[Corollary~3.5.4]{samuelleibniz}: in the notation of that work $\partial^w_u=\sum_vc^w_{u,v}\partial_v$, the coefficients $c^w_{u_0,u_1,\ldots,u_k}$ of the $k$-fold rule factor as $\sum_{v_1,\ldots,v_{k-1}}c^{v_1}_{u_0,u_1}c^{v_2}_{v_1,u_2}\cdots c^{w}_{v_{k-1},u_k}$ by associativity, and $\sum_{u_j}c^{v_j}_{v_{j-1},u_j}\partial_{u_j}(f_j)=\partial^{v_j}_{v_{j-1}}(f_j)$. By (ii) only chains $u_0\leq u_1\leq\cdots\leq u_k$ contribute.
\end{proof}

\begin{lemma}[Linear forms] \label{lemma:skewlinear}
Let $u,w\in S_\infty$ and $i\geq1$. Then
$$\partial^w_u(z_i)=\begin{cases}z_{u(i)}&\text{if }w=u,\\ 1&\text{if }w=ut_{ib}\text{ with }b>i\text{ and }\ell(w)=\ell(u)+1,\\ -1&\text{if }w=ut_{bi}\text{ with }b<i\text{ and }\ell(w)=\ell(u)+1,\\ 0&\text{otherwise.}\end{cases}$$
\end{lemma}
\begin{proof}
The case $w=u$ is Lemma \ref{lemma:skewfacts}(ii), and by the same lemma $\partial^w_u(z_i)=0$ unless $u<w$ and $\ell(w)=\ell(u)+1$. So let $w=ut_{ab}$ with $a<b$ and $\ell(w)=\ell(u)+1$, so that $u(a)<u(b)$. Fix a reduced word $(a_1,\ldots,a_m)$ for $w$ and put $\sigma_p=s_{a_{p+1}}\cdots s_{a_m}$. The reduced subwords for $u$ omit exactly one letter, and omitting the $p$-th letter gives the element $w\,\sigma_p^{-1}s_{a_p}\sigma_p=w\,t_{\sigma_p^{-1}(a_p),\,\sigma_p^{-1}(a_p+1)}$. The reflections $\sigma_p^{-1}s_{a_p}\sigma_p$, $1\leq p\leq m$, are pairwise distinct (they are the $m$ reflections $t$ with $\ell(wt)<\ell(w)$, see \cite[Section~1.6--1.7]{humphreys1992reflection}), so exactly one $p$ yields $u$, namely the one with $\{\sigma_p^{-1}(a_p),\sigma_p^{-1}(a_p+1)\}=\{a,b\}$, and by Lemma \ref{lemma:skewfacts}(i)
$$\partial^w_u(z_i)=s_{a_1}\cdots s_{a_{p-1}}\,\partial^{a_p}\,\sigma_p(z_i)=s_{a_1}\cdots s_{a_{p-1}}\,\partial^{a_p}\big(z_{\sigma_p(i)}\big),$$
which is $1$ if $\sigma_p(i)=a_p$, $-1$ if $\sigma_p(i)=a_p+1$, and $0$ otherwise. It remains to see that $\sigma_p^{-1}(a_p)=a$ and $\sigma_p^{-1}(a_p+1)=b$. The word $(a_1,\ldots,a_p)$ is reduced and ends with $a_p$, so $\tau=s_{a_1}\cdots s_{a_p}$ satisfies $\tau(a_p)>\tau(a_p+1)$; since $w=\tau\sigma_p$ this says $w(\sigma_p^{-1}(a_p))>w(\sigma_p^{-1}(a_p+1))$. As $w(a)=u(b)>u(a)=w(b)$, the claim follows.
\end{proof}

\subsection{Factorial elementary polynomials and dominant permutations}

\begin{lemma} \label{lemma:Ehe}
For $0\leq p\leq k$,
\begin{align*}
E_{p,k}(x;z)&=\sum_{m=0}^{p}(-1)^{p-m}h_{p-m}(z_1,\ldots,z_{k-p+1})\,e_m(x_1,\ldots,x_k),\\
\qelem_{p,k}(x;z)&=\sum_{m=0}^{p}(-1)^{p-m}h_{p-m}(z_1,\ldots,z_{k-p+1})\,e^{(q)}_m(x_1,\ldots,x_k).
\end{align*}
In particular $E_{p,k}(x;z)$ and $\qelem_{p,k}(x;z)$ are symmetric in $z_1,\ldots,z_{k-p+1}$, and $\qelem_{p,k}(x;z)=\mathcal Q\big(E_{p,k}(x;z)\big)$.
\end{lemma}
\begin{proof}
Lemma \ref{lemma:quantumtoordinary} with $y=0$ gives $\qsch_w(x;z)=\sum_{v_1\cdot v_2=w}\sch_{v_1}(0;z)\qsch_{v_2}(x)$. The permutation $c[p,k]=s_{k-p+1}\cdots s_k$ has a unique reduced word, so its reduced factorizations are $v_1=s_{k-p+1}\cdots s_j$, $v_2=s_{j+1}\cdots s_k$ with $k-p\leq j\leq k$. Here $v_2$ is the Grassmannian permutation with descent at $k$ whose diagram is a single column of height $m=k-j$, so $\sch_{v_2}(x)=e_m(x_1,\ldots,x_k)$ and $\qsch_{v_2}(x)=e^{(q)}_m(x_1,\ldots,x_k)$ \cite{fomin1997quantum}; and $\sch_{v_1}(0;z)=(-1)^{\ell(v_1)}\sch_{v_1^{-1}}(z)$, where $v_1^{-1}=s_j\cdots s_{k-p+1}$ is the Grassmannian permutation with descent at $k-p+1$ whose diagram is a single row of length $\ell(v_1)=p-m$, so that $\sch_{v_1^{-1}}(z)=h_{p-m}(z_1,\ldots,z_{k-p+1})$. Summing over $j$ gives the second formula, and the first is its specialization at $q=0$. The symmetry is clear from the formulas, and the last assertion holds because $\mathcal Q$ fixes $\mathbb Z[z]$ and sends the standard elementary monomial $e_m(x_1,\ldots,x_k)$ to $e^{(q)}_m(x_1,\ldots,x_k)$.
\end{proof}

\begin{lemma}[Dominant permutations] \label{lemma:dominantfactor}
Let $\mu$ be a dominant permutation with $\code(\mu)=\theta$, $\theta_1\geq\cdots\geq\theta_k>0$. Then $\mu^{-1}$ is dominant, the $i$-th column of its diagram has height $\theta_i$, and
$$\sch_{\mu^{-1}}(x;z)=\prod_{i=1}^kE_{\theta_i,\theta_i}(x;z_i)=\prod_{i=1}^k\prod_{j=1}^{\theta_i}(x_j-z_i).$$
If moreover $\theta_1>\cdots>\theta_k$, then $\qsch_{\mu^{-1}}(x;z)=\prod_{i=1}^k\qelem_{\theta_i,\theta_i}(x;z_i)$.
\end{lemma}
\begin{proof}
A permutation is dominant if and only if its diagram is the Young diagram of its code, placed in the upper left corner, and for such a permutation $\rho$ one has $\sch_\rho(x;z)=\prod_{(a,b)\in D(\rho)}(x_a-z_b)$ \cite{notes}. The diagram of $\mu^{-1}$ is the transpose of that of $\mu$, so $D(\mu^{-1})=\{(j,i)\mid 1\leq i\leq k,\ 1\leq j\leq\theta_i\}$, which is again a Young diagram in the upper left corner; this gives the first assertions and the factorization. For the quantum statement, expand each factor by Lemma \ref{lemma:qelemkk} at $q=0$: $E_{\theta_i,\theta_i}(x;z_i)=\sum_{m}(-z_i)^{\theta_i-m}e_m(x_1,\ldots,x_{\theta_i})$. Multiplying out, each term $\prod_ie_{m_i}(x_1,\ldots,x_{\theta_i})$ is a standard elementary monomial because the alphabets $x_1,\ldots,x_{\theta_i}$ have pairwise distinct lengths, and $\mathcal Q$ sends it to $\prod_ie^{(q)}_{m_i}(x_1,\ldots,x_{\theta_i})$. Hence $\qsch_{\mu^{-1}}(x;z)=\mathcal Q(\sch_{\mu^{-1}}(x;z))=\prod_i\sum_m(-z_i)^{\theta_i-m}e^{(q)}_m(x_1,\ldots,x_{\theta_i})=\prod_i\qelem_{\theta_i,\theta_i}(x;z_i)$ by Lemma \ref{lemma:qelemkk}.
\end{proof}

\subsection{The skew divided differences of $E_{k,k}(x;z_i)$}

For $i\geq1$ and an integer $b\neq i$ put $\epsilon_i(b)=1$ if $b<i$ and $\epsilon_i(b)=-1$ if $b>i$. If $u\Tom{i}w$ with chain $b_1,\ldots,b_k$ as in Definition \ref{def:itom}, then $\prod_j\epsilon_i(b_j)=(-1)^{\ell(u,w)-\sigma_i(u,w)}$. The heart of the matter is the following lemma, which governs the cancellation of the chains through position $i$ that revisit a position.

\begin{lemma}[Rotation] \label{lemma:rotation}
Let $\psi\in S_\infty$, $i\geq1$, and let $\gamma_1,\ldots,\gamma_m$ ($m\geq2$) be distinct integers different from $i$, with indices read modulo $m$. For $1\leq r\leq m$ let
$$R_r=(\gamma_r,\gamma_{r+1},\ldots,\gamma_{r+m-1},\gamma_r),$$
a sequence of length $m+1$, and call $R_r$ \emph{valid} if the sequence of permutations $\psi,\ \psi t_{i\gamma_r},\ \psi t_{i\gamma_r}t_{i\gamma_{r+1}},\ \ldots$ obtained by applying its entries in turn is a chain of Bruhat covers. Then:
\begin{enumerate}
\item[(i)] $t_{i\gamma_r}t_{i\gamma_{r+1}}\cdots t_{i\gamma_{r+m-1}}t_{i\gamma_r}=C$ for every $r$, where $C$ is the cycle $\gamma_j\mapsto\gamma_{j-1}$ ($1\leq j\leq m$), which fixes $i$. In particular all valid $R_r$ lead from $\psi$ to the same permutation $\psi C$.
\item[(ii)] Listed in cyclic order of $r$, the valid $R_r$ alternate between those with $\gamma_r<i$ and those with $\gamma_r>i$. In particular the number of valid $R_r$ with $\gamma_r<i$ equals the number of valid $R_r$ with $\gamma_r>i$.
\end{enumerate}
\end{lemma}
\begin{proof}
(i) Since $t_{ib}t_{ic}t_{ib}=t_{bc}$ for distinct $b,c\neq i$, the product equals $t_{\gamma_r\gamma_{r+1}}t_{\gamma_r\gamma_{r+2}}\cdots t_{\gamma_r\gamma_{r+m-1}}$, which in the notation \eqref{eq:pivot} is the cycle $(\gamma_{r+m-1},\ldots,\gamma_{r+1},\gamma_r)$, sending $\gamma_{r+m-1}\mapsto\gamma_{r+m-2}\mapsto\cdots\mapsto\gamma_{r+1}\mapsto\gamma_r\mapsto\gamma_{r+m-1}$.

(ii) Put $c=\psi(i)$ and $v_j=\psi(\gamma_j)$. Fix $r$ and follow $R_r$. Before the $s$-th step ($1\leq s\leq m$) the position $i$ carries $v_{r+s-2}$ (that is, $c$ when $s=1$), and the position $\gamma_{r+s-1}$, which has not been moved yet, carries $v_{r+s-1}$; before the last step the position $i$ carries $v_{r-1}$ and the position $\gamma_r$ carries $c$. Call a step exchanging the position $i$ with a position $g$ \emph{ascending} if the value at $g$ exceeds the value at $i$ when $g>i$, and is smaller than it when $g<i$. By Lemma \ref{lemma:lengthchange} an ascending step increases the length by at least $1$, a non-ascending step decreases it, and a step is a Bruhat cover exactly when it is ascending and increases the length by exactly $1$. The total change in length along $R_r$ is $\ell(\psi C)-\ell(\psi)$ by (i), independently of $r$. Consequently $R_r$ is valid if and only if all of its $m+1$ steps are ascending and $\ell(\psi C)-\ell(\psi)=m+1$.

By the description of the values before each step, the steps $2,\ldots,m$ of $R_r$ are all ascending if and only if
\begin{equation}\label{eq:DC}
v_j>v_{j-1}\iff\gamma_j>i\qquad\text{for all }j\neq r,
\end{equation}
while the first step is ascending if and only if ($v_r>c\iff\gamma_r>i$) and the last if and only if ($c>v_{r-1}\iff\gamma_r>i$). The last two conditions together say precisely that $c$ lies strictly between $v_{r-1}$ and $v_r$ and that \eqref{eq:DC} holds for $j=r$ as well. Therefore $R_r$ is valid if and only if
$$\ell(\psi C)-\ell(\psi)=m+1,\qquad\text{\eqref{eq:DC} holds for all }j,\qquad\text{and }c\text{ lies strictly between }v_{r-1}\text{ and }v_r.$$
The first two conditions do not involve $r$. If they hold, the valid indices $r$ are exactly those at which the cyclic sequence $v_1,\ldots,v_m$, none of whose terms equals $c$, crosses the level $c$; and consecutive crossings of a cyclic sequence are alternately upward ($v_{r-1}<c<v_r$) and downward ($v_{r-1}>c>v_r$). By \eqref{eq:DC} the upward crossings are those with $\gamma_r>i$ and the downward ones those with $\gamma_r<i$.
\end{proof}

\begin{theorem} \label{theorem:skewtop}
Let $u,w\in S_\infty$, $i\geq1$ and $k\geq0$. Then $\partial^w_u\big(E_{k,k}(x;z_i)\big)=0$ unless $u\Tom{i}w$, in which case
$$\partial^w_u\big(E_{k,k}(x;z_i)\big)=(-1)^{\ell(u,w)-\sigma_i(u,w)}\,E_{k-\ell(u,w),\,k}\big(x;z_{D_i(u,w)}\big).$$
\end{theorem}
\begin{proof}
Write $E_{k,k}(x;z_i)=f_1\cdots f_k$ with $f_j=x_j-z_i$. By Lemma \ref{lemma:skewfacts}(iii),
$$\partial^w_u\big(E_{k,k}(x;z_i)\big)=\sum_{u=u_0\leq\cdots\leq u_k=w}\ \prod_{j=1}^k\partial^{u_j}_{u_{j-1}}(x_j-z_i),$$
and by Lemmas \ref{lemma:skewlinear} and \ref{lemma:skewfacts}(ii) the factor $\partial^{u_j}_{u_{j-1}}(x_j-z_i)$ equals $x_j-z_{u_j(i)}$ if $u_j=u_{j-1}$, equals $\epsilon_i(b)$ if $u_j=u_{j-1}t_{ib}$ with $\ell(u_j)=\ell(u_{j-1})+1$, and vanishes otherwise. Call a sequence $(u_0,\ldots,u_k)$ in which every step is of one of these two kinds an \emph{$i$-walk} from $u$ to $w$. Thus the sum is over the $i$-walks, and a walk contributes the product of $x_j-z_{u_j(i)}$ over its trivial steps and of $\epsilon_i(b)$ over its nontrivial steps.

An \emph{$i$-chain} from $u$ to $w$ is a sequence $\beta=(\beta_1,\ldots,\beta_\ell)$ of integers different from $i$ such that, with $u^{(0)}=u$ and $u^{(s)}=u^{(s-1)}t_{i\beta_s}$, every step $u^{(s-1)}\lessdot u^{(s)}$ is a Bruhat cover and $u^{(\ell)}=w$; repetitions among the $\beta_s$ are allowed. Put $c_s(\beta)=u^{(s)}(i)$, the value at position $i$ after $s$ steps, and $\epsilon(\beta)=\prod_s\epsilon_i(\beta_s)$. Every $i$-walk determines the $i$-chain formed by its nontrivial steps, and conversely an $i$-chain of length $\ell\leq k$ arises from exactly the $\binom{k}{\ell}$ walks obtained by choosing the positions $j_1<\cdots<j_p$ ($p=k-\ell$) of the trivial steps. In such a walk the trivial step at $j_r$ is preceded by $j_r-r$ nontrivial steps, so it contributes the factor $x_{j_r}-z_{c_{j_r-r}(\beta)}$. Summing over the choices of $j_1<\cdots<j_p$ and comparing with the explicit formula for $E_{p,k}$ in Definition \ref{definition:pieri}, the walks with chain $\beta$ contribute in total
$$\epsilon(\beta)\sum_{1\leq j_1<\cdots<j_p\leq k}\ \prod_{r=1}^p\big(x_{j_r}-z_{c_{j_r-r}(\beta)}\big)=\epsilon(\beta)\,E_{k-\ell,k}\big(x;z_{c_0(\beta)},z_{c_1(\beta)},\ldots,z_{c_\ell(\beta)}\big).$$
Hence
\begin{equation}\label{eq:chainsum}
\partial^w_u\big(E_{k,k}(x;z_i)\big)=\sum_{\beta}\epsilon(\beta)\,E_{k-\ell(\beta),k}\big(x;z_{c_0(\beta)},\ldots,z_{c_{\ell(\beta)}(\beta)}\big),
\end{equation}
the sum being over all $i$-chains $\beta$ from $u$ to $w$ of length $\ell(\beta)\leq k$. By Lemma \ref{lemma:Ehe} the summand depends on $c_0(\beta),\ldots,c_\ell(\beta)$ only through their multiset.

\emph{Chains without repetitions.} If $\beta_1,\ldots,\beta_\ell$ are distinct, then $\beta$ is a witness for $u\Tom{i}w$, and by Definition \ref{def:itom} it is the unique one; moreover $\ell=\ell(u,w)$, $\epsilon(\beta)=(-1)^{\ell(u,w)-\sigma_i(u,w)}$, and since the position $\beta_s$ is untouched before the $s$-th step, $c_s(\beta)=u(\beta_s)$, so that $\{c_0(\beta),\ldots,c_\ell(\beta)\}=D_i(u,w)$. So the chains without repetitions contribute exactly the right-hand side of the theorem if $u\Tom{i}w$ and $\ell(u,w)\leq k$, and nothing otherwise; if $u\Tom{i}w$ and $\ell(u,w)>k$ then both sides of the theorem vanish, the left by degree. It remains to show that the chains with a repetition contribute zero in total.

\emph{Chains with repetitions.} Let $\beta$ be an $i$-chain from $u$ to $w$ with a repetition, and let $s$ be the least index for which $\beta_s\in\{\beta_1,\ldots,\beta_{s-1}\}$, say $\beta_s=\beta_t$ with $t<s$. Put $m=s-t$ and $(\gamma_1,\ldots,\gamma_m)=(\beta_t,\ldots,\beta_{s-1})$; these are distinct, and $m\geq2$ because two consecutive covers cannot use the same transposition. The segment $(\beta_t,\ldots,\beta_s)$ of $\beta$ is the sequence $R_1$ of Lemma \ref{lemma:rotation} for $\psi=u^{(t-1)}$, and it is valid. The alternation in Lemma \ref{lemma:rotation}(ii) pairs each valid $R_r$ with $\gamma_r>i$ with the next valid sequence in cyclic order, whose entry lies below $i$, and each valid $R_r$ with $\gamma_r<i$ with the previous one; this is a fixed-point-free involution on the valid sequences which exchanges the two sides of $i$. Let $R_{r^*}$ be the sequence paired with $R_1$, and let $\beta^*$ be obtained from $\beta$ by replacing the segment $(\beta_t,\ldots,\beta_s)$ with $R_{r^*}$. By Lemma \ref{lemma:rotation}(i) the new segment again leads from $u^{(t-1)}$ to $u^{(t-1)}C=u^{(s)}$, so $\beta^*$ is an $i$-chain from $u$ to $w$ of the same length as $\beta$. The entries of $\beta^*$ preceding the index $s$ form the same set as those of $\beta$ and are distinct, and $\beta^*_s=\beta^*_t$, so the first repetition of $\beta^*$ occurs at the same index $s$ with the same $t$, the same $\psi$, and the same cyclic sequence $\gamma$; hence $\beta^{**}=\beta$, and $\beta\mapsto\beta^*$ is a fixed-point-free involution on the set of $i$-chains with repetitions. The values $c_0,\ldots,c_{t-1}$ and $c_s,\ldots,c_\ell$ are the same for $\beta$ and $\beta^*$, since the permutations $u^{(0)},\ldots,u^{(t-1)}$ and $u^{(s)},\ldots,u^{(\ell)}$ are; and $c_t,\ldots,c_{s-1}$ are the values $\psi(\gamma_1),\ldots,\psi(\gamma_m)$ in some order for both chains. So the multisets of values coincide, and the $E$-factors in \eqref{eq:chainsum} agree. Finally $\epsilon(\beta^*)=-\epsilon(\beta)$: the two segments consist of the entries $\gamma_1,\ldots,\gamma_m$, each once, together with one additional entry, which is $\gamma_1$ for $\beta$ and $\gamma_{r^*}$ for $\beta^*$, and these lie on opposite sides of $i$. Hence the chains with repetitions cancel in pairs in \eqref{eq:chainsum}, and the theorem follows.
\end{proof}

\subsection{Conclusion of the proofs}

\begin{proof}[Proof of Theorem \ref{theorem:orddouble}]
From the definition of the double Schubert polynomials, for every $w$ and $a$,
\begin{equation}\label{eq:leftdd}
\partial^a\sch_w(x;z)=\begin{cases}-\sch_{s_aw}(x;z)&\text{if }s_aw<w,\\ 0&\text{if }s_aw>w:\end{cases}
\end{equation}
indeed $\sch_w=(-1)^{\ell(w_0)-\ell(w)}\partial^{ww_0}\sch_{w_0}$, and $\partial^a\partial^{ww_0}$ equals $\partial^{s_aww_0}$ if $\ell(s_aww_0)>\ell(ww_0)$, which is the case $s_aw<w$, and equals $0$ otherwise. Since $\mu=v^{-1}\cdot y$ is reduced, so is $\mu^{-1}=y^{-1}\cdot v$; writing $y=s_{a_1}\cdots s_{a_r}$ reduced, the permutations $s_{a_j}s_{a_{j+1}}\cdots s_{a_r}\mu^{-1}$, $j=r+1,r,\ldots,1$, therefore decrease in length by one at each step, from $\mu^{-1}$ down to $y\mu^{-1}=v$. Applying \eqref{eq:leftdd} $r$ times and then Lemma \ref{lemma:dominantfactor},
$$\sch_v(x;z)=(-1)^{\ell(y)}\,\partial^{y}\sch_{\mu^{-1}}(x;z)=(-1)^{\ell(y)}\,\partial^{y}\Big(\prod_{i=1}^kE_{\theta_i,\theta_i}(x;z_i)\Big).$$
Now expand by the ascending Leibniz formula, Lemma \ref{lemma:skewfacts}(iii) with $u=1$ and $w=y$, and evaluate each factor by Theorem \ref{theorem:skewtop}:
\begin{align*}
\partial^{y}\Big(\prod_{i=1}^kE_{\theta_i,\theta_i}(x;z_i)\Big)&=\sum_{1=v_0\leq\cdots\leq v_k=y}\ \prod_{i=1}^k\partial^{v_i}_{v_{i-1}}\big(E_{\theta_i,\theta_i}(x;z_i)\big)\\
&=\sum_{1=v_0\Tom{1}v_1\Tom{2}\cdots\Tom{k}v_k=y}\ \prod_{i=1}^k(-1)^{\ell(v_{i-1},v_i)-\sigma_i(v_{i-1},v_i)}E_{\theta_i-\ell(v_{i-1},v_i),\,\theta_i}\big(x;z_{D_i(v_{i-1},v_i)}\big).
\end{align*}
Since $\sum_i\ell(v_{i-1},v_i)=\ell(y)$, the signs combine to $\prod_i(-1)^{\sigma_i(v_{i-1},v_i)}$; and $y=v\mu$.
\end{proof}

\begin{proof}[Proof of Theorem \ref{theorem:quandouble}]
A strict dominant permutation is dominant, so Theorem \ref{theorem:orddouble} applies to $v$ and $\mu$ and expresses $\sch_v(x;z)$ as the sum, over the chains $1=v_0\Tom{1}\cdots\Tom{k}v_k=v\mu$, of $\prod_i(-1)^{\sigma_i}E_{\thop_i-\ell_i,\thop_i}(x;z_{D_i})$, where $\ell_i=\ell(v_{i-1},v_i)$, $\sigma_i=\sigma_i(v_{i-1},v_i)$ and $D_i=D_i(v_{i-1},v_i)$. Apply $\mathcal Q$ to both sides. The left-hand side becomes $\qsch_v(x;z)$. On the right, by Lemma \ref{lemma:Ehe} each factor is $E_{\thop_i-\ell_i,\thop_i}(x;z_{D_i})=\sum_{m}c_{i,m}\,e_m(x_1,\ldots,x_{\thop_i})$ with coefficients $c_{i,m}\in\mathbb Z[z]$, and multiplying out, each term $\prod_ie_{m_i}(x_1,\ldots,x_{\thop_i})$ of the product over $i$ is a standard elementary monomial because $\thop_1>\cdots>\thop_k$. Hence $\mathcal Q$ sends the product to $\prod_i\sum_mc_{i,m}\,e^{(q)}_m(x_1,\ldots,x_{\thop_i})=\prod_i\qelem_{\thop_i-\ell_i,\thop_i}(x;z_{D_i})$, again by Lemma \ref{lemma:Ehe}.
\end{proof}

We prove Theorem \ref{theorem:quandoublemed} by reaching the medium shape from the strict shape by divided differences, one for each retained pair; this also explains the form of the doubled factor. For $k\geq1$ and $i\geq1$ write $\qsch_{d_2[k]}(x;z_i,z_{i+1})$ for the quantum double Schubert polynomial of the $k\times2$ rectangle (Section \ref{section:additional}) with $(z_1,z_2)$ replaced by $(z_i,z_{i+1})$, i.e.\ $\qsch_{d_2[k]}(x;z_i,z_{i+1})=\mathcal Q\big(E_{k,k}(x;z_i)E_{k,k}(x;z_{i+1})\big)=\qelem_{k,k;k}(x;z_i,z_{i+1})$.

\begin{corollary}[Factorization for the medium shape] \label{corollary:mediumfactor}
Let $\mu$ be a dominant permutation of medium shape with $\medtheta=\code(\mu)$, paired indices $\mathcal P$ and single indices $\mathcal S$ as in Theorem \ref{theorem:quandoublemed}, and let $\nu=\mu\prod_{i\in\mathcal P}s_i$. Then $\nu$ is strict dominant, $\code(\nu)$ is $\medtheta$ with the first entry of each retained pair increased by one, $\mu=\nu\prod_{i\in\mathcal P}s_i$ with $\ell(\mu)=\ell(\nu)-|\mathcal P|$, and
$$\qsch_{\mu^{-1}}(x;z)=(-1)^{|\mathcal P|}\Big(\prod_{i\in\mathcal P}\partial^{i}\Big)\qsch_{\nu^{-1}}(x;z)=\prod_{i\in\mathcal S}\qelem_{\medtheta_i,\medtheta_i}(x;z_i)\ \prod_{i\in\mathcal P}\qsch_{d_2[\medtheta_i]}(x;z_i,z_{i+1}).$$
In particular, for $\mu=\medperm(v^{-1})$ one has $\nu=\mom(v^{-1})$ by Proposition \ref{prop:weakorder}(iii).
\end{corollary}
\begin{proof}
For $i\in\mathcal P$ with $k=\medtheta_i=\medtheta_{i+1}$ we have $\mu(i)<\mu(i+1)$, so $\ell(\mu s_i)=\ell(\mu)+1$ and $\code(\mu s_i)$ is $\code(\mu)$ with $(\code_i,\code_{i+1})=(k,k)$ replaced by $(k+1,k)$; the indices in $\mathcal P$ are pairwise at distance $\geq2$, so the $s_i$ commute and act on disjoint pairs of positions, which gives the assertions about $\nu$ (its code is strictly decreasing on its nonzero part because the entry preceding a retained pair is $\geq k+2$ and the entry following it is $\leq k-1$). Enumerate $\mathcal P=\{i_1<\cdots<i_r\}$ and put $\mu_0=\mu$, $\mu_t=\mu_{t-1}s_{i_t}$, so $\mu_r=\nu$. Since $s_{i_t}\mu_t^{-1}=(\mu_ts_{i_t})^{-1}=\mu_{t-1}^{-1}$ is shorter than $\mu_t^{-1}$, \eqref{eq:leftdd} gives $\sch_{\mu_{t-1}^{-1}}(x;z)=-\partial^{i_t}\sch_{\mu_t^{-1}}(x;z)$, and iterating, $\sch_{\mu^{-1}}(x;z)=(-1)^r\partial^{i_1}\cdots\partial^{i_r}\sch_{\nu^{-1}}(x;z)$. The map $\mathcal Q$ acts on the $x$-variables and is $\mathbb Z[q][z]$-linear, while the $\partial^a$ act on the $z$-variables and are $\mathbb Z[q][x]$-linear; hence they commute, and applying $\mathcal Q$ gives the first equality. For the second, Lemma \ref{lemma:dominantfactor} gives $\qsch_{\nu^{-1}}(x;z)=\prod_j\qelem_{\code_j(\nu),\code_j(\nu)}(x;z_j)$, in which the factor of index $j\notin\mathcal P\cup(\mathcal P+1)$ is $\qelem_{\medtheta_j,\medtheta_j}(x;z_j)$ and the factors of indices $i,i+1$ for $i\in\mathcal P$, $k=\medtheta_i$, are $\qelem_{k+1,k+1}(x;z_i)\qelem_{k,k}(x;z_{i+1})$. The operator $\partial^i$ involves only $z_i,z_{i+1}$, so it acts on this pair of factors alone, and the operators for distinct $i\in\mathcal P$ act on disjoint sets of variables. Finally, with $P_m(t)=\prod_{j\leq m}(x_j-t)$,
$$P_{k+1}(z_i)P_k(z_{i+1})-P_{k+1}(z_{i+1})P_k(z_i)=P_k(z_i)P_k(z_{i+1})\big[(x_{k+1}-z_i)-(x_{k+1}-z_{i+1})\big]=(z_{i+1}-z_i)\,P_k(z_i)P_k(z_{i+1}),$$
so $-\partial^i\big(P_{k+1}(z_i)P_k(z_{i+1})\big)=P_k(z_i)P_k(z_{i+1})$. Applying $\mathcal Q$, which is multiplicative on $P_{k+1}(z_i)P_k(z_{i+1})$ (distinct alphabets, Lemma \ref{lemma:dominantfactor}) and commutes with $\partial^i$, gives $-\partial^i\big(\qelem_{k+1,k+1}(x;z_i)\qelem_{k,k}(x;z_{i+1})\big)=\mathcal Q\big(P_k(z_i)P_k(z_{i+1})\big)=\qsch_{d_2[k]}(x;z_i,z_{i+1})$.
\end{proof}

\begin{proof}[Proof of Theorem \ref{theorem:quandoublemed}]
Exactly as in the proof of Theorem \ref{theorem:orddouble}, $\qsch_v(x;z)=(-1)^{\ell(y)}\partial^y\qsch_{\mu^{-1}}(x;z)$ (the identity \eqref{eq:leftdd} holds for the quantum polynomials, $\mathcal Q$ commuting with the $z$-divided differences). Insert the factorization of Corollary \ref{corollary:mediumfactor} and expand by the ascending Leibniz formula, Lemma \ref{lemma:skewfacts}(iii), treating each doubled factor as a single factor. The skew operators $\partial^w_u$ are $\mathbb Z[q][x]$-linear and so commute with $\mathcal Q$; hence, by Lemma \ref{lemma:skewfacts}(iii) once more and Theorem \ref{theorem:skewtop},
\begin{multline*}
\partial^{w}_{u}\,\qsch_{d_2[k]}(x;z_i,z_{i+1})=\mathcal Q\Big(\sum_{u\leq u'\leq w}\partial^{u'}_u\big(E_{k,k}(x;z_i)\big)\,\partial^{w}_{u'}\big(E_{k,k}(x;z_{i+1})\big)\Big)\\
=\sum_{u\Tom iu'\Tom{i+1}w}(-1)^{\ell_i+\ell_{i+1}-\sigma_i-\sigma_{i+1}}\,\qelem_{k-\ell_i,\,k-\ell_{i+1};\,k}\big(x;z_{D_i},z_{D_{i+1}}\big),
\end{multline*}
with $\ell_i,\sigma_i,D_i$ attached to $u\Tom iu'$ and $\ell_{i+1},\sigma_{i+1},D_{i+1}$ to $u'\Tom{i+1}w$, while the single factors are evaluated by Theorem \ref{theorem:skewtop} and Lemma \ref{lemma:Ehe} as in the proof of Theorem \ref{theorem:quandouble}. Collecting the signs as in the proof of Theorem \ref{theorem:orddouble} gives the theorem. In this form it is literally Theorem \ref{theorem:quandouble} for $\nu=\mu\prod_{i\in\mathcal P}s_i$ with the divided differences $\partial^i$, $i\in\mathcal P$, applied first: the two Pieri steps at columns $i,i+1$ of the strict shape, in $k+1$ and $k$ variables, are replaced by the doubled step in $k$ variables.
\end{proof}

\section{General-degree and doubled Pieri formulas} \label{section:additional}

In this section we record the two forms of the Pieri rule that the quantum multiplication algorithm of Section \ref{section:algorithms} actually uses: the rule for a quantum factorial elementary polynomial of arbitrary degree $p\leq k$, and a ``doubled'' rule for the quantization of a product of two such polynomials in the \emph{same} $k$ variables. Recall from Definition \ref{definition:pieri} that $n_k(u,w)=k-|\fx_k(u,w)|$ is the number of positions $i\leq k$ at which $u$ and $w$ differ, equivalently the number of transpositions in a chain witnessing $u\tom{k}w$, and that $E_{p,k}(y_A;z)$ denotes the factorial elementary polynomial in the variables $y_a$, $a\in A$ (with $|A|=k$), and $z_1,z_2,\ldots$.

\begin{theorem}[Quantum Pieri formula, general degree] \label{theorem:pierigeneral}
Suppose $u\in S_\infty$ and $0\leq p\leq k$. Then
$$\qsch_u(x;y)\,\qelem_{p,k}(x;z)=\sum_{\substack{u\tom{k} w\\ n_k(u,w)\leq p}}q_k(u,w)\,E_{p-n_k(u,w),\,k-n_k(u,w)}\big(y_{\fx_k(u,w)};z\big)\,\qsch_w(x;y).$$
\end{theorem}

For $p=k$ the coefficient is $q_k(u,w)\wte(y_{\fx_k(u,w)};z)$ and this is Theorem \ref{theorem:pieri}; at $q=0$ it is the Pieri formula of \cite[Theorem~7.1]{samuelmolev}. This is the form implemented in \cite{schubmult} (the routines \texttt{elem\_sym\_perms\_q} and \texttt{elem\_sym\_func\_q}). We deduce the general case from the case $p=k$.

\begin{proof}
Let $X=(X_1,\ldots,X_k)$ be indeterminates. The factorial elementary polynomial
$$E_{p,k}(X;z)=\sum_{1\leq i_1<\cdots<i_p\leq k}\ \prod_{j=1}^p(X_{i_j}-z_{i_j-j+1})$$
is symmetric in $X$ and of degree at most one in each $X_i$, hence is a $\mathbb Z[z]$-linear combination of the elementary symmetric polynomials:
\begin{equation}\label{eq:Eexpansion}
E_{p,k}(X;z)=\sum_{j=0}^{p}c_{p,k,j}(z)\,e_j(X_1,\ldots,X_k),\qquad c_{p,k,j}(z)\in\mathbb Z[z].
\end{equation}
Two consequences. First, applying the quantization map $\mathcal Q$ (which is $\mathbb Z[z]$-linear and sends $e_j(x_1,\ldots,x_k)$ to $e^{(q)}_j(x_1,\ldots,x_k)$) gives
\begin{equation}\label{eq:qEexpansion}
\qelem_{p,k}(x;z)=\sum_{j=0}^{p}c_{p,k,j}(z)\,e^{(q)}_j(x_1,\ldots,x_k).
\end{equation}
Second, for $0\leq n\leq k$ extract from \eqref{eq:Eexpansion} the coefficient of the monomial $X_{k-n+1}X_{k-n+2}\cdots X_k$. On the right it is $\sum_jc_{p,k,j}(z)\,e_{j-n}(X_1,\ldots,X_{k-n})$. On the left, a term of $E_{p,k}(X;z)$ contains this monomial exactly when $\{i_{p-n+1},\ldots,i_p\}=\{k-n+1,\ldots,k\}$, and then the remaining factors $\prod_{j\leq p-n}(X_{i_j}-z_{i_j-j+1})$ with $i_1<\cdots<i_{p-n}\leq k-n$ are precisely the terms of $E_{p-n,k-n}(X_1,\ldots,X_{k-n};z)$. Hence, for every set $A$ of $k-n$ indices,
\begin{equation}\label{eq:Ereduction}
\sum_{j=0}^{p}c_{p,k,j}(z)\,e_{j-n}(y_A)=E_{p-n,\,k-n}(y_A;z),
\end{equation}
both sides being $0$ when $n>p$.

Now put $z_1=-t$ in Theorem \ref{theorem:pieri}. By Lemma \ref{lemma:qelemkk} the left side becomes $\sum_{j=0}^k t^{k-j}\,\qsch_u(x;y)\,e^{(q)}_j(x_1,\ldots,x_k)$, while on the right $\wte(y_{\fx_k(u,w)};z)=\prod_{a\in\fx_k(u,w)}(y_a+t)=\sum_i t^{\,k-n_k(u,w)-i}e_i(y_{\fx_k(u,w)})$, since $|\fx_k(u,w)|=k-n_k(u,w)$. Comparing the coefficients of $t^{k-j}$ gives, for $0\leq j\leq k$,
\begin{equation}\label{eq:nonfactorialpieri}
\qsch_u(x;y)\,e^{(q)}_j(x_1,\ldots,x_k)=\sum_{u\tom{k}w}q_k(u,w)\,e_{j-n_k(u,w)}\big(y_{\fx_k(u,w)}\big)\,\qsch_w(x;y).
\end{equation}
Multiplying \eqref{eq:nonfactorialpieri} by $c_{p,k,j}(z)$, summing over $j$ and using \eqref{eq:qEexpansion} and then \eqref{eq:Ereduction} with $A=\fx_k(u,w)$ and $n=n_k(u,w)$ yields the theorem.
\end{proof}

\subsection*{The doubled formula}

Recall that the quantization map $\mathcal Q$ of Fomin, Gelfand and Postnikov (Lemma \ref{lemma:quantumtoordinary}) is linear over the coefficient variables, sends each standard elementary monomial $\prod_j e_{i_j}(x_1,\ldots,x_j)$ to $\prod_j e^{(q)}_{i_j}(x_1,\ldots,x_j)$, and satisfies $\mathcal Q(\sch_w(x;z))=\qsch_w(x;z)$. In particular $\mathcal Q(E_{p,k}(x;z))=\qelem_{p,k}(x;z)$, and more generally $\mathcal Q$ is multiplicative on a product of factorial elementary polynomials $E_{p_j,k_j}(x;z^{(j)})$, $j=1,2,\ldots$, whose numbers of variables $k_j$ are \emph{pairwise distinct}, since such a product expands into standard elementary monomials. It is \emph{not} multiplicative when two factors have the same number of variables: in general$$\mathcal Q\big(E_{p_1,k}(x;z)E_{p_2,k}(x;z')\big)\neq\qelem_{p_1,k}(x;z)\,\qelem_{p_2,k}(x;z').$$ The following definition isolates the terms of the iterated Pieri formula \eqref{eq:allk} below which survive in the quantization of the product.

\begin{definition}[Compatible pairs of chains] \label{def:compatible}
Let $u\tom{k}v\tom{k}w$. We say that the pair of chains is \emph{compatible} if for every nontrivial cycle $d$ of $v^{-1}w$ and every nontrivial cycle $c$ of $u^{-1}v$ with $\max c=\max d$, one has $c\cap d=\{\max d\}$. In words: a cycle of the second step which shares its top with a cycle of the first step may meet that cycle only in the top. (A cycle of $v^{-1}w$ whose top is not the top of any cycle of $u^{-1}v$ is unrestricted.)
\end{definition}

Throughout this subsection we write $\qelem_{k,k}(x;z_1)$ and $\qelem_{k,k}(x;z_2)$ for the quantum factorial elementary polynomial $\qelem_{k,k}(x;z)=\mathcal Q\big(\prod_{i\leq k}(x_i-z_1)\big)$ and for the same polynomial with $z_1$ replaced by $z_2$; thus $\qelem_{k,k}(x;z_2)$ is $\qelem_{k,k}(x;z')$ for the shifted alphabet $z'=(z_2,z_3,\ldots)$. For a pair of chains $u\tom{k}v\tom{k}w$ put
$$\omega_k(u,v,w)=q_k(u,v)\,q_k(v,w)\prod_{a\in\fx_k(u,v)}(y_a-z_1)\prod_{b\in\fx_k(v,w)}(y_b-z_2),$$
so that two applications of Theorem \ref{theorem:pieri}, the second with $z_1$ replaced by $z_2$, give
\begin{equation}\label{eq:allk}
\qsch_u(x;y)\,\qelem_{k,k}(x;z_1)\,\qelem_{k,k}(x;z_2)=\sum_{u\tom{k}v\tom{k}w}\omega_k(u,v,w)\,\qsch_w(x;y).
\end{equation}
Recall that $d_2[k]$ denotes the dominant permutation with $\code(d_2[k])=(2^k)$, whose diagram is a $k\times2$ rectangle, so that
$$\sch_{d_2[k]}(x;z)=\prod_{i\leq k}(x_i-z_1)(x_i-z_2)=E_{k,k}(x;z_1)E_{k,k}(x;z_2),\qquad \qsch_{d_2[k]}(x;z)=\mathcal Q\big(E_{k,k}(x;z_1)E_{k,k}(x;z_2)\big);$$
we set $d_2[0]=\mathrm{id}$.

\begin{theorem}[Doubled quantum Pieri formula] \label{theorem:doubledpieri}
For $u\in S_\infty$ and $k\geq1$,
\begin{multline*}
\qsch_u(x;y)\,\qsch_{d_2[k]}(x;z)=\sum_{\substack{u\tom{k}v\tom{k}w\\ \text{compatible}}}\omega_k(u,v,w)\,\qsch_w(x;y)\\
=\sum_{\substack{u\tom{k}v\tom{k}w\\ \text{compatible}}}q_k(u,v)q_k(v,w)\prod_{a\in\fx_k(u,v)}(y_a-z_1)\prod_{b\in\fx_k(v,w)}(y_b-z_2)\,\qsch_w(x;y).
\end{multline*}
\end{theorem}

Thus the doubled formula is obtained from the iterated formula \eqref{eq:allk} by \emph{omitting} the incompatible pairs, exactly as the quantization $\qsch_{d_2[k]}$ of the product $E_{k,k}(x;z_1)E_{k,k}(x;z_2)$ omits terms from the product $\qelem_{k,k}(x;z_1)\qelem_{k,k}(x;z_2)$ of the quantizations. This is the rule implemented by the routine \texttt{double\_elem\_sym\_q} of \cite{schubmult}, and it is the ``paired step'' of the quantum algorithm of Section \ref{section:algorithms}. We have verified the theorem by computer for all $u\in S_4$ and all $k\leq3$; in these cases the unrestricted sum \eqref{eq:allk} fails to equal the left-hand side in roughly half of the instances, so the compatibility condition is essential.

\begin{example}[Compatible and incompatible pairs] \label{example:doubled}
(a) $k=1$. A step $u\tom1v$ is either trivial or a single cover $t_{1b}$, so a cycle of $u^{-1}v$ and a cycle of $v^{-1}w$ with a common top $B$ sharing the element $1$ are both equal to $(1,B)$, i.e.\ $v=ut_{1B}$ and $w=u$; and $t_{1B}$ can be a cover of both $u$ and $ut_{1B}$ only when $B=2$ (if $B\geq3$, the classical cover among the two requires no intermediate value in the interval between $u(1)$ and $u(B)$, while the quantum cover requires every intermediate value to lie in it). Thus the incompatible pairs are exactly the round trips $u\tom1us_1\tom1u$, in which one step is a classical and the other a quantum cover, both fixed sets are empty, and $\omega_1(u,us_1,u)=q_1$. Theorem \ref{theorem:doubledpieri} for $k=1$ therefore says that $\qsch_u\,\qsch_{312}(x;z)$ is obtained from $\qsch_u\,\qelem_{1,1}(x;z_1)\,\qelem_{1,1}(x;z_2)$ by deleting the term $q_1\qsch_u$, in agreement with $\qsch_{312}(x;z)=\mathcal Q\big((x_1-z_1)(x_1-z_2)\big)=(x_1-z_1)(x_1-z_2)-q_1$, which is the case $k=1$ of the recurrence of Lemma \ref{lemma:thetaPP} below. For $u=21$ the seven compatible pairs and the rejected round trip $21\tom112\tom121$ are listed in Example \ref{example:algorithm}.

(b) $k=2$, $u=1324$. The steps $u\tom2v$ are: the trivial one; the single covers $t_{13}$ (classical, $v=2314$), $t_{23}$ (quantum, weight $q_2$, $v=1234$) and $t_{24}$ (classical, $v=1423$); and the two-step chains $t_{24}t_{13}$ ($v=2413$) and $t_{13}t_{23}$ (weight $q_2$, $v=2134$). Consider first $w=3124$. Exactly two pairs $u\tom2v\tom2w$ exist, namely
$$1324\xrightarrow{t_{13}}2314\xrightarrow{t_{23}t_{13}}3124\qquad\text{and}\qquad1324\xrightarrow{t_{13}t_{23}}2134\xrightarrow{t_{13}}3124,$$
with $\omega_2=q_2(y_3-z_1)$ and $\omega_2=q_2(y_1-z_2)$ respectively ($\fx_2(1324,2314)=\{3\}$, $\fx_2(2134,3124)=\{1\}$, the other fixed sets being empty). In the first, the cycles $(1,3)$ of $u^{-1}v$ and $(1,2,3)$ of $v^{-1}w$ share the top $3=k+1$ and the element $1$; in the second, the cycles $(2,1,3)$ and $(1,3)$ do. Both pairs are incompatible, so the coefficient of $\qsch_{3124}$ in $\qsch_{1324}\,\qsch_{3412}(x;z)$ is $0$, whereas the coefficient in $\qsch_{1324}\,\qelem_{2,2}(x;z_1)\qelem_{2,2}(x;z_2)$ is $q_2\big[(y_3-z_1)+(y_1-z_2)\big]$. The recurrence of Lemma \ref{lemma:thetaPP} accounts for this difference by $q_2$ times the coefficient of $\qsch_{3124}$ in $\qsch_{1324}\,\qsch_{d_2[1]}(x;z)$, which by (a) is the sum over the two (compatible) level-$1$ pairs $1324\tom11324\tom13124$ and $1324\tom13124\tom13124$, namely $(y_1-z_1)+(y_3-z_2)$; and indeed $(y_3-z_1)+(y_1-z_2)=(y_1-z_1)+(y_3-z_2)$, although the two sums do not agree term by term. For $w=2413$, by contrast, all four pairs
$$1324\tom21324\xrightarrow{t_{24}t_{13}}2413,\quad 1324\xrightarrow{t_{13}}2314\xrightarrow{t_{24}}2413,\quad 1324\xrightarrow{t_{24}}1423\xrightarrow{t_{13}}2413,\quad 1324\xrightarrow{t_{24}t_{13}}2413\tom22413$$
are compatible (the cycles of the two steps have distinct tops, or one step is trivial), and the coefficient of $\qsch_{2413}$ in both products is
$$(y_1-z_1)(y_3-z_1)+(y_3-z_1)(y_2-z_2)+(y_1-z_1)(y_4-z_2)+(y_2-z_2)(y_4-z_2).$$
Both coefficients were checked against the output of \cite{schubmult}.
\end{example}

\begin{corollary}[Doubled formula in lower degrees] \label{corollary:lowerdegrees}
For $u\in S_\infty$, $0\leq p_1,p_2\leq k$ and two alphabets $z,z'$, with $\qelem_{p_1,p_2;k}(x;z,z')=\mathcal Q\big(E_{p_1,k}(x;z)E_{p_2,k}(x;z')\big)$ as in \eqref{eq:qelemdouble}, and writing $n_1=n_k(u,v)$, $n_2=n_k(v,w)$ for a pair $u\tom kv\tom kw$,
\begin{multline*}
\qsch_u(x;y)\,\qelem_{p_1,p_2;k}(x;z,z')\\
=\sum_{\substack{u\tom{k}v\tom{k}w\ \text{compatible}\\ n_1\leq p_1,\ n_2\leq p_2}}q_k(u,v)\,q_k(v,w)\,E_{p_1-n_1,\,k-n_1}\big(y_{\fx_k(u,v)};z\big)\,E_{p_2-n_2,\,k-n_2}\big(y_{\fx_k(v,w)};z'\big)\,\qsch_w(x;y).
\end{multline*} For $p_1=p_2=k$ this is Theorem \ref{theorem:doubledpieri}. This is the form implemented in the routine \texttt{double\_elem\_sym\_q} of \cite{schubmult} and used by the paired step of Method \ref{method:schubmultq}.
\end{corollary}
\begin{proof}
Theorem \ref{theorem:doubledpieri} is an identity of polynomials in $z_1,z_2$, in which $\qsch_{d_2[k]}(x;z)=\mathcal Q\big(\prod_{i\leq k}(x_i-z_1)\prod_{i\leq k}(x_i-z_2)\big)$ and $\mathcal Q$ is linear over $\mathbb Z[z_1,z_2]$. Substitute $z_1=-t$, $z_2=-t'$ for independent indeterminates $t,t'$. On the left, $\prod_{i\leq k}(x_i+t)\prod_{i\leq k}(x_i+t')=\sum_{a,b=0}^kt^{k-a}t'^{\,k-b}e_a(x_1,\ldots,x_k)e_b(x_1,\ldots,x_k)$. On the right, for a pair $u\tom kv\tom kw$ with $n_1=n_k(u,v)$ and $n_2=n_k(v,w)$ we have $\prod_{a\in\fx_k(u,v)}(y_a+t)=\sum_it^{k-n_1-i}e_i(y_{\fx_k(u,v)})$ and similarly for the second step, since $|\fx_k(u,v)|=k-n_1$. Comparing the coefficients of $t^{k-a}t'^{\,k-b}$ gives, for $0\leq a,b\leq k$,
\begin{equation}\label{eq:nonfactorialdoubled}
\qsch_u(x;y)\,\mathcal Q\big(e_a(x_{\leq k})e_b(x_{\leq k})\big)=\sum_{\substack{u\tom kv\tom kw\\\text{compatible}}}q_k(u,v)q_k(v,w)\,e_{a-n_k(u,v)}\big(y_{\fx_k(u,v)}\big)\,e_{b-n_k(v,w)}\big(y_{\fx_k(v,w)}\big)\,\qsch_w(x;y).
\end{equation}
Now expand both factors by \eqref{eq:Eexpansion}: $E_{p_1,k}(x;z)E_{p_2,k}(x;z')=\sum_{a,b}c_{p_1,k,a}(z)c_{p_2,k,b}(z')\,e_a(x_{\leq k})e_b(x_{\leq k})$. Apply $\mathcal Q$, which is linear over $\mathbb Z[z,z']$, multiply by $\qsch_u(x;y)$, insert \eqref{eq:nonfactorialdoubled}, and sum first over $a$ and then over $b$ using \eqref{eq:Ereduction} with $(A,n)=(\fx_k(u,v),n_k(u,v))$ and $(\fx_k(v,w),n_k(v,w))$ respectively; this gives the asserted coefficient, which vanishes unless $n_k(u,v)\leq p_1$ and $n_k(v,w)\leq p_2$.
\end{proof}

We first reduce Theorem \ref{theorem:doubledpieri} to a single combinatorial identity about pairs of chains, \eqref{eq:Jk} below (Proposition \ref{prop:doubledreduction}). To prove \eqref{eq:Jk} we pass to a termwise reformulation \eqref{eq:Bk} based on the quantum Giambelli formula, prove that the pairs of chains organize into blocks along which \eqref{eq:Bk} telescopes (Lemma \ref{lemma:blocks}), reduce \eqref{eq:Bk} to an identity between the ends of the blocks and the mixed chains that move the position $k+1$ (Proposition \ref{prop:junction}), and prove that identity by an explicit injection, built from a single local lemma which generalizes the up/down pairing of Lemma \ref{lemma:updown} (Theorem \ref{theorem:junction}). This completes the proof of Theorem \ref{theorem:doubledpieri} for all $k$.

\subsubsection*{The recurrence}

The quantization of two equal columns is governed by the following recurrence, in which only powers of $q$ appear.

\begin{lemma}[Quantization of two equal columns] \label{lemma:thetaPP}
For all $k\geq1$,
$$\qelem_{k,k}(x;z_1)\,\qelem_{k,k}(x;z_2)=\qsch_{d_2[k]}(x;z)+q_k\,\qsch_{d_2[k-1]}(x;z).$$
\end{lemma}
\begin{proof}
Write $P_k(z_1)=\prod_{i\leq k}(x_i-z_1)=E_{k,k}(x;z)$ and $Q_k(z_1)=\mathcal Q(P_k(z_1))=\qelem_{k,k}(x;z)$, and similarly with $z_2$, with $P_0=Q_0=1$ and $P_{-1}=Q_{-1}=0$; thus $\qsch_{d_2[k]}(x;z)$ equals $\mathcal Q(P_k(z_1)P_k(z_2))$. Since $P_{k+1}(z_2)=(x_{k+1}-z_2)P_k(z_2)$, we have
$$P_k(z_1)P_{k+1}(z_2)-P_{k+1}(z_1)P_k(z_2)=(z_1-z_2)\,P_k(z_1)P_k(z_2).$$
The products $P_k(z_1)P_{k+1}(z_2)$ and $P_{k+1}(z_1)P_k(z_2)$ are $\mathbb Z[z_1,z_2]$-linear combinations of the standard elementary monomials $e_a(x_{\leq k})e_b(x_{\leq k+1})$, so $\mathcal Q$ acts on them multiplicatively, and
$$(z_1-z_2)\,\qsch_{d_2[k]}(x;z)=W_k,\qquad W_k:=Q_k(z_1)Q_{k+1}(z_2)-Q_{k+1}(z_1)Q_k(z_2).$$
The three-term recursion $Q_{k+1}(z_2)=(x_{k+1}-z_2)Q_k(z_2)+q_kQ_{k-1}(z_2)$ of Lemma \ref{lemma:elemrecursion} gives
\begin{multline*}
W_k=Q_k(z_1)\big[(x_{k+1}-z_2)Q_k(z_2)+q_kQ_{k-1}(z_2)\big]-\big[(x_{k+1}-z_1)Q_k(z_1)+q_kQ_{k-1}(z_1)\big]Q_k(z_2)\\
=(z_1-z_2)Q_k(z_1)Q_k(z_2)-q_kW_{k-1}.
\end{multline*}
Dividing by $z_1-z_2$ (all terms are divisible, $W_k$ being antisymmetric in $z_1,z_2$) gives $\qsch_{d_2[k]}=Q_kQ_k-q_k\qsch_{d_2[k-1]}$.
\end{proof}

Lemma \ref{lemma:thetaPP} reduces Theorem \ref{theorem:doubledpieri} to a recurrence in $k$. For $u,w\in S_\infty$ let $\mathrm{All}_k(u,w)$, $\mathrm{Comp}_k(u,w)$ and $\mathrm{Inc}_k(u,w)$ denote the sum of $\omega_k(u,v,w)$ over all pairs $u\tom{k}v\tom{k}w$, over the compatible pairs, and over the incompatible pairs, respectively, so that $\mathrm{All}_k=\mathrm{Comp}_k+\mathrm{Inc}_k$ and, by \eqref{eq:allk}, $\mathrm{All}_k(u,w)=[\qsch_w]\big(\qsch_u\,\qelem_{k,k}(x;z_1)\qelem_{k,k}(x;z_2)\big)$; here $\mathrm{Comp}_0(u,w)=\mathrm{All}_0(u,w)=\delta_{u,w}$.

\begin{proposition}[Reduction] \label{prop:doubledreduction}
Suppose Theorem \ref{theorem:doubledpieri} holds with $k$ replaced by $k-1$ (for $k=1$ this is the trivial statement $\qsch_u\qsch_{d_2[0]}=\qsch_u$). Then Theorem \ref{theorem:doubledpieri} holds for $k$ if and only if
\begin{equation}\label{eq:Jk}
\mathrm{Inc}_k(u,w)=q_k\,\mathrm{Comp}_{k-1}(u,w)\qquad\text{for all }u,w\in S_\infty.\tag{$\mathrm J_k$}
\end{equation}
Consequently Theorem \ref{theorem:doubledpieri} holds for all $k\leq K$ if and only if \eqref{eq:Jk} holds for all $k\leq K$.
\end{proposition}
\begin{proof}
By Lemma \ref{lemma:thetaPP}, \eqref{eq:allk} and the hypothesis,
$$[\qsch_w]\big(\qsch_u\,\qsch_{d_2[k]}\big)=\mathrm{All}_k(u,w)-q_k[\qsch_w]\big(\qsch_u\,\qsch_{d_2[k-1]}\big)=\mathrm{All}_k(u,w)-q_k\,\mathrm{Comp}_{k-1}(u,w),$$
while $\mathrm{Comp}_k=\mathrm{All}_k-\mathrm{Inc}_k$. The two sides of Theorem \ref{theorem:doubledpieri} therefore differ by $\mathrm{Inc}_k(u,w)-q_k\mathrm{Comp}_{k-1}(u,w)$. The last assertion follows by induction on $k$.
\end{proof}

The identity \eqref{eq:Jk} says that the terms omitted from the iterated Pieri formula at level $k$ are exactly accounted for, with an extra factor $q_k$, by the doubled formula one level down. Its left-hand side is a sum over pairs of chains whose cycles overlap, and the following lemma locates the overlap.

\begin{lemma}[Overlapping cycles have top $k+1$] \label{lemma:overlap}
Let $u\tom{k}v\tom{k}w$, and suppose that a cycle $c$ of $u^{-1}v$ and a cycle $d$ of $v^{-1}w$ have the same top $B$ and share an element $a\leq k$. Then $B=k+1$, and moreover $k\in c\cup d$, i.e.\ one of the two cycles contains the transposition $t_{k,k+1}$. Consequently a pair of chains in which either step contains no transposition $t_{\cdot,k+1}$ is compatible, and in an incompatible pair both steps contain a cycle with top $k+1$, these two cycles overlap, no other pair of cycles does, and the letter $t_{k,k+1}$ occurs in one of the two chains.
\end{lemma}
\begin{proof}
We show more generally that there is no position $r\notin c\cup d$ with $a'<r<B$ for every $a'\in(c\cup d)\setminus\{B\}$; applied to $r=k+1$ this gives $B=k+1$ (as $k+1\notin c\cup d$ when $B\geq k+2$), and applied to $r=k$ it gives $k\in c\cup d$. Suppose such an $r$ exists. By Lemma \ref{lemma:tight}(iv) and Lemma \ref{lemma:piericycle} the fixed sets and weights of the two chains do not depend on the order in which their cycles are processed, and every order gives a chain of covers; we process $c$ last in the first step and $d$ first in the second. Write $u'$ for the permutation obtained from $u$ by applying the other cycles of $u^{-1}v$, so that $v=u'c$ and $c$ is $u'$-tight, and write $c=(a_p,\ldots,a_1,B)$, $\beta_0=u'(B)=u(B)$, $\beta_l=u'(a_l)=u(a_l)$, and $d=(b_r,\ldots,b_1,B)$, $\gamma_0=v(B)$, $\gamma_m=v(b_m)$. By Lemma \ref{lemma:tight}(i), $v(B)=\beta_p$ and $v(a_l)=\beta_{l-1}$ for $1\leq l\leq p$; as $a=a_i=b_j$ with $i,j\geq1$, we get
$$\gamma_0=\beta_p,\qquad \gamma_j=\beta_{i-1},\qquad\text{and}\qquad \beta_{i-1}\neq\beta_p.$$
Let $\eta$ be the value at the position $r$ when the letters of $c$ are applied. Since $r\notin c\cup d$ and nothing happens between the last letter of $c$ and the first letter of $d$, $\eta$ is also the value at position $r$ when the letters of $d$ are applied. Moreover $a_l<r<B$ for every $l$ and $b_m<r<B$ for every $m$, so the conditions \eqref{eq:Cj} and \eqref{eq:Qj} apply to the position $r$ at every step of $c$ and of $d$; for $c$ they read: $\eta\notin(\beta_l,\beta_{l-1})$ if $\beta_l<\beta_{l-1}$, and $\eta\in(\beta_{l-1},\beta_l)$ if $\beta_l>\beta_{l-1}$, for $1\leq l\leq p$. By Corollary \ref{corollary:cyclic} the cyclic sequence $\Gamma=(\beta_0,\ldots,\beta_p)$ is unimodal, and the displayed conditions say precisely that $\eta$ lies in the gap $(\beta_p,\beta_0)$ of $\Gamma$: if $\Gamma$ has no descent then $\eta\notin(\beta_p,\beta_0)$, i.e.\ $\eta<\beta_p=\min\Gamma$ or $\eta>\beta_0=\max\Gamma$; and if the ascent of $\Gamma$ is $(\beta_{q-1},\beta_q)$ then $\eta\in(\min\Gamma,\max\Gamma)$ and $\eta$ avoids all the descent gaps $(\beta_l,\beta_{l-1})$ with $l\neq q$, which together with $(\beta_0,\beta_p)$ cover $(\min\Gamma,\max\Gamma)\setminus\Gamma$, so $\beta_0<\eta<\beta_p$. In other words, in the cyclically decreasing arrangement of $\Gamma\cup\{\eta\}$ the value $\eta$ sits between $\beta_p$ and $\beta_0$. The same argument for $d$ shows that in the cyclically decreasing arrangement of $\Gamma'\cup\{\eta\}$, $\Gamma'=(\gamma_0,\ldots,\gamma_r)$, the value $\eta$ sits between $\gamma_r$ and $\gamma_0=\beta_p$.

If $\eta>\beta_p$, the first statement forces $\beta_p=\min\Gamma$ and $\eta>\max\Gamma$, while the second forces $\Gamma'\cap(\beta_p,\eta)=\emptyset$; but $\beta_{i-1}\in\Gamma\setminus\{\beta_p\}$ lies in $(\beta_p,\eta)$ and belongs to $\Gamma'$, a contradiction. If $\eta<\beta_p$, the second statement forces $\eta=\min(\Gamma'\cup\{\eta\})$ and $\beta_p=\max\Gamma'$, while the first forces $\Gamma\cap(\eta,\beta_p)=\emptyset$; but $\beta_{i-1}\in\Gamma'\setminus\{\beta_p\}$ lies in $(\eta,\beta_p)$ and belongs to $\Gamma$, again a contradiction. This proves the claim, hence $B=k+1$ and $k\in c\cup d$. The consequences follow because a step contains at most one cycle with a given top.
\end{proof}

The identity \eqref{eq:Jk} will be established for all $k$ in Theorem \ref{theorem:junction}. Before that we record two features of it which any proof must account for, and which explain the detour through \eqref{eq:Bk}. First, by Lemma \ref{lemma:overlap}, in every incompatible pair the overlapping cycles have top $B=k+1$ and one of them contains $k$, so that both steps contain a transposition $t_{a,k+1}$ and one of them contains $t_{k,k+1}$; the incompatibility is a ``detour'' of the level-$(k-1)$ transposition $t_{ak}$ through the position $k+1$, in accordance with $t_{a,k+1}t_{k,k+1}t_{a,k+1}=t_{ak}$. Second, \eqref{eq:Jk} is not a termwise identity: already for $u=1324$, $k=2$, $w=3124$ (Example \ref{example:doubled}(b)) the left-hand side is $q_2\big[(y_3-z_1)+(y_1-z_2)\big]$ and the right-hand side is $q_2\big[(y_1-z_1)+(y_3-z_2)\big]$. In all examples the pairs with given $(u,w)$ and $q$-weight organize into chains $v_1,\ldots,v_{m+1}$ whose fixed sets are nested, with contributions $\prod_{j<i}(y_{c_j}-z_1)\prod_{j>i}(y_{c_j}-z_2)$, which sum to the divided difference $\big(\prod_j(y_{c_j}-z_1)-\prod_j(y_{c_j}-z_2)\big)/(z_2-z_1)$ and so depend only on the set $\{c_j\}$; the identity \eqref{eq:Jk} then amounts to a matching of these sets between the two levels. We now explain the origin of this divided-difference structure and use it to recast \eqref{eq:Jk} as a termwise identity, which we then prove.

\subsubsection*{A termwise reformulation via the Giambelli formula}

The proof of Lemma \ref{lemma:thetaPP} shows that
$$(z_1-z_2)\,\qsch_{d_2[k]}(x;z)=\qelem_{k,k}(x;z_1)\,\qelem_{k+1,k+1}(x;z_2)-\qelem_{k+1,k+1}(x;z_1)\,\qelem_{k,k}(x;z_2),$$
which is the quantum Giambelli formula \cite{fomin1997quantum} for the $k\times2$ rectangle written as a divided difference. Expanding both products on the right by Theorem \ref{theorem:pieri} at the two levels $k$ and $k+1$ gives the following reformulation. For $u,w\in S_\infty$ put
\begin{align*}
M_k(u,w)&=\sum_{u\tom{k}v\tom{k+1}w}q_k(u,v)\,q_{k+1}(v,w)\prod_{a\in\fx_k(u,v)}(y_a-z_1)\prod_{b\in\fx_{k+1}(v,w)}(y_b-z_2),\\
M'_k(u,w)&=\sum_{u\tom{k+1}v\tom{k}w}q_{k+1}(u,v)\,q_k(v,w)\prod_{a\in\fx_{k+1}(u,v)}(y_a-z_1)\prod_{b\in\fx_k(v,w)}(y_b-z_2).
\end{align*}

\begin{proposition}[Giambelli form] \label{prop:giambelliform}
For every $k\geq1$ and all $u,w\in S_\infty$,
$$(z_1-z_2)\,[\qsch_w]\big(\qsch_u(x;y)\,\qsch_{d_2[k]}(x;z)\big)=M_k(u,w)-M'_k(u,w).$$
Consequently Theorem \ref{theorem:doubledpieri} holds for $k$ if and only if $(z_1-z_2)\,\mathrm{Comp}_k(u,w)=M_k(u,w)-M'_k(u,w)$ for all $u,w$.
\end{proposition}
\begin{proof}
Multiply the displayed Giambelli identity by $\qsch_u(x;y)$ and expand $\qsch_u\,\qelem_{k,k}(x;z_1)\qelem_{k+1,k+1}(x;z_2)$ by two applications of Theorem \ref{theorem:pieri}, first at level $k$ with the variable $z_1$ and then at level $k+1$ with the variable $z_2$; this gives $M_k$. The product $\qsch_u\,\qelem_{k+1,k+1}(x;z_1)\qelem_{k,k}(x;z_2)$, expanded first at level $k+1$ and then at level $k$, gives $M'_k$.
\end{proof}

The right-hand side of Proposition \ref{prop:giambelliform} is a sum of products in which the levels of the two steps differ, and the factor $z_1-z_2$ on the left can be absorbed into the compatible pairs as follows. For a pair $u\tom{k}v\tom{k}w$ let $c(v)=v(k+1)$. If the first step does not use the position $k+1$ then $c(v)=u(k+1)$, and otherwise $c(v)=u(a)$ for a position $a\leq k$ moved by the first step; in either case $c(v)\notin\fx_k(u,v)$, and $c(v)\notin\fx_k(v,w)\subseteq\{v(1),\ldots,v(k)\}$. Hence, writing $z_1-z_2=(y_{c(v)}-z_2)-(y_{c(v)}-z_1)$,
\begin{multline*}
(z_1-z_2)\,\omega_k(u,v,w)\\
=q_k(u,v)q_k(v,w)\Big[\prod_{a\in\fx_k(u,v)}(y_a-z_1)\prod_{b\in\fx_k(v,w)\cup\{c(v)\}}(y_b-z_2)-\prod_{a\in\fx_k(u,v)\cup\{c(v)\}}(y_a-z_1)\prod_{b\in\fx_k(v,w)}(y_b-z_2)\Big],
\end{multline*}
and each of the two products on the right has the same shape as a term of $M_k$ or $M'_k$. Writing $\mathcal T(u,w)$ for the multiset of triples $\big(q_k(u,v)q_k(v,w),\ \fx_k(u,v),\ \fx_k(v,w)\cup\{c(v)\}\big)$ over the compatible pairs $u\tom{k}v\tom{k}w$, $\mathcal S(u,w)$ for the multiset of $\big(q_k(u,v)q_k(v,w),\ \fx_k(u,v)\cup\{c(v)\},\ \fx_k(v,w)\big)$ over the same pairs, $\mathcal M_k(u,w)$ for the multiset of $\big(q_k(u,v)q_{k+1}(v,w),\fx_k(u,v),\fx_{k+1}(v,w)\big)$ over the chains $u\tom{k}v\tom{k+1}w$ and $\mathcal M'_k(u,w)$ for the corresponding multiset over the chains $u\tom{k+1}v\tom{k}w$, we are led to the following termwise statement, from which Theorem \ref{theorem:doubledpieri} follows immediately by Proposition \ref{prop:giambelliform}.

\begin{theorem}[Termwise form] \label{theorem:termwise}
For all $k\geq1$ and $u,w\in S_\infty$ there is an equality of multisets
\begin{equation}\label{eq:Bk}
\mathcal T(u,w)\uplus\mathcal M'_k(u,w)=\mathcal S(u,w)\uplus\mathcal M_k(u,w).\tag{$\mathrm B_k$}
\end{equation}
\end{theorem}
\begin{proof}
Proposition \ref{prop:junction} and Theorem \ref{theorem:junction} below.
\end{proof}

Unlike \eqref{eq:Jk}, the identity \eqref{eq:Bk} is an equality of multisets, so it calls for a bijection. We first show that the greater part of this bijection is forced: the compatible pairs with given $(u,w)$ fall into explicit \emph{blocks}, inside which the terms of $\mathcal T$ and $\mathcal S$ cancel telescopically, and the unblocked ends of the blocks are matched with the chains of $\mathcal M_k$ and $\mathcal M'_k$ that do not move the position $k+1$. This reduces \eqref{eq:Bk} to an identity between the blocked ends of the blocks and the mixed chains that do move $k+1$ (Proposition \ref{prop:junction}); in the simplest situation that identity is exactly the up/down pairing of Lemma \ref{lemma:updown}, and we then prove it in general (Theorem \ref{theorem:junction}).

\subsubsection*{Blocks, and the reduction of \eqref{eq:Bk} to its junctions}

For a pair $u\tom{k}v\tom{k}w$ we write $[A\mid B]=\prod_{a\in A}(y_a-z_1)\prod_{b\in B}(y_b-z_2)$ and
$$T(v)=q_k(u,v)q_k(v,w)\,[\fx_k(u,v)\mid\fx_k(v,w)\cup\{c(v)\}],\qquad S(v)=q_k(u,v)q_k(v,w)\,[\fx_k(u,v)\cup\{c(v)\}\mid\fx_k(v,w)],$$
so that $(z_1-z_2)\,\omega_k(u,v,w)=T(v)-S(v)$ by the display above, and $\mathcal T(u,w)$, $\mathcal S(u,w)$ are the multisets of the $T(v)$, $S(v)$ over the compatible pairs (we identify the monomial $q[A\mid B]$ with the triple $(q,A,B)$). Recall from Lemma \ref{lemma:tight}(iv) and Lemma \ref{lemma:piericycle} that the cycles of a step may be processed in any order.

\begin{definition}[Block data, cut, block] \label{def:blockdata}
Let $u\tom{k}v\tom{k}w$ be compatible, $\pi_1=u^{-1}v$, $\pi_2=v^{-1}w$. Let $c$ be the cycle of $\pi_1$ with top $k+1$ and $d$ the cycle of $\pi_2$ with top $k+1$ (either may be trivial). By compatibility $c\cap d\subseteq\{k+1\}$, so writing $c=(x_i,\ldots,x_1,k+1)$ and $d=(x_n,\ldots,x_{i+1},k+1)$ we have, by \eqref{eq:pivot},
$$X:=cd=t_{x_1,k+1}\cdots t_{x_n,k+1}=(x_n,\ldots,x_1,k+1),$$
a cycle with top $k+1$ which we call the \emph{merged cycle} of the pair. Let $\rho_1=\pi_1c^{-1}$ and $\rho_2=d^{-1}\pi_2$ be the products of the remaining cycles of $\pi_1$ and of $\pi_2$; all of these have tops $\geq k+2$, and $u^{-1}w=\rho_1X\rho_2$. We call $(\rho_1,X,\rho_2)$ the \emph{block data} and $i$ the \emph{cut} of the pair, and we put
$$i_{\min}=\max\big(\{j:x_j\in\supp\rho_2\}\cup\{0\}\big),\qquad i_{\max}=\min\big(\{j:x_j\in\supp\rho_1\}\cup\{n+1\}\big)-1.$$
Since the cycles of $\pi_1$ are disjoint, $x_1,\ldots,x_i\notin\supp\rho_1$, and likewise $x_{i+1},\ldots,x_n\notin\supp\rho_2$; hence $i_{\min}\leq i\leq i_{\max}$. The \emph{block} of the pair is the set of compatible pairs $u\tom{k}v'\tom{k}w$ with the same block data; the blocks partition the compatible pairs with given $(u,w)$.
\end{definition}

\begin{lemma}[Blocks] \label{lemma:blocks}
Let $u\tom{k}v\tom{k}w$ be compatible with block data $(\rho_1,X,\rho_2)$, $X=(x_n,\ldots,x_1,k+1)$, and put $u'=u\rho_1$, $c_j=(x_j,\ldots,x_1,k+1)$, $d_j=(x_n,\ldots,x_{j+1},k+1)$ and $v_j=u'c_j$ for $0\leq j\leq n$, so that $v=v_i$ where $i$ is the cut. Then the block of the pair is $\{v_j: i_{\min}\leq j\leq i_{\max}\}$, and for every $j$ in this range:
\begin{enumerate}
\item[(i)] $u\tom{k}v_j\tom{k}w$ is a compatible pair with block data $(\rho_1,X,\rho_2)$ and cut $j$, and $q_k(u,v_j)\,q_k(v_j,w)=q_k(u,v)\,q_k(v,w)$;
\item[(ii)] $\fx_k(u,v_j)=F_1\setminus\{u(x_1),\ldots,u(x_j)\}$ and $\fx_k(v_j,w)=F_2\setminus\{w(x_{j+1}),\ldots,w(x_n)\}$, where $F_1=\{u(a):a\leq k,\ a\notin\supp\rho_1\}$ and $F_2=\{w(a):a\leq k,\ a\notin\supp\rho_2\}$;
\item[(iii)] $c(v_j)=u(x_j)$ if $j\geq1$ and $c(v_0)=u(k+1)$; moreover $w(x_j)=c(v_{j-1})$ for $i_{\min}<j\leq i_{\max}$.
\end{enumerate}
Consequently $S(v_j)=T(v_{j-1})$ for $i_{\min}<j\leq i_{\max}$, and the contribution of the block telescopes:
$$(z_1-z_2)\sum_{j=i_{\min}}^{i_{\max}}\omega_k(u,v_j,w)=T(v_{i_{\max}})-S(v_{i_{\min}}).$$
\end{lemma}
\begin{proof}
(i) For $j\leq i_{\max}$ the positions $x_1,\ldots,x_j$ avoid $\supp\rho_1$, so $\rho_1c_j$ is a product of disjoint cycles, each with exactly one element greater than $k$; similarly $d_j\rho_2$ for $j\geq i_{\min}$. By Lemma \ref{lemma:piericycle} and Lemma \ref{lemma:tight}(iv), $u\tom{k}v_j\tom{k}w$ holds if and only if $\rho_1$ is $u$-tight, $c_j$ is $u'$-tight, $d_j$ is $v_j$-tight, and $\rho_2$ is $v_jd_j$-tight, where $v_jd_j=u'X$ does not depend on $j$. The first and the last condition do not involve $j$, and they hold because they hold for $j=i$ (as $\pi_1=\rho_1c$ is $u$-tight and $\pi_2=d\rho_2$ is $v$-tight). For the two middle conditions write $\beta_0=u'(k+1)$ and $\beta_l=u'(x_l)$ for $1\leq l\leq n$. The pivot factorization of $X$ is the concatenation of those of $c_j$ and $d_j$, so by Lemma \ref{lemma:tight}(i) applied to $X$ and $u'$, the $l$-th letter $t_{x_l,k+1}$ exchanges the value $\beta_l$ at $x_l$ with the value $\beta_{l-1}$ at $k+1$, for every $l$ and independently of $j$; in particular $v_j(k+1)=\beta_j$ and $v_j(x_l)=\beta_l$ for $l>j$. Hence
$$\desc_{u'}(c_j)=\{x_l:l\leq j,\ \beta_l>\beta_{l-1}\},\qquad \desc_{v_j}(d_j)=\{x_l:l>j,\ \beta_l>\beta_{l-1}\},$$
and by Definition \ref{def:tight}
$$\Lw_{u'}(c_j)+\Lw_{v_j}(d_j)=n-2\sum_{l:\,\beta_l>\beta_{l-1}}(k+1-x_l)=\Lw_{u'}(X),\qquad \qw_{u'}(c_j)\,\qw_{v_j}(d_j)=\prod_{l:\,\beta_l>\beta_{l-1}}q_{x_l,k+1}=\qw_{u'}(X),$$
both independent of $j$. By Lemma \ref{lemma:tight}(ii),
$$\ell(u'X)-\ell(u')=\big(\ell(v_j)-\ell(u')\big)+\big(\ell(u'X)-\ell(v_j)\big)\geq\Lw_{u'}(c_j)+\Lw_{v_j}(d_j)=\Lw_{u'}(X),$$
with equality if and only if $c_j$ is $u'$-tight and $d_j$ is $v_j$-tight. Equality holds for $j=i$, hence for every $j$. By Corollary \ref{corollary:cycleweight} and Lemma \ref{lemma:tight}(iv), $q_k(u,v_j)q_k(v_j,w)=\qw_u(\rho_1)\,\qw_{u'}(c_j)\,\qw_{v_j}(d_j)\,\qw_{u'X}(\rho_2)=\qw_u(\rho_1)\qw_{u'}(X)\qw_{u'X}(\rho_2)$, independent of $j$. The pair is compatible, since its cycles with tops $\geq k+2$ are those of $\rho_1$ and $\rho_2$, as for the given pair, and $c_j\cap d_j=\{k+1\}$; its block data is $(\rho_1,X,\rho_2)$ and its cut is $j$ by construction. Conversely every pair in the block has cut in $[i_{\min},i_{\max}]$ and is therefore one of the $v_j$, since the block data and the cut determine $\pi_1=\rho_1c_j$.

(ii) The positions $a\leq k$ moved by $\pi_1^{(j)}=\rho_1c_j$ are those of $\supp\rho_1$ together with $x_1,\ldots,x_j$, which gives the first formula (note $u(x_l)\in F_1$ for $l\leq i_{\max}$). A position $a\leq k$ is fixed by $d_j\rho_2$ if and only if $a\notin\supp\rho_2\cup\{x_{j+1},\ldots,x_n\}$, and then $v_j(a)=w(a)$; this gives the second formula (note $w(x_l)\in F_2$ for $l>i_{\min}$).

(iii) We saw that $v_j(k+1)=\beta_j$, which is $u'(x_j)=u(x_j)$ for $1\leq j\leq i_{\max}$ (as $x_j\notin\supp\rho_1$) and $u'(k+1)=u(k+1)$ for $j=0$ (as $k+1\notin\supp\rho_1$). Now let $i_{\min}<j\leq i_{\max}$ and consider the pair $v_{j-1}$. The position $x_j$ is not moved by its first step $\rho_1c_{j-1}$. Its second step is $d_{j-1}\rho_2$, whose first letter $t_{x_j,k+1}$ moves the value $v_{j-1}(k+1)=c(v_{j-1})$ to the position $x_j$, and $x_j$ is moved by no later letter: not by the remaining letters of $d_{j-1}$, and not by $\rho_2$ because $j>i_{\min}$. Hence $w(x_j)=c(v_{j-1})$.

For the consequence write $A_j=\fx_k(u,v_j)$ and $B_j=\fx_k(v_j,w)$. For $i_{\min}<j\leq i_{\max}$, (ii) and (iii) give $A_j\cup\{c(v_j)\}=A_j\cup\{u(x_j)\}=A_{j-1}$ and $B_j=B_{j-1}\cup\{w(x_j)\}=B_{j-1}\cup\{c(v_{j-1})\}$, while the $q$-weights agree by (i); thus $S(v_j)=T(v_{j-1})$. Summing $(z_1-z_2)\omega_k(u,v_j,w)=T(v_j)-S(v_j)$ over the block, all terms cancel except $T(v_{i_{\max}})$ and $S(v_{i_{\min}})$.
\end{proof}

We call $v_{i_{\max}}$ the \emph{head} and $v_{i_{\min}}$ the \emph{tail} of the block; the head is \emph{blocked} if $i_{\max}<n$ (some $x_j$ lies in $\supp\rho_1$) and the tail is \emph{blocked} if $i_{\min}>0$ (some $x_j$ lies in $\supp\rho_2$). Thus a block is the set of all cuts of the merged transposition word $t_{x_1,k+1}\cdots t_{x_n,k+1}$ that are not obstructed by a spectator cycle of $\rho_1$ or $\rho_2$ sharing a position with $X$, and the transit value $c(v_j)$ runs through the values of $u$ at $x_n,x_{n-1},\ldots,x_1,k+1$ as the cut moves from right to left.

\begin{proposition}[Reduction of \eqref{eq:Bk} to junctions] \label{prop:junction}
Let $u,w\in S_\infty$. Let $\mathcal M_k^{\mathrm{mov}}(u,w)\subseteq\mathcal M_k(u,w)$ be the sub-multiset of the chains $u\tom{k}v\tom{k+1}w$ whose second step moves the position $k+1$, let $\mathcal M'^{\,\mathrm{mov}}_k(u,w)\subseteq\mathcal M'_k(u,w)$ consist of the chains $u\tom{k+1}v\tom{k}w$ whose first step moves $k+1$, and let $\mathcal T^{\mathrm{bl}}(u,w)$, respectively $\mathcal S^{\mathrm{bl}}(u,w)$, be the multiset of $T(v)$ over the blocked heads $v$, respectively of $S(v)$ over the blocked tails $v$, of the blocks of compatible pairs $u\tom{k}v\tom{k}w$. Then \eqref{eq:Bk} holds for $(u,w)$ if and only if
\begin{equation}\label{eq:junction}
\mathcal T^{\mathrm{bl}}(u,w)\uplus\mathcal M'^{\,\mathrm{mov}}_k(u,w)=\mathcal S^{\mathrm{bl}}(u,w)\uplus\mathcal M^{\mathrm{mov}}_k(u,w),\tag{$\mathrm J^{\mathrm B}_k$}
\end{equation}
and in that case $(z_1-z_2)\,\mathrm{Comp}_k(u,w)=M_k(u,w)-M'_k(u,w)$. In particular, Theorem \ref{theorem:doubledpieri} holds for $k$ as soon as \eqref{eq:junction} holds for all $u,w\in S_\infty$.
\end{proposition}
\begin{proof}
By Lemma \ref{lemma:blocks} the compatible pairs with given $(u,w)$ are partitioned into blocks, and inside each block the elements $S(v_j)$, $i_{\min}<j\leq i_{\max}$, of $\mathcal S(u,w)$ coincide with the elements $T(v_{j-1})$ of $\mathcal T(u,w)$. Cancelling them, \eqref{eq:Bk} is equivalent to
$$\{T(v):v\text{ a head}\}\uplus\mathcal M'_k(u,w)=\{S(v):v\text{ a tail}\}\uplus\mathcal M_k(u,w).$$
Now let $v=v_n$ be an unblocked head, so $d=1$ and $v(k+1)=w(k+1)$. By Lemma \ref{lemma:leveldichotomy}(i), applied with $k+1$ in place of $k$, $v\tom{k+1}w$ with $q_{k+1}(v,w)=q_k(v,w)$ and $\fx_{k+1}(v,w)=\fx_k(v,w)\cup\{v(k+1)\}$; thus $u\tom{k}v\tom{k+1}w$ is a chain of $\mathcal M_k(u,w)\setminus\mathcal M_k^{\mathrm{mov}}(u,w)$, and its triple is $T(v)$. Conversely, let $u\tom{k}v\tom{k+1}w$ be a chain not moving $k+1$. By the same lemma $v\tom{k}w$ with the same weight, the second step has no cycle with top $k+1$, so the pair $u\tom{k}v\tom{k}w$ is compatible with cut $i=n$; hence it is the head of its block, this head is unblocked ($x_1,\ldots,x_n$ avoid $\supp\rho_1$ as the cycles of $\pi_1$ are disjoint), and its triple is $T(v)$. Since distinct blocks have distinct heads, $v\mapsto v$ is a bijection between the unblocked heads and the chains of $\mathcal M_k(u,w)$ not moving $k+1$, preserving the triples. Symmetrically, the unblocked tails $v_0$ (for which $c=1$ and $u(k+1)=v_0(k+1)$) correspond bijectively to the chains $u\tom{k+1}v\tom{k}w$ not moving $k+1$, with $S(v_0)=\big(q,\fx_{k+1}(u,v_0),\fx_k(v_0,w)\big)$. Cancelling these two families gives the equivalence with \eqref{eq:junction}. The remaining assertions follow from Proposition \ref{prop:giambelliform}.
\end{proof}

The identity \eqref{eq:junction} concerns only the \emph{junctions} of the blocks: a blocked head $v_{i_{\max}}$, whose next cut is obstructed by a cycle $e$ of $\rho_1$ with top $B\geq k+2$ through $x_{i_{\max}+1}$, and a blocked tail $v_{i_{\min}}$, obstructed by a cycle of $\rho_2$ through $x_{i_{\min}}$, together with the mixed chains in which the position $k+1$ is moved by a cycle with top $B\geq k+2$. The correspondence realising \eqref{eq:junction} \emph{transports} the obstructing cycle $e$ across $X$: in the simplest case a blocked tail with block data $(\rho_1,X,\rho_2)$ and $e\subseteq\rho_2$ through $x_j$, $j=i_{\min}$, is matched with the blocked head of the block with data $(\rho_1\cdot XeX^{-1},X,\rho_2e^{-1})$, in which $e$ has been conjugated into the first step and is attached to $x_{j-1}$, so that the cut $j-1$ is skipped; in general pivots of $X$ migrate into the transported cycle or out of it, exactly as in the transposition words $F_m$ of Lemma \ref{lemma:updown}, and a chain of $\mathcal M_k^{\mathrm{mov}}$ is matched in the same way with the head obtained by transporting the cycle through $k+1$ into the first step. The simplest junction, in which the obstructing cycle is a single transposition $t_{k+1,B}$, is exactly Lemma \ref{lemma:updown} (applied with $k+1,B$ in the roles of $k,B$): its transposition words $F_{\mathrm{up}}$, $F_{\mathrm{std}}$ and $F_m$ are a chain of $\mathcal M_k^{\mathrm{mov}}$, a chain of $\mathcal M'^{\,\mathrm{mov}}_k$ with trivial second step, and the blocked heads with block data $\big((x_n,\ldots,x_m,B),(x_m,\ldots,x_1,k+1),1\big)$, and they reappear below as the junction word $J=F_{\mathrm{up}}$ and its resolutions $J_0=F_{\mathrm{std}}$, $J_m\sim F_m$.

\subsection*{Proof of the junction identity}

We now prove \eqref{eq:junction}, and hence \eqref{eq:Bk} and Theorem \ref{theorem:doubledpieri}, for all $k$. Throughout we write $K=k+1$. The proof constructs, for every $u,w$ and every $q$-monomial $q$, an injective map from the \emph{sources} (the blocked tails and the chains of $\mathcal M_k^{\mathrm{mov}}$ of weight $q$) to the \emph{targets} (the blocked heads and the chains of $\mathcal M'^{\,\mathrm{mov}}_k$ of weight $q$) which preserves the triples, and a symmetry which exchanges sources and targets; the two together give \eqref{eq:junction}. The map is a sequence of local rewritings of the pair of chains, each governed by a single lemma (Lemma \ref{lemma:token}) which generalizes Lemma \ref{lemma:updown}; the lemma is proved not by tracking covers, but by comparing the formal $q$-weights of transposition words.

\subsubsection*{Overview}
For orientation we summarize the four kinds of configurations in \eqref{eq:junction} and the triples $(q,A,B)$, $A,B\subseteq\{1,2,\ldots\}$ sets of values, that they carry; here $M_1$, $M_2$ denote the sets of positions $\leq k$ moved by the first and by the second step, $\overline S=[1,k]\setminus S$, and $q$ is the product of the $q$-weights of the two steps.
\begin{center}
\begin{tabular}{llll}
configuration & steps & role & triple\\\hline
blocked tail $v$ & $u\tom kv\tom kw$ & source & $S(v)=\big(q,\ u(\overline{M_1})\cup\{v(K)\},\ w(\overline{M_2})\big)$\\
chain of $\mathcal M_k^{\mathrm{mov}}$ & $u\tom kv\tom{k+1}w$, $K$ moved by step 2 & source & $\big(q,\ u(\overline{M_1}),\ w(\overline{M_2})\big)$\\
blocked head $v$ & $u\tom kv\tom kw$ & target & $T(v)=\big(q,\ u(\overline{M_1}),\ w(\overline{M_2})\cup\{v(K)\}\big)$\\
chain of $\mathcal M'^{\,\mathrm{mov}}_k$ & $u\tom{k+1}v\tom kw$, $K$ moved by step 1 & target & $\big(q,\ u(\overline{M_1}),\ w(\overline{M_2})\big)$
\end{tabular}
\end{center}
(For a blocked tail, $v(K)=u(x)$ where $x\leq k$ is the position whose next cut is obstructed, so $u(\overline{M_1})\cup\{v(K)\}=u(\overline{M_1\setminus\{x\}})$; for a blocked head, $v(K)=w(x_g)$ for the obstructed position $x_g$, so the last entry is $w(\overline{M_2\setminus\{x_g\}})$. In the mixed chains one of the two fixed sets is taken at level $k+1$, which accounts for the missing extra element.) In every case the triple is thus determined by the weight and by the two sets of positions $\leq k$ moved by the two steps, with the obstructed position removed from the appropriate set. The junction move of Proposition \ref{prop:process} rewrites a source as a transposition word with a cut and transports the obstructing cycle from the second step to the first across the top-$K$ letters; it preserves the product $u^{-1}w$, the weight (Lemma \ref{lemma:token}), and the two sets of moved positions up to the obstructed one (the invariants $N_1$ and $M_2$ below), and therefore the triple. After finitely many moves the result is a target, and the symmetry $v\mapsto w_0v$ of Lemma \ref{lemma:reversal}, which exchanges sources and targets, gives the reverse inclusion. A complete instance is worked out in Example \ref{example:junction}.

\subsubsection*{Transposition words, cover paths and formal weights}
A \emph{transposition word} is a finite sequence $W=(t_1,\ldots,t_N)$ of transpositions $t_i=t_{a_ib_i}$, $a_i<b_i$ (not of simple reflections: the pivot words and junction words below are transposition words, and the only words in the simple reflections used in this paper are the reduced words of Section \ref{section:proofvside}); its product is $t_1\cdots t_N$ and its length is $N$. Given a permutation $\psi$ put $\psi_0=\psi$ and $\psi_i=\psi_{i-1}t_i$. The $i$-th letter \emph{exchanges} the values $\psi_{i-1}(a_i)$ and $\psi_{i-1}(b_i)$, and we call it a \emph{descent letter} (from $\psi$) if $\psi_{i-1}(a_i)>\psi_{i-1}(b_i)$. The \emph{formal weight} of $W$ from $\psi$ is
$$\mathrm{wt}_\psi(W)=\prod_{i\text{ a descent letter}}q_{a_ib_i},$$
and $W$ is a \emph{cover path} from $\psi$ if $\psi_{i-1}\lup\psi_i$ for every $i$; in that case $\mathrm{wt}_\psi(W)=\prod_i\kappa(\psi_{i-1},t_i)$. Two transposition words are \emph{equivalent} if they have the same product and the same length.

\begin{lemma}[Degree criterion] \label{lemma:degree}
Let $W$ be a transposition word of length $N$ with product $\pi$. Then $\ell(\psi\pi)-\ell(\psi)\geq N-\deg\mathrm{wt}_\psi(W)$, with equality if and only if $W$ is a cover path from $\psi$. Consequently, if $W$ is a cover path from $\psi$ and $W'$ is equivalent to $W$, then $W'$ is a cover path from $\psi$ if and only if $\deg\mathrm{wt}_\psi(W')=\deg\mathrm{wt}_\psi(W)$; in particular this is the case whenever $\mathrm{wt}_\psi(W')=\mathrm{wt}_\psi(W)$.
\end{lemma}
\begin{proof}
By Lemma \ref{lemma:lengthchange}, the $i$-th step changes the length by at least $1$ if the letter is not a descent letter and by at least $1-2(b_i-a_i)=1-\deg q_{a_ib_i}$ if it is, with equality in both cases exactly when $\psi_{i-1}\lup\psi_i$. Summing over $i$ gives the inequality and the characterization of equality, and the rest follows because $\ell(\psi\pi)-\ell(\psi)$ and $N$ are the same for $W$ and $W'$.
\end{proof}

\begin{lemma}[Inert reordering] \label{lemma:inert}
Let $W=W_1UVW_2$ be a concatenation of transposition words such that no position is moved both by a letter of $U$ and by a letter of $V$. Then $W'=W_1VUW_2$ is equivalent to $W$, each letter of $W'$ exchanges the same two values from $\psi$ as the corresponding letter of $W$, $\mathrm{wt}_\psi(W')=\mathrm{wt}_\psi(W)$, and $W'$ is a cover path from $\psi$ if and only if $W$ is.
\end{lemma}
\begin{proof}
The letters of $U$ and $V$ commute, so the products agree. A letter of $U$ only moves positions untouched by $V$, so it finds the same values at its two positions whether or not $V$ has been applied before it, and symmetrically for $V$; the letters of $W_1$ and $W_2$ see the same permutations in both transposition words. The last assertion is Lemma \ref{lemma:degree}.
\end{proof}

\begin{lemma}[Reversal] \label{lemma:reversal}
Work in $S_n$ and let $w_0$ be the longest element. Let $W=(t_1,\ldots,t_N)$ be a transposition word with product $\pi$, and let $W^{\mathrm{rev}}=(t_N,\ldots,t_1)$. Then the $i$-th letter of $W$ is a descent letter from $\psi$ if and only if the $(N+1-i)$-th letter of $W^{\mathrm{rev}}$ is a descent letter from $w_0\psi\pi$; hence $\mathrm{wt}_{w_0\psi\pi}(W^{\mathrm{rev}})=\mathrm{wt}_\psi(W)$, and $W$ is a cover path from $\psi$ if and only if $W^{\mathrm{rev}}$ is a cover path from $w_0\psi\pi$.
\end{lemma}
\begin{proof}
The $(N+1-i)$-th step of $W^{\mathrm{rev}}$ goes from $w_0\psi_i$ to $w_0\psi_{i-1}=(w_0\psi_i)t_i$ and exchanges the values $w_0\psi_i(a_i)=n+1-\psi_{i-1}(b_i)$ and $w_0\psi_i(b_i)=n+1-\psi_{i-1}(a_i)$; the first is larger than the second if and only if $\psi_{i-1}(a_i)>\psi_{i-1}(b_i)$. Since $\ell(w_0v)=\ell(w_0)-\ell(v)$, the length change of this step is $\ell(\psi_i)-\ell(\psi_{i-1})$, the length change of the $i$-th step of $W$, so the two steps are covers of the same kind simultaneously (Definition \ref{def:cover}).
\end{proof}

\subsubsection*{The token lemma}
Fix $C\geq k+2$ and write $T_x=t_{xK}$, $S_x=t_{xC}$ for $x\leq k$, and $Z=t_{KC}$; thus for every $x\leq k$ the three products $T_xS_x=S_xZ=ZT_x$ coincide (they are the $3$-cycle $(x,K,C)$ in some order), and likewise $T_xZ=ZS_x=S_xT_x$. Let $x_1,\ldots,x_{j-1}$ be distinct positions $\leq k$, let $x_j$ be either a further position $\leq k$ (the \emph{tail case}) or $x_j=K$ (the \emph{$\mathcal M$ case}), and let $e_1,\ldots,e_{s-1}$ be distinct positions $\leq k$ different from $x_j$ (they may coincide with some of $x_1,\ldots,x_{j-1}$). The \emph{junction word} attached to these data is the transposition word
$$J=T_{x_1}\cdots T_{x_j}\cdot S_{e_1}\cdots S_{e_{s-1}}\cdot S_{x_j},$$
where in the $\mathcal M$ case the letter $T_{x_j}=T_K$ is omitted and $S_{x_j}=Z$. For $0\leq m\leq j-1$ its \emph{$m$-th resolution} is
\begin{align*}
J_m&=S_{x_m}S_{x_{m+1}}\cdots S_{x_{j-1}}\cdot T_{x_1}\cdots T_{x_{m-1}}\cdot T_{x_m}\cdot T_{e_1}\cdots T_{e_{s-1}}\cdot T_{x_j}\qquad(1\leq m\leq j-1),\\
J_0&=Z\,S_{x_1}\cdots S_{x_{j-1}}\cdot T_{e_1}\cdots T_{e_{s-1}}\cdot T_{x_j},
\end{align*}
again with $T_{x_j}$ omitted in the $\mathcal M$ case. (In the language of the detours: the letter $S_{x_j}$ is converted into the token $Z$, which travels to the left through the letters $S_{e_{s-1}},\ldots,S_{e_1}$, turning them into $T_{e_{s-1}},\ldots,T_{e_1}$, and then through $T_{x_{j-1}},\ldots,T_{x_{m+1}}$, turning them into $S$'s, until it is absorbed by $T_{x_m}$, which it turns into $S_{x_m}T_{x_m}$; for $m=0$ it is not absorbed.) Lemma \ref{lemma:updown} is the case $s=1$ of the $\mathcal M$ case: there $J=F_{\mathrm{up}}$, $J_0=F_{\mathrm{std}}$ and $J_m$ is $F_m$ with its first two blocks interchanged.

\begin{lemma} \label{lemma:resolutionproduct}
Each $J_m$ is equivalent to $J$.
\end{lemma}
\begin{proof}
The lengths agree. In the tail case, $T_{x_j}$ commutes with $S_{e_1},\ldots,S_{e_{s-1}}$, and $T_{x_j}S_{x_j}=ZT_{x_j}$; in the $\mathcal M$ case the last letter is $Z$ already. In both cases $J=T_{x_1}\cdots T_{x_{j-1}}\,S_{e_1}\cdots S_{e_{s-1}}\,Z\,T_{x_j}$ (with $T_{x_j}$ omitted in the $\mathcal M$ case). Using $S_eZ=ZT_e$ repeatedly, $J=T_{x_1}\cdots T_{x_{j-1}}\,Z\,T_{e_1}\cdots T_{e_{s-1}}T_{x_j}$. Using $T_xZ=ZS_x$ for $x=x_{j-1},\ldots,x_{m+1}$ and then $T_{x_m}Z=S_{x_m}T_{x_m}$, we get $T_{x_1}\cdots T_{x_{j-1}}Z=T_{x_1}\cdots T_{x_{m-1}}\,S_{x_m}T_{x_m}\,S_{x_{m+1}}\cdots S_{x_{j-1}}$ for $m\geq1$, and $=Z\,S_{x_1}\cdots S_{x_{j-1}}$ for $m=0$. Finally the letters $S_{x_{m+1}},\ldots,S_{x_{j-1}}$ commute with $T_{x_1},\ldots,T_{x_m}$ and may be moved to the front; this gives $J_m$.
\end{proof}

\begin{lemma}[Token lemma] \label{lemma:token}
Let $\psi\in S_\infty$ and suppose that the junction word $J$ is a cover path from $\psi$. Put $\beta_0=\psi(K)$, $\beta_l=\psi(x_l)$ for $1\leq l\leq j-1$, and $\gamma=\psi(C)$. Then the cyclic sequence $\Gamma=(\beta_0,\beta_1,\ldots,\beta_{j-1})$ is unimodal; let $g\in\{0,1,\ldots,j-1\}$ be the insertion gap of $\gamma$ in $\Gamma$ (gap $g\geq1$ is $(\beta_{g-1},\beta_g)$, gap $0$ is the wrap $(\beta_{j-1},\beta_0)$; $g=0$ if $j=1$).
\begin{enumerate}
\item[(i)] For $0\leq m\leq j-1$, the resolution $J_m$ is a cover path from $\psi$ with $\mathrm{wt}_\psi(J_m)=\mathrm{wt}_\psi(J)$ if and only if $m=g$.
\item[(ii)] If $e_t=x_l$ for some $t\leq s-1$ and $l\leq j-1$, then $l>g$. If moreover the transposition word $J\,S_{e_{s+1}}\cdots S_{e_r}$ is a cover path from $\psi$ for some further distinct positions $e_{s+1},\ldots,e_r\leq k$ (different from $x_j$ and from $e_1,\ldots,e_{s-1}$), and $e_t=x_l$ for some $t\geq s+1$ and $l\leq j-1$, then $l\leq g$.
\end{enumerate}
\end{lemma}
\begin{proof}
Write $T=T_{x_1}\cdots T_{x_{j'}}$, where $j'=j$ in the tail case and $j'=j-1$ in the $\mathcal M$ case, and $\beta_j=\psi(x_j)$ in the tail case. The transposition word $T$ is the pivot chain of the cycle $(x_{j'},\ldots,x_1,K)$, so by Lemma \ref{lemma:tight}(i) its letter $T_{x_l}$ exchanges $\beta_l$ (at $x_l$) with $\beta_{l-1}$ (at $K$), and after $T$ the position $x_l$ carries $\beta_{l-1}$ and $K$ carries $\beta_{j'}$. The remaining letters $S_{e_1},\ldots,S_{e_{s-1}},S_{x_j}$ of $J$ form the pivot chain of the cycle $(x_j,e_{s-1},\ldots,e_1,C)$ (in the $\mathcal M$ case, of $(K,e_{s-1},\ldots,e_1,C)$) from $\psi'=\psi T$; put $\varepsilon_0=\gamma$, $\varepsilon_t=\psi'(e_t)$ for $t\leq s-1$, and $\varepsilon_s=\psi'(x_j)=\beta_{j-1}$ (in both cases). By Lemma \ref{lemma:tight}(i) the letter $S_{e_t}$ exchanges $\varepsilon_t$ with $\varepsilon_{t-1}$ and $S_{x_j}$ exchanges $\varepsilon_s$ with $\varepsilon_{s-1}$. Since $J$ is a cover path, both cycles are tight (Lemma \ref{lemma:tight}(iii)), so by Corollary \ref{corollary:cyclic} the cyclic sequences $(\beta_0,\ldots,\beta_{j'})$ and $E=(\varepsilon_0,\ldots,\varepsilon_s)$ are unimodal; in particular $\Gamma$, a block of the first, is unimodal. We record
\begin{equation}\label{eq:wtJ}
\mathrm{wt}_\psi(J)=\prod_{l=1}^{j'}q_{x_lK}^{[\beta_l>\beta_{l-1}]}\ \prod_{t=1}^{s}q_{e_tC}^{[\varepsilon_t>\varepsilon_{t-1}]},
\end{equation}
where $[P]\in\{0,1\}$ is the truth value of $P$ and $e_s$ stands for $x_j$ (tail case) or $K$ ($\mathcal M$ case).

Now track $J_m$ for $m\geq1$. Its first block is the pivot chain of $(x_{j-1},\ldots,x_m,C)$ from $\psi$: $S_{x_m}$ exchanges $\beta_m$ with $\gamma$, and $S_{x_l}$ ($m<l<j$) exchanges $\beta_l$ with $\beta_{l-1}$; afterwards $x_m$ carries $\gamma$, $x_l$ carries $\beta_{l-1}$ for $m<l<j$, and $C$ carries $\beta_{j-1}$. The block $T_{x_1}\cdots T_{x_{m-1}}$ exchanges $\beta_l$ with $\beta_{l-1}$ and leaves $\beta_{m-1}$ at $K$. Then $T_{x_m}$ exchanges $\gamma$ (at $x_m$) with $\beta_{m-1}$ (at $K$), after which $K$ carries $\gamma$ and $x_m$ carries $\beta_{m-1}$. Next, $T_{e_t}$ ($t\leq s-1$) exchanges the value at $e_t$ with the value at $K$. The value at $e_t$ is $\varepsilon_t$: if $e_t$ is none of the $x_l$ it is $\psi(e_t)=\psi'(e_t)$, and if $e_t=x_l$ ($l\leq j-1$) it is $\beta_{l-1}$ in the present transposition word (by the above, in each of the cases $l<m$, $l=m$, $l>m$) and also $\psi'(x_l)=\beta_{l-1}$. The value at $K$ is $\varepsilon_{t-1}$ by induction. Finally, in the tail case, $T_{x_j}$ exchanges $\beta_j$ (the position $x_j$ has not been touched) with $\varepsilon_{s-1}$. For $m=0$ the tracking is the same except that $Z$ exchanges $\beta_0$ with $\gamma$ and then $S_{x_l}$ ($1\leq l\leq j-1$) exchanges $\beta_l$ with $\beta_{l-1}$, after which $K$ carries $\gamma$. Hence, with $q_{aC}=q_{aK}q_{KC}$ for $a\leq k$ and \eqref{eq:wtJ},
\begin{equation}\label{eq:wtJm}
\mathrm{wt}_\psi(J_m)=\mathrm{wt}_\psi(J)\cdot q_{x_mK}^{\,P_m}\;q_{x_jK}^{\,Q}\;q_{KC}^{\,a_m-a_E},
\end{equation}
where $P_0=0$, $P_m=[\beta_m>\gamma]+[\gamma>\beta_{m-1}]-[\beta_m>\beta_{m-1}]$ for $m\geq1$, $Q=[\beta_j>\varepsilon_{s-1}]-[\beta_j>\beta_{j-1}]-[\beta_{j-1}>\varepsilon_{s-1}]$ in the tail case and $Q=0$ in the $\mathcal M$ case, $a_m$ is the number of ascents of the sequence $(\gamma,\beta_m,\beta_{m+1},\ldots,\beta_{j-1})$ (for $m=0$: of $(\gamma,\beta_0,\ldots,\beta_{j-1})$), and $a_E$ is the number of ascents of $E=(\varepsilon_0,\ldots,\varepsilon_s)$.

\emph{$Q=0$.} In the tail case, the letter $S_{x_j}$ of $J$ is a cover of the permutation reached before it, in which $x_j$ carries $\beta_{j-1}$, $C$ carries $\varepsilon_{s-1}$, and the position $K$, which lies strictly between $x_j$ and $C$, carries $\beta_j$. By Definition \ref{def:cover}, if $\beta_{j-1}<\varepsilon_{s-1}$ then $\beta_j\notin(\beta_{j-1},\varepsilon_{s-1})$, and if $\beta_{j-1}>\varepsilon_{s-1}$ then $\beta_j\in(\varepsilon_{s-1},\beta_{j-1})$; in each case $[\beta_j>\varepsilon_{s-1}]=[\beta_j>\beta_{j-1}]+[\beta_{j-1}>\varepsilon_{s-1}]$.

\emph{$a_E=[\beta_{j-1}>\gamma]$.} The cyclic sequence $E$ is unimodal, so as a linear sequence it has at most one ascent; if it has none then $\varepsilon_s<\varepsilon_0$, and if its ascent is $(\varepsilon_{t-1},\varepsilon_t)$ then by Corollary \ref{corollary:cyclic} $\varepsilon_t>\cdots>\varepsilon_s>\varepsilon_0>\cdots>\varepsilon_{t-1}$, so $\varepsilon_s>\varepsilon_0$. Thus $a_E=[\varepsilon_s>\varepsilon_0]=[\beta_{j-1}>\gamma]$. The same argument applies to any linear block of a unimodal cyclic sequence: \emph{the block has an ascent if and only if its last entry exceeds its first.}

\emph{$P_m\geq0$, with equality if and only if $m=g$ $(1\leq m\leq j-1)$.} If $\beta_m<\beta_{m-1}$ then $P_m=0$ exactly when $\beta_m<\gamma<\beta_{m-1}$, and $P_m=1$ otherwise. If $\beta_m>\beta_{m-1}$, then $(\beta_{m-1},\beta_m)$ is the ascent of $\Gamma$, so $\beta_{m-1}=\min\Gamma$ and $\beta_m=\max\Gamma$, and $P_m=0$ exactly when $\gamma\notin(\beta_{m-1},\beta_m)$, i.e.\ when $\gamma$ lies outside $[\min\Gamma,\max\Gamma]$; again $P_m=1$ otherwise. By the definition of the insertion gap, in both cases $P_m=0$ if and only if gap $m$ is the insertion gap of $\gamma$ in $\Gamma$.

(i) Since the monomials $q_{x_mK}$, $q_{x_jK}$ involve only $q_i$ with $i<K$ and $q_{KC}$ only $q_i$ with $i\geq K$, and $Q=0$, \eqref{eq:wtJm} shows that $\mathrm{wt}_\psi(J_m)=\mathrm{wt}_\psi(J)$ if and only if $P_m=0$ and $a_m=a_E$; by Lemmas \ref{lemma:resolutionproduct} and \ref{lemma:degree} this is also equivalent to $J_m$ being a cover path of weight $\mathrm{wt}_\psi(J)$. By the previous step, $P_m=0$ holds for $m=0$ and for $m=g$ (if $g\geq1$) only. Let $\Sigma$ be the unimodal cyclic sequence obtained by inserting $\gamma$ into its insertion gap of $\Gamma$.

If $g\geq1$: the sequence $(\gamma,\beta_g,\ldots,\beta_{j-1})$ is a block of $\Sigma$, so $a_g=[\beta_{j-1}>\gamma]=a_E$, and $J_g$ qualifies. For $m=0$ we show $a_0\neq a_E$. Let $d\in\{0,1\}$ be the number of linear ascents of $(\beta_0,\ldots,\beta_{j-1})$, so that $a_0=[\beta_0>\gamma]+d$. If $d=0$ then $\beta_0>\cdots>\beta_{j-1}$ and the wrap is the ascent of $\Gamma$; as $g\neq0$, $\gamma$ lies in $(\beta_{j-1},\beta_0)$, so $a_0=1$ while $a_E=[\beta_{j-1}>\gamma]=0$. If $d=1$, the wrap $(\beta_{j-1},\beta_0)$ is a descent gap, and $g\neq0$ means $\gamma\notin(\beta_0,\beta_{j-1})$: if $\gamma<\beta_0$ then $a_0=2\neq a_E$, and if $\gamma>\beta_{j-1}$ then $a_0=1$ while $a_E=0$. So $J_0$ does not qualify, and $J_g$ is the unique resolution that does.

If $g=0$: only $m=0$ is possible, and $(\gamma,\beta_0,\ldots,\beta_{j-1})$ is the whole of $\Sigma$ read from $\gamma$, a linear block of a unimodal cyclic sequence, so $a_0=[\beta_{j-1}>\gamma]=a_E$ and $J_0$ qualifies.

(ii) Recall that a cyclic subsequence of a unimodal cyclic sequence is unimodal, and that a unimodal cyclic sequence of three distinct numbers is cyclically decreasing, i.e.\ reads $(\max,\mathrm{mid},\min)$ up to rotation; so two unimodal cyclic sequences containing the same three distinct numbers list them in the same cyclic order. Let $e_t=x_l$ with $t\leq s-1$, $l\leq j-1$, so that $\varepsilon_t=\beta_{l-1}$, and suppose $l\leq g$ (so $g\geq1$). In $\Sigma=(\beta_0,\ldots,\beta_{g-1},\gamma,\beta_g,\ldots,\beta_{j-1})$ the three distinct numbers $\beta_{l-1},\gamma,\beta_{j-1}$ occur in the cyclic order $(\beta_{l-1},\gamma,\beta_{j-1})$, because $l-1\leq g-1$ and $l-1<j-1$. In the unimodal cyclic sequence $E$ they occur as $\varepsilon_0=\gamma$, $\varepsilon_t=\beta_{l-1}$, $\varepsilon_s=\beta_{j-1}$ with $0<t<s$, i.e.\ in the cyclic order $(\gamma,\beta_{l-1},\beta_{j-1})$, which is the opposite one; this contradiction proves $l>g$. For the second assertion, the whole cycle $(e_r,\ldots,e_1,C)$ is tight from $\psi'$ and its cyclic value sequence $(\varepsilon_0,\ldots,\varepsilon_r)$, with $\varepsilon_t=\psi'(e_t)=\beta_{l-1}$ for $e_t=x_l$, is unimodal. Suppose $e_t=x_l$ with $t\geq s+1$ and $g<l\leq j-1$. In $\Sigma$ the numbers $\gamma,\beta_{l-1},\beta_{j-1}$ occur in the cyclic order $(\gamma,\beta_{l-1},\beta_{j-1})$, because $g\leq l-1<j-1$; in the $e$-sequence they occur as $\varepsilon_0,\varepsilon_s,\varepsilon_t$ with $0<s<t$, i.e.\ in the order $(\gamma,\beta_{j-1},\beta_{l-1})$, the opposite one. Hence $l\leq g$.
\end{proof}

\subsubsection*{States and the junction process}
Fix $u,w\in S_\infty$ and a monomial $q$. By Lemma \ref{lemma:tight}(iii), a level-$k$ step $v\tom{k}v'$ may be represented by any \emph{pivot word} of $v^{-1}v'$, i.e.\ any transposition word obtained by interleaving the pivot factorizations of its cycles; all of them are cover paths from $v$ of weight $q_k(v,v')$, and the same holds for level-$(k+1)$ steps. A \emph{state} is a transposition word $U=W_1\,|\,W_2$ with a cut, which is a cover path from $u$ of weight $q$ with product $u^{-1}w$, together with an \emph{active position} $x$, such that:
\begin{enumerate}
\item[(S1)] $W_1=R_1\,T_{x_1}\cdots T_{x_j}$, where $R_1$ is a pivot word of a product of pairwise disjoint cycles, each consisting of its top, which is $\geq k+2$, together with positions $\leq k$, and $x_1,\ldots,x_j$ are distinct positions $\leq k$ with $x_j=x$ and $x_1,\ldots,x_{j-1}\notin\supp R_1$ (the active position may or may not lie in $\supp R_1$); or, in the \emph{$\mathcal M$ case}, $x=K$ and $W_1=R_1T_{x_1}\cdots T_{x_{j-1}}$ with $x_1,\ldots,x_{j-1}$ as before.
\item[(S2)] $W_2=E\,W_2'$, where $E=S_{e_1}\cdots S_{e_r}$ is the pivot factorization of a cycle $e=(e_r,\ldots,e_1,C)$ with top $C\geq k+2$ containing $x$, say $x=e_s$, whose other elements $e_t$ ($t\neq s$) are positions $\leq k$; and $W_2'$ is a pivot word of a product of pairwise disjoint cycles, each consisting of its top, which is $\geq K$, together with positions $\leq k$, disjoint from $e$ and not containing $x$.
\item[(S3)] The cycle of $W_2'$ with top $K$, if there is one, is disjoint from $\{x_1,\ldots,x_{j-1}\}$.
\end{enumerate}
Here $T_y=t_{yK}$ and $S_y=t_{yC}$ as before (with $C$ the top of $e$). Conditions (S1) and (S2) say that every cycle of $R_1$ and of $W_2'$, and in the tail case also $e$, has exactly one element greater than $k$, as the cycles of a level-$k$ step do, while in the $\mathcal M$ case $e$ has exactly one element greater than $K$; the requirement on the tops alone would not suffice, since a cycle such as $(K,C')$ with $C'\geq k+2$ has top $\geq K$ but is not a cycle of any level-$k$ or level-$(k+1)$ step. Two states are identified if they differ by reorderings as in Lemma \ref{lemma:inert}, i.e.\ we regard a state as the data (cycles of $R_1$, the sequence $x_1,\ldots,x_j$, the cycle $e$, the cycles of $W_2'$) rather than as a particular interleaving. A state is a \emph{source} if $x\notin\supp R_1$ (and the active position is in $[1,k]$ or $K$ accordingly); the sources are exactly the blocked tails and the chains of $\mathcal M_k^{\mathrm{mov}}$: a blocked tail $u\tom{k}v\tom{k}w$ with block data $(\rho_1,X,\rho_2)$ and cut $i=i_{\min}>0$ is the state with $R_1$ a pivot word of $\rho_1$, $(x_1,\ldots,x_j)=(x_1,\ldots,x_i)$, $e$ the cycle of $\rho_2$ through $x_i$, and $W_2'$ a pivot word of $d_i\rho_2e^{-1}$; a chain $u\tom{k}v\tom{k+1}w$ of $\mathcal M_k^{\mathrm{mov}}$ with $X=(x_n,\ldots,x_1,K)$ is the state with $(x_1,\ldots,x_{j-1})=(x_1,\ldots,x_n)$, $x=K$, $e$ the cycle of $v^{-1}w$ through $K$, and $W_2'$ a pivot word of the other cycles. In both cases (S1) and (S2) hold because the two steps are $k$-shaped, respectively $k$- and $(k+1)$-shaped (the cycles of $\rho_1$, $\rho_2$ and $v^{-1}w$ other than those through $K$ have tops $\geq k+2$); (S3) holds for a blocked tail because its top-$K$ cycle $d_i=(x_n,\ldots,x_{i+1},K)$ of $W_2'$ avoids $x_1,\ldots,x_{i-1}$, and vacuously for a chain of $\mathcal M_k^{\mathrm{mov}}$, whose $W_2'$ contains no cycle through $K$. Conversely such a state is a blocked tail (resp.\ a chain of $\mathcal M_k^{\mathrm{mov}}$) of weight $q$: since $x\notin\supp R_1$, (S1) and (S2) make $W_1$ and $W_2$ pivot words of $k$-shaped permutations (in the $\mathcal M$ case, of a $k$-shaped and a $(k+1)$-shaped permutation), which are tight because $U$ is a cover path (Lemma \ref{lemma:tight}(iii)), so $u\tom{k}v\tom{k}w$ (resp.\ $u\tom{k}v\tom{k+1}w$ with $K$ moved) for $v=uW_1$; in the tail case the pair is compatible, its two cycles with top $K$, namely $(x_j,\ldots,x_1,K)$ and the top-$K$ cycle of $W_2'$, meeting only in $K$ by (S3) and $x\notin\supp W_2'$, and cycles with a common top $\geq k+2$ sharing no element $\leq k$ by Lemma \ref{lemma:overlap}; and its cut $j$ is $i_{\min}$, which is positive, because $x_j=x\in\supp e\subseteq\supp\rho_2$. For a state $U$ with active position $x$ put
$$N_1(U)=\supp R_1\cup\{x_1,\ldots,x_{j-1}\},\qquad M_2(U)=\{a\leq k: a\text{ is moved by }W_2\},$$
the set of positions $\leq k$ moved by $W_1$ other than (possibly) the active one, and the set of positions $\leq k$ moved by $W_2$.

The \emph{junction move} transforms a state $U$ as follows. The transposition word $J=T_{x_1}\cdots T_{x_j}\,S_{e_1}\cdots S_{e_{s-1}}\,S_x$ formed by the $K$-letters of $W_1$ and the first $s$ letters of $E$ is a junction word in the sense of the token lemma (tail case if $x\leq k$, $\mathcal M$ case if $x=K$), it occupies a contiguous segment of $U$, and it is a cover path from $\psi=uR_1$. Let $g$ be the index provided by Lemma \ref{lemma:token}, and replace $J$ by $J_g$; by Lemmas \ref{lemma:degree}, \ref{lemma:resolutionproduct} and \ref{lemma:token}, the resulting transposition word
$$U^+=\begin{cases}R_1\,S_{x_g}\cdots S_{x_{j-1}}\,T_{x_1}\cdots T_{x_{g-1}}\ \big|\ T_{x_g}\,T_{e_1}\cdots T_{e_{s-1}}\,T_x\,S_{e_{s+1}}\cdots S_{e_r}\,W_2'&(g\geq1),\\[3pt]
R_1\,Z\,S_{x_1}\cdots S_{x_{j-1}}\ \big|\ T_{e_1}\cdots T_{e_{s-1}}\,T_x\,S_{e_{s+1}}\cdots S_{e_r}\,W_2'&(g=0),\end{cases}$$
with $T_x$ omitted when $x=K$ and with the cut placed as indicated, is a cover path from $u$ of weight $q$ with product $u^{-1}w$.

\begin{proposition}[The junction process] \label{prop:process}
Let $U$ be a state with active position $x$ and let $U^+$ be obtained by the junction move, with index $g$.
\begin{enumerate}
\item[(a)] If $g=0$, then $U^+$ is a chain $u\tom{k+1}v^+\tom{k}w$ of $\mathcal M'^{\,\mathrm{mov}}_k$ of weight $q$, with triple $\big(q,\ u(\overline{N_1(U)}),\ w(\overline{M_2(U)})\big)$, where $\overline S=[1,k]\setminus S$.
\item[(b)] If $g\geq1$ and $x_g$ is not moved by $S_{e_{s+1}}\cdots S_{e_r}W_2'$, then $U^+$ is a compatible pair $u\tom{k}v^+\tom{k}w$ of weight $q$ which is the blocked head of its block, with $T(v^+)=\big(q,\ u(\overline{N_1(U)}),\ w(\overline{M_2(U)})\big)$.
\item[(c)] If $g\geq1$ and $x_g$ is moved by $S_{e_{s+1}}\cdots S_{e_r}W_2'$, let $f$ be the cycle of $S_{e_{s+1}}\cdots S_{e_r}W_2'$ through $x_g$ (so $f$ is either $(e_r,\ldots,e_{s+1},C)$ or a cycle of $W_2'$). Then the transposition word $U'$ obtained from $U^+$ by moving the cut one letter to the right, past $T_{x_g}$, and moving the letters of $f$ to the front of the suffix, is a state with active position $x_g$, with $N_1(U')=N_1(U)$ and $M_2(U')=M_2(U)$, and the index $j$ of its junction word is $g<j$.
\end{enumerate}
Consequently, starting from a source and iterating the junction move, one reaches after at most $j$ moves a target (a blocked head or a chain of $\mathcal M'^{\,\mathrm{mov}}_k$) of weight $q$ whose triple is the triple $S(v)$ of the source if the source is a blocked tail $v$, and the triple of the source if the source is a chain of $\mathcal M_k^{\mathrm{mov}}$. We write $\Phi(U)$ for this target.
\end{proposition}
\begin{proof}
We first collect what Lemma \ref{lemma:token}(ii) says about the positions $e_t$ which are among $x_1,\ldots,x_{j-1}$: if $t<s$ then $e_t=x_l$ with $l>g$, and if $t>s$ then $e_t=x_l$ with $l\leq g$. (The second assertion of (ii) applies because $JS_{e_{s+1}}\cdots S_{e_r}$ is a segment of the cover path $U$.) In particular $x_g\notin\{e_1,\ldots,e_{s-1}\}$, and $\{e_{s+1},\ldots,e_r\}\cap\{x_{g+1},\ldots,x_{j-1}\}=\emptyset$.

\emph{The prefix of $U^+$.} For $g\geq1$ it is $R_1$ followed by the pivot factorization of the cycle $(x_{j-1},\ldots,x_g,C)$ and the pivot factorization of $(x_{g-1},\ldots,x_1,K)$. The positions $x_1,\ldots,x_{j-1}$ avoid $\supp R_1$ by (S1), so if $R_1$ has no cycle with top $C$ these are three pairwise disjoint pieces, and if $R_1$ contains a cycle $(z_t,\ldots,z_1,C)$, then its letters followed by $S_{x_g}\cdots S_{x_{j-1}}$ form the pivot factorization of the cycle $(x_{j-1},\ldots,x_g,z_t,\ldots,z_1,C)$, disjoint from the other cycles. Thus the prefix is a pivot word of a product $\pi_1^+$ of pairwise disjoint cycles, each with exactly one element $>k$ (by (S1), the cycles of $R_1$ have this property, and the two new cycles have elements $x_l\leq k$ besides their tops), i.e.\ $\pi_1^+$ is $k$-shaped; being an initial segment of a cover path it is $u$-tight (Lemma \ref{lemma:tight}(iii)), so $u\tom{k}v^+:=u\pi_1^+$. The positions $\leq k$ it moves form $N_1(U)$. For $g=0$ the prefix is $R_1$ followed by $Z S_{x_1}\cdots S_{x_{j-1}}$, which merges with the cycle of $R_1$ with top $C$ (if any) into the cycle $(x_{j-1},\ldots,x_1,K,z_t,\ldots,z_1,C)$; so the prefix is a pivot word of a $(k+1)$-shaped, $u$-tight permutation moving $K$, i.e.\ $u\tom{k+1}v^+$ with $K$ moved, and it moves the positions $N_1(U)\cup\{K\}$.

\emph{The suffix of $U^+$.} Its letters $T$ are $T_{x_g}$ (if $g\geq1$), $T_{e_1},\ldots,T_{e_{s-1}}$, $T_x$ (if $x\leq k$), and the letters $T$ of $W_2'$, which form the pivot factorization of a cycle $D$ with top $K$ (trivial if $x=K$, since then $K\in e$ and $W_2'$ is disjoint from $e$ by (S2)). The positions $x_g,e_1,\ldots,e_{s-1},x$ and those of $D$ are pairwise distinct ($x_g\notin\{e_t\}_{t<s}$ by the first paragraph; $x\notin\{e_t\}_{t\neq s}$; $D$ is disjoint from $e$ and does not contain $x$ by (S2), and does not contain $x_g$ by (S3)), so these letters form, in the order in which they occur, the pivot factorization of the cycle $D^+=(D,\ x,\ e_{s-1},\ldots,e_1,\ x_g,\ K)$ (we list the elements following the cycle from the top), respectively $(D,x,e_{s-1},\ldots,e_1,K)$ when $g=0$ and $(D,e_{s-1},\ldots,e_1,x_g,K)$ in the $\mathcal M$ case; by (S2) all its elements other than $K$ are $\leq k$. The remaining letters of the suffix are those of $(e_r,\ldots,e_{s+1},C)$ and of the other cycles of $W_2'$, which are pairwise disjoint, disjoint from $D$, $x$ and $e_1,\ldots,e_{s-1}$, have, by (S2), exactly one element greater than $k$, and, under the hypothesis of (a) or (b), do not contain $x_g$. Hence in cases (a) and (b) the suffix is a pivot word of a $k$-shaped permutation $\pi_2^+$, which is $v^+$-tight because the suffix is a terminal segment of a cover path; so $v^+\tom{k}w$. The positions $\leq k$ moved by the suffix are $M_2(U)$ in case (a) and $M_2(U)\sqcup\{x_g\}$ in case (b).

(a) By the above, $U^+$ is a chain $u\tom{k+1}v^+\tom{k}w$ whose first step moves $K$; its weight is $q$ and its triple is $(q,\fx_{k+1}(u,v^+),\fx_k(v^+,w))=(q,u(\overline{N_1(U)}),w(\overline{M_2(U)}))$.

(b) We have $u\tom{k}v^+\tom{k}w$. The two steps are compatible: the cycles with top $K$ are $(x_{g-1},\ldots,x_1,K)$ and $D^+$, which meet only in $K$ (as $e_t=x_l$ with $t<s$ forces $l>g$); and by Lemma \ref{lemma:overlap} two cycles of the two steps with a common top $\geq k+2$ cannot share an element $\leq k$. The block data of the pair are $\big(\rho_1^+,\ X^+=(x_{g-1},\ldots,x_1,K)D^+,\ \rho_2^+\big)$ with cut $g-1$, where $\rho_1^+$ contains the cycle through $x_g$ with top $C$; since $x_g$ is the $g$-th pivot of $X^+$, $i_{\max}=g-1$, so the pair is the head of its block and this head is blocked. Finally $\fx_k(u,v^+)=u(\overline{N_1(U)})$, $\fx_k(v^+,w)=w(\overline{M_2(U)\sqcup\{x_g\}})$ and $c(v^+)=v^+(K)=\beta_{g-1}$, the value at $K$ after the prefix; the letter $T_{x_g}$ moves this value to $x_g$, which is not moved by any later letter, so $w(x_g)=c(v^+)$ and $T(v^+)=(q,u(\overline{N_1(U)}),w(\overline{M_2(U)}))$.

(c) Since $x_g\notin D$ by (S3), $f\neq D$: so $f$ is either $(e_r,\ldots,e_{s+1},C)$ or a cycle of $W_2'$ other than $D$, and in either case, by (S2), $f$ consists of its top, which is $\geq k+2$, together with positions $\leq k$, one of them $x_g$. The prefix of $U'$ is the prefix of $U^+$ followed by $T_{x_g}$: it has the form (S1) with $R_1'$ the prefix of $U^+$ minus its letters $T$ (a pivot word of the cycles of $R_1$ other than the one with top $C$, if any, together with $(x_{j-1},\ldots,x_g,z_t,\ldots,z_1,C)$; each of these consists of a top $\geq k+2$ and positions $\leq k$, by (S1) for $R_1$), with $K$-letters $T_{x_1},\ldots,T_{x_{g-1}},T_{x_g}$, and active position $x_g\in\supp R_1'$; the positions $x_1,\ldots,x_{g-1}$ avoid $\supp R_1'$. The letters of $f$ move no position moved by the letters of the suffix of $U^+$ which they pass (these letters are $T_{x_g}$---not passed---$T_{e_1},\ldots,T_{e_{s-1}},T_x$, the letters of $(e_r,\ldots,e_{s+1},C)$ if $f\neq$ this cycle, and letters of cycles of $W_2'$ other than $f$; all are disjoint from $f$ because the cycles of the suffix are pairwise disjoint apart from the common element $x_g$ of $T_{x_g}$ and $f$), so by Lemma \ref{lemma:inert} $U'$ is a cover path of weight $q$. Its suffix is $f$ followed by a pivot word of pairwise disjoint cycles: the cycle $D'=(D,x,e_{s-1},\ldots,e_1,K)$ formed by the letters $T_{e_1},\ldots,T_{e_{s-1}},T_x$ and those of $D$ (the letter $T_{x_g}$ having moved to the prefix), the cycle $(e_r,\ldots,e_{s+1},C)$ if $f$ is not this cycle, and the cycles of $W_2'$ other than $f$; all of these are disjoint from $f$, none of them contains $x_g$, and each consists of a top $\geq K$ together with positions $\leq k$. This is (S2) for $U'$, with $e'=f$ and $x'=x_g$. For (S3), the top-$K$ cycle of the new $W_2'$ is $D'$, and it avoids $x_1,\ldots,x_{g-1}$: $D$ does by (S3) for $U$, $x=x_j$ does, and $e_t=x_l$ with $t<s$ forces $l>g$ by the first paragraph. So $U'$ is a state with active position $x_g\leq k$, whose junction word has $g$ letters $T$ and index $g<j$. Its sets are $N_1(U')=\supp R_1'\cup\{x_1,\ldots,x_{g-1}\}=\supp R_1\cup\{x_1,\ldots,x_{j-1}\}=N_1(U)$ and $M_2(U')=M_2(U)\cup\{x_g\}=M_2(U)$, since $x_g\in\supp f\subseteq M_2(U)$.

For the consequence, note that a source $U$ has $N_1(U)=M_1\setminus\{x\}$ where $M_1$ is the set of positions $\leq k$ moved by its first step, and $M_2(U)=M_2$ is the set moved by its second step. If $U$ is a blocked tail $v$, then $S(v)=(q,\fx_k(u,v)\cup\{c(v)\},\fx_k(v,w))$ with $c(v)=u(x)$ (Lemma \ref{lemma:blocks}(iii)), which is $(q,u(\overline{M_1\setminus\{x\}}),w(\overline{M_2}))$; if $U$ is a chain of $\mathcal M_k^{\mathrm{mov}}$ its triple is $(q,\fx_k(u,v),\fx_{k+1}(v,w))=(q,u(\overline{M_1}),w(\overline{M_2}))$, and $N_1(U)=M_1$ as $x=K$. Since the junction move strictly decreases the index $j$ of the junction word, the iteration terminates, by (a) or (b), after at most $j$ moves, and the invariance of $N_1$ and $M_2$ in (c) gives the triple.
\end{proof}

\begin{example}[A junction move] \label{example:junction}
Let $k=1$, $K=2$, $u=123$ and $w=312$, so that $u^{-1}w=312=t_{12}t_{13}$. The steps $u\tom1v$ are the trivial one and $t_{12}$ ($v=213$; $t_{13}$ is not a cover of $123$ because of the intermediate value $2$), and the only compatible pair $u\tom1v\tom1w$ is
$$123\xrightarrow{t_{12}}213\xrightarrow{t_{13}}312,$$
with weight $q=1$ and $\omega_1=[\emptyset\mid\emptyset]=1$. Its block data are $\rho_1=1$, $X=(x_1,K)=(1,2)$, $\rho_2=(1,3)$, with cut $i=1=n$; since $x_1=1\in\supp\rho_2$ we have $i_{\min}=1$, so $v=213$ is a \emph{blocked tail}, and since $\supp\rho_1=\emptyset$, $i_{\max}=1$ and $v$ is also an unblocked head. Its triples are $S(v)=(1,\{u(1)\},\emptyset)=(1,\{1\},\emptyset)$ and $T(v)=(1,\emptyset,\{v(2)\})=(1,\emptyset,\{1\})$.

As a state, $v$ is the transposition word $U=T_1\,|\,S_1=t_{12}\,|\,t_{13}$ with $R_1$ empty, $(x_1)=(1)$, active position $x=1$, $e=(1,3)$ with $C=3$, $r=s=1$, and $W_2'$ empty; its junction word is $J=T_1S_1$ with $j=1$, so the only resolution is $J_0=Z\,T_1=t_{23}\,t_{12}$ (indeed $T_xS_x=ZT_x$), and $g=0$. From $\psi=u=123$ the word $J$ exchanges $1,2$ and then $2,3$, both ascents, and $J_0$ exchanges $2,3$ and then $1,3$, both ascents; so $\mathrm{wt}_u(J_0)=\mathrm{wt}_u(J)=1$, as Lemma \ref{lemma:token}(i) asserts, and $J_0$ is a cover path by Lemma \ref{lemma:degree}. The junction move gives
$$U^+=Z\ \big|\ T_1=t_{23}\ \big|\ t_{12},\qquad\text{i.e.}\qquad 123\xrightarrow{t_{23}}132\xrightarrow{t_{12}}312,$$
a chain $u\tom2v^+\tom1w$ whose first step moves $K=2$: a chain of $\mathcal M'^{\,\mathrm{mov}}_1$, as in Proposition \ref{prop:process}(a). Its triple is $\big(1,\fx_2(123,132),\fx_1(132,312)\big)=(1,\{1\},\emptyset)=S(v)$, as required. Here $N_1(U)=\emptyset$ and $M_2(U)=\{1\}$, and the triple is $(1,u(\overline{N_1}),w(\overline{M_2}))=(1,\{1\},\emptyset)$ in the notation of the proposition.

For completeness, the other side of \eqref{eq:Bk} for this $(u,w)$: $\mathcal M_1(u,w)$ consists of the single chain $123\xrightarrow{t_{12}}213\xrightarrow{t_{13}}312$ read with the second step at level $2$ (it does not move $K$, so it is not in $\mathcal M_1^{\mathrm{mov}}$), with triple $(1,\emptyset,\fx_2(213,312))=(1,\emptyset,\{1\})=T(v)$, matching the unblocked head as in the proof of Proposition \ref{prop:junction}. Thus \eqref{eq:junction} reads $\{S(v)\}=\{(1,\{1\},\emptyset)\}$ on both sides, \eqref{eq:Bk} reads $\{T(v)\}\uplus\{(1,\{1\},\emptyset)\}=\{S(v)\}\uplus\{(1,\emptyset,\{1\})\}$, and Proposition \ref{prop:giambelliform} reads $(z_1-z_2)\cdot1=(y_1-z_2)-(y_1-z_1)$. The reversal of Lemma \ref{lemma:reversal} sends this configuration to the pair $w_0w=132\tom1 w_0v=231\tom1 w_0u=321$ for $(u,w)=(132,321)$, which is a blocked head, and the chain $U^+$ to the chain $132\xrightarrow{t_{12}}312\xrightarrow{t_{23}}321$ of $\mathcal M_1^{\mathrm{mov}}$ for that pair.
\end{example}

\subsubsection*{Injectivity}

\begin{lemma} \label{lemma:reversestate}
Let $U$ be a state and $U^+$ the result of the junction move, and work in $S_n$. Then the reversed transposition word $(U^+)^{\mathrm{rev}}$, read from $w_0w$ with the reversed cut, is (after reorderings as in Lemma \ref{lemma:inert}) a state for the pair $(w_0w,w_0u)$ and the weight $q$, with active position $x_g$ if $g\geq1$ and $K$ if $g=0$, and the junction move applied to it produces $U^{\mathrm{rev}}$.
\end{lemma}
\begin{proof}
By Lemma \ref{lemma:reversal}, $(U^+)^{\mathrm{rev}}$ and $U^{\mathrm{rev}}$ are cover paths from $w_0w$ of weight $q$ with product $w^{-1}u$. Let $g\geq1$. The prefix of $(U^+)^{\mathrm{rev}}$ is the reversed suffix of $U^+$: the reversal of the pivot word of $W_2'$ and of $S_{e_{s+1}}\cdots S_{e_r}$, followed by $T_x\,T_{e_{s-1}}\cdots T_{e_1}\,T_{x_g}$. Its letters $T$ are those of $D$ (reversed), then $T_x,T_{e_{s-1}},\ldots,T_{e_1},T_{x_g}$; the other letters move positions disjoint from those of $D$, $x$ and $e_1,\ldots,e_{s-1}$, so by Lemma \ref{lemma:inert} they may be moved to the front, giving the form (S1) with $R_1^{\mathrm{rev}}$ consisting of the reversed pivot words of the cycles of $W_2'$ other than $D$ and of $(e_r,\ldots,e_{s+1},C)$, the $K$-letters being $T_{y_1},\ldots,T_{y_{p}}$ with $(y_1,\ldots,y_p)=(D\text{ reversed},x,e_{s-1},\ldots,e_1,x_g)$, and active position $x_g$ (which may lie in $\supp R_1^{\mathrm{rev}}$). The suffix of $(U^+)^{\mathrm{rev}}$ is $T_{x_{g-1}}\cdots T_{x_1}\,S_{x_{j-1}}\cdots S_{x_g}\,R_1^{\mathrm{rev}}$; here $S_{x_{j-1}}\cdots S_{x_g}$ followed by the reversed letters of the cycle of $R_1$ with top $C$ (if any) is the pivot factorization of a cycle $e^{\mathrm{rev}}$ with top $C$ through $x_g$, in which $x_g$ is preceded by $x_{j-1},\ldots,x_{g+1}$; moving $T_{x_{g-1}}\cdots T_{x_1}$ past it (Lemma \ref{lemma:inert}) gives the form (S2), the remaining cycles being $(x_1,\ldots,x_{g-1},K)$ and the other cycles of $R_1$, which do not contain $x_g$. The shape requirements of (S1) and (S2) are inherited from those of $U$: the cycles of $R_1^{\mathrm{rev}}$ are cycles of $W_2'$ other than $D$, which by (S2) for $U$ have tops $\geq k+2$ (their tops are $\geq K$ and differ from $K$) and other elements $\leq k$, together with $(e_r,\ldots,e_{s+1},C)$; the $K$-letters $y_1,\ldots,y_{p-1}$ are positions $\leq k$ not in $\supp R_1^{\mathrm{rev}}$ (the cycles of $W_2'$ are pairwise disjoint and disjoint from $e$, and $x\notin\supp W_2'$); the cycle $e^{\mathrm{rev}}$ has top $C\geq k+2$ and other elements $x_{j-1},\ldots,x_{g+1}$ and the $z$'s of the top-$C$ cycle of $R_1$, all $\leq k$; and the remaining cycles $(x_1,\ldots,x_{g-1},K)$ and those of $R_1$ consist of a top $\geq K$ and positions $\leq k$ by (S1) for $U$. Finally (S3) holds for $(U^+)^{\mathrm{rev}}$: its top-$K$ cycle $(x_1,\ldots,x_{g-1},K)$ avoids $y_1,\ldots,y_{p-1}$, because $D$ avoids $x_1,\ldots,x_{j-1}$ by (S3) for $U$, $x=x_j$, and $e_t=x_l$ with $t<s$ forces $l>g$ (Lemma \ref{lemma:token}(ii)). So $(U^+)^{\mathrm{rev}}$ is a state. Its junction word is $T_{y_1}\cdots T_{y_p}\,S_{x_{j-1}}\cdots S_{x_{g+1}}\,S_{x_g}$. Suppose first that $x\leq k$. The resolution with stop at the letter $T_x$ (the letter $y_m$ with $y_m=x$) is
$$S_x\,S_{e_{s-1}}\cdots S_{e_1}\ \cdot\ (D\text{ reversed})\ \big|\ T_x\,T_{x_{j-1}}\cdots T_{x_{g+1}}\,T_{x_g}\ \cdot\ (\text{rest of the suffix}),$$
which, after moving the block $S_xS_{e_{s-1}}\cdots S_{e_1}$ past the letters of $D$ and the reversed letters of the cycle of $R_1$ with top $C$ past $T_{x_{g-1}}\cdots T_{x_1}$ (both inert), is exactly
$$U^{\mathrm{rev}}=\big(W_2'\big)^{\mathrm{rev}}\,S_{e_r}\cdots S_{e_{s+1}}\,S_xS_{e_{s-1}}\cdots S_{e_1}\ \big|\ T_xT_{x_{j-1}}\cdots T_{x_1}\,R_1^{\mathrm{rev}}.$$
Since $U^{\mathrm{rev}}$ is a cover path of weight $q$, Lemma \ref{lemma:token}(i) shows that this resolution is the one selected by the junction move. If $x=K$, there is no letter $T_x$, the sequence $(y_1,\ldots,y_p)$ is $(e_{s-1},\ldots,e_1,x_g)$, and the same comparison shows that the resolution $J_0=Z\,S_{e_{s-1}}\cdots S_{e_1}\cdot T_{x_{j-1}}\cdots T_{x_g}$ of the junction word yields $U^{\mathrm{rev}}=(W_2')^{\mathrm{rev}}S_{e_r}\cdots S_{e_{s+1}}ZS_{e_{s-1}}\cdots S_{e_1}\,|\,T_{x_{j-1}}\cdots T_{x_1}R_1^{\mathrm{rev}}$ up to inert reorderings. The case $g=0$ is identical, with the prefix of $(U^+)^{\mathrm{rev}}$ having no active letter after $T_{e_1}$, the suffix containing the cycle $(x_{j-1},\ldots,x_1,K,\ldots,C)$ through $K$, so that the reversed state is of $\mathcal M$ type (its $e^{\mathrm{rev}}$ has top $C\geq k+2$, contains $K$, and has all other elements $\leq k$; its $W_2'$ consists of the cycles of $R_1$ other than the one with top $C$, which by (S1) have tops $\geq k+2$ and positions $\leq k$, and (S3) is vacuous), and the resolution with stop at $T_x$ (resp.\ $J_0$ if $x=K$) again giving $U^{\mathrm{rev}}$.
\end{proof}

\begin{proposition} \label{prop:injective}
The map $\Phi$ from the sources of weight $q$ for $(u,w)$ to the targets of weight $q$ for $(u,w)$ is injective.
\end{proposition}
\begin{proof}
Let $U_0$ be a source and $U_0,U_1,\ldots,U_m$ the states traversed by the process, $U_{i+1}$ being obtained from $U_i$ by a junction move followed, in case (c), by the normalization of Proposition \ref{prop:process}(c), and $U_m^+=\Phi(U_0)$ the final target. We show that $U_i$ is determined by $U_{i+1}$ (resp.\ by $U_m^+$ for $i=m$). The transposition word $U_i^+$ (the result of the move, before normalization) is determined by $U_{i+1}$: one moves the active letter $T$ of $U_{i+1}$, the last letter of its prefix, back to the front of the suffix (the normalization of (c) only reorders inert letters otherwise). By Lemma \ref{lemma:reversestate}, $(U_i^+)^{\mathrm{rev}}$ is a state whose junction move produces $U_i^{\mathrm{rev}}$; since the junction move is deterministic (Lemma \ref{lemma:token}(i) selects a unique resolution), $U_i^{\mathrm{rev}}$, hence $U_i$, is determined by $U_i^+$. Finally, whether $U_i$ is the source or an intermediate state is read off from $U_i$ itself: a state is a source if and only if its active position does not lie in $\supp R_1$. Hence the whole sequence, and in particular $U_0$, is determined by $\Phi(U_0)$.
\end{proof}

\begin{theorem}[The junction identity] \label{theorem:junction}
The identity \eqref{eq:junction} holds for all $k\geq1$ and all $u,w\in S_\infty$. Consequently \eqref{eq:Bk}, Theorem \ref{theorem:doubledpieri} and \eqref{eq:Jk} hold for all $k\geq1$.
\end{theorem}
\begin{proof}
Fix $k$, $u$, $w$, a monomial $q$, and $n$ so large that everything takes place in $S_n$. Let $\mathrm{Src}(u,w)$ be the multiset of triples of the sources of weight $q$ for $(u,w)$, i.e.\ $S(v)$ for the blocked tails $v$ and the triples of the chains of $\mathcal M_k^{\mathrm{mov}}$, and $\mathrm{Tgt}(u,w)$ the multiset of $T(v)$ for the blocked heads and of the triples of the chains of $\mathcal M'^{\,\mathrm{mov}}_k$, all of weight $q$. By Propositions \ref{prop:process} and \ref{prop:injective}, $\Phi$ is an injective, triple-preserving map from sources to targets; hence $\mathrm{Src}(u,w)\subseteq\mathrm{Tgt}(u,w)$ as multisets.

Now apply the reversal of Lemma \ref{lemma:reversal} to whole configurations: it sends a compatible pair $u\tom{k}v\tom{k}w$ to the pair $w_0w\tom{k}w_0v\tom{k}w_0u$ (the reversed pivot words are pivot words of the inverse permutations, which are tight by Lemma \ref{lemma:tight}(iii)), with the same weight, with block data $(\rho_2^{-1},X^{-1},\rho_1^{-1})$ and cut $n-i$, so that $i_{\min}$ and $i_{\max}$ are exchanged with $n-i_{\max}$ and $n-i_{\min}$: blocked tails go to blocked heads and conversely. Likewise it exchanges the chains of $\mathcal M_k^{\mathrm{mov}}$ for $(u,w)$ with those of $\mathcal M'^{\,\mathrm{mov}}_k$ for $(w_0w,w_0u)$. On triples it acts by $\sigma(q,A,B)=(q,w_0B,w_0A)$, where $w_0$ acts on values: indeed $\fx_k(w_0w,w_0v)=w_0\fx_k(v,w)$, $\fx_k(w_0v,w_0u)=w_0\fx_k(u,v)$ and $(w_0v)(K)=w_0(v(K))$, so that $S(w_0v)=\sigma(T(v))$ and $T(w_0v)=\sigma(S(v))$, and similarly for the mixed chains. Hence $\sigma(\mathrm{Tgt}(u,w))=\mathrm{Src}(w_0w,w_0u)$ and $\sigma(\mathrm{Src}(u,w))=\mathrm{Tgt}(w_0w,w_0u)$. Applying the first paragraph to $(w_0w,w_0u)$ gives $\mathrm{Src}(w_0w,w_0u)\subseteq\mathrm{Tgt}(w_0w,w_0u)$, i.e.\ $\mathrm{Tgt}(u,w)\subseteq\mathrm{Src}(u,w)$. Therefore $\mathrm{Src}(u,w)=\mathrm{Tgt}(u,w)$ for every $q$, which is \eqref{eq:junction}. The consequences follow from Propositions \ref{prop:junction} and \ref{prop:doubledreduction}.
\end{proof}

\section{Proof of the parabolic reduction} \label{section:pw}

In this section we prove Proposition \ref{prop:parabolic} in finite rank. Peterson's comparison formula expresses the Gromov--Witten invariants of $G/P$ in terms of those of $G/B$; it was proved by Woodward \cite{woodward2005d} and extended to equivariant quantum cohomology by Huang and Li \cite{huang2015equivariant}, whose form we use. We use the notation of Section \ref{section:definitions}: $W_P$ is a parabolic subgroup of $S_N$ generated by $P\subseteq\{s_1,\ldots,s_{N-1}\}$, $\omega_P=w_P$ is its longest element, and for $\lambda_P\in Q^\vee/Q^\vee_P$ we write $\lambda_B\in Q^\vee$ for its Peterson--Woodward lift. In the language of roots, $\lambda_B$ is the unique lift with $\langle\alpha,\lambda_B\rangle\in\{0,-1\}$ for all positive roots $\alpha$ of the Levi subgroup of $P$; in type $A$ these are the roots $\alpha_i+\cdots+\alpha_{j-1}$ with $s_i,\ldots,s_{j-1}\in P$, and $\langle\alpha_i+\cdots+\alpha_{j-1},d\rangle=\langle d\rangle_{ij}$, so this is the condition of Definition \ref{def:pwlift}.

\begin{proposition}[{\cite[Proposition~2.10]{huang2015equivariant}}] \label{prop:huangli}
Let $W_P$ be a parabolic subgroup of $W$ and $u,v,w\in W^P$. Then for any $\lambda_P\in Q^\vee/Q_P^\vee$, we have that
$$c_{u,v}^{w,\lambda_P}(y;y) = c_{u,v}^{w\omega_P\omega_{P'},\lambda_B}(y;y),$$
where $\omega_P$ is the longest element in $W_P$, $\omega_{P'}$ is the longest element in $W_{P'}$, where $P'$ is the parabolic subgroup generated by $\{s_\alpha\mid \alpha\in\Delta_P,\langle \alpha, \lambda_B\rangle=0\}$, and $\lambda_B$ is the Peterson-Woodward lift of $\lambda_P$.
\end{proposition}

Here the left-hand side is a structure constant of $QH_T^*(G/P)$, i.e.\ of the parabolic quantum double Schubert polynomials $\qsch^P_w(x;y)$ of Definition \ref{def:parabolic}, which represent the Schubert classes of $QH_T^*(G/P)$ \cite{lam2014quantum}, and the right-hand side is a structure constant of $QH_T^*(G/B)$.

\begin{corollary} \label{corollary:pwz}
With the notation of Proposition \ref{prop:huangli}, for all $u,v,w\in W^P$,
$$c_{u,v}^{w,\lambda_P}(y;z) = c_{u,v}^{w\omega_P\omega_{P'},\lambda_B}(y;z).$$
\end{corollary}
\begin{proof}
The parabolic polynomials satisfy the same quantum Cauchy identity as in Lemma \ref{lemma:quantumtoordinary}. Indeed, by \cite[Lemma~18]{lam2014quantum}, $\qsch^P_w(x;y)=\theta^P(\sch_w(x;y))$ for $w\in W^P$, where $\theta^P$ is the parabolic quantization map of Ciocan-Fontanine \cite{ciocan1999partial}, extended linearly over the coefficient variables; applying $\theta^P$ to the classical identity $\sch_v(x;z)=\sum_{v_1\cdot v_2=v}\sch_{v_1}(y;z)\sch_{v_2}(x;y)$ gives $\qsch^P_v(x;z)=\sum_{v_1\cdot v_2=v}\sch_{v_1}(y;z)\qsch^P_{v_2}(x;y)$, where we note that $v_2\in W^P$ for every reduced factorization (see below), so that $\theta^P(\sch_{v_2}(x;y))$ is indeed $\qsch^P_{v_2}(x;y)$. Hence
$$c_{u,v}^{w,\lambda_P}(y;z) = \sum_{v_1\cdot v_2=v}c_{u,v_2}^{w,\lambda_P}(y;y)\,\sch_{v_1}(y;z).$$
Here, if $v=v_1\cdot v_2$ is reduced and $v\in W^P$, then $v_2\in W^P$: for $s_i\in P$, $v_2s_i<v_2$ would give $\ell(vs_i)\leq\ell(v_1)+\ell(v_2s_i)<\ell(v)$, contradicting $v\in W^P$. Thus Proposition \ref{prop:huangli} applies to each term, and resumming by Lemma \ref{lemma:quantumtoordinary} gives
$$\sum_{v_1\cdot v_2=v}c_{u,v_2}^{w\omega_P\omega_{P'},\lambda_B}(y;y)\,\sch_{v_1}(y;z) = c_{u,v}^{w\omega_P\omega_{P'},\lambda_B}(y;z).\qedhere$$
\end{proof}

\begin{proof}[Proof of Proposition \ref{prop:parabolic}]
Let $d_P\in Q^\vee/Q^\vee_P$ with Peterson--Woodward lift $d=\lambda_B$. In type $A$, $\langle\alpha_i,d\rangle=2d_i-d_{i-1}-d_{i+1}=\langle d\rangle_i$, so the parabolic subgroup $P'$ of Proposition \ref{prop:huangli} is generated by $\{s_i\in P\mid\langle d\rangle_i=0\}$, that is, by $\{s_i\mid i\in z_P(d)\}$, and $\omega_P\omega_{P'}=w_Pw_{z_P(d)}$. The claim is therefore Corollary \ref{corollary:pwz}. Its application to stabilized polynomial coefficients is explained immediately after Proposition \ref{prop:parabolic}.
\end{proof}

\appendix
\section{Proof of the positivity theorem} \label{appendix:positivity}
\begin{center}
\textsc{by Tongren Xiao and Yihua Yue}\footnote{School of Mathematics, Sun Yat-sen University, Guangzhou, China}
\end{center}

This appendix proves Theorem \ref{theorem:positivity}. The argument combines the geometric argument of Mihalcea \cite{mihalcea2006positivity} with the Graham-positivity theorem of Gao and Xiong \cite{gao2025grahampositivitytripleschubert} for $B^-$-stable effective cycles.

\begin{proof}[Proof of Theorem \ref{theorem:positivity}]
Our product is written with the output basis carrying the parameters of the second factor, so commutativity of multiplication identifies the coefficient $C_{v,u}^{w,d}(z,y)$ below with $c_{u,v}^{w,d}(y;z)$. It is therefore enough to prove that $C_{u,v}^{w,d}(y,t)$ is in $\mathbb{Z}_{\geq 0}[t_i-y_j]$.

Fix $u,v,w\in S_n$, increasing $n$ if necessary so that all permutations and parameters occurring in the coefficient under consideration are included. Write
$\mathbf a=(t_1,\ldots,t_n,y_1,\ldots,y_n)$, and let $\tau\in S_{2n}$ interchange the two blocks, so that $\tau\mathbf a=(y_1,\ldots,y_n,t_1,\ldots,t_n)$. Regard the product expansion in $S_{2n}$:
\begin{equation}\label{eq:xiaoyue-product}
\qsch_u(x;\tau\mathbf a)\qsch_v(x;\mathbf a)
=\sum_{w\in S_{2n}}C_{u,v}^{w,\tau}(\mathbf a;q)\qsch_w(x;\mathbf a).
\end{equation}
We recall the geometric interpretation of these coefficients. Let $e_i(x)$ and $E_i(x,q)$ denote, respectively, the elementary symmetric polynomials and their quantum analogues in $2n$ variables, and set
\[
J_{\mathbf a}^{2n}=\langle e_i(x)-e_i(\mathbf a):1\leq i\leq 2n\rangle,\qquad
J_{q\mathbf a}^{2n}=\langle E_i(x,q)-e_i(\mathbf a):1\leq i\leq 2n\rangle.
\]
The Fomin--Gelfand--Postnikov deformation map \cite{fomin1997quantum}, extended $\mathbb Z[\mathbf a]$-linearly, induces
\[
\overline{\mathcal Q}:\mathbb Z[x,\mathbf a]/J_{\mathbf a}^{2n}
\longrightarrow \mathbb Z[x,\mathbf a,q]/J_{q\mathbf a}^{2n}.
\]
The classical Borel presentation identifies the source with $H_T^*(X)$, and Kim's presentation \cite{kim1996equivariant} identifies the target with $QH_T^*(X)$, where $X=\mathrm{Fl}_{2n}$. Under the Lam--Shimozono identification \cite{lam2014quantum} of quantum double Schubert polynomials with quantum Schubert classes, \eqref{eq:xiaoyue-product} is the quantum equivariant Littlewood--Richardson expansion for the classes obtained by applying $\overline{\mathcal Q}$ to the first factor. Consequently,
\begin{equation}\label{eq:xiaoyue-eqlr}
C_{u,v}^{w,\tau}(\mathbf a;q)
=\sum_{d\geq0}q^d(\pi_{M_d})^T_*
\left(
\operatorname{ev}_1^*[\tau X_u]^T\,
\operatorname{ev}_2^*[X_v]^T\,
\operatorname{ev}_3^*[Y_w]^T
\right).
\end{equation}
Here $M_d=\overline M_{0,3}(X,d)$, $\operatorname{ev}_i:M_d\to X$ are the evaluation maps, $\pi_{M_d}:M_d\to\mathrm{pt}$ is the structure map, and $X_r$ and $Y_r$ denote the Schubert and opposite Schubert varieties indexed by $r$.

For each effective degree $d$, put
\[
F=(\operatorname{ev}_1,\operatorname{ev}_2):M_d\longrightarrow X\times X,\qquad
I_{u,v}^d=F^{-1}(\tau X_u\times X_v),\qquad
\Gamma_{u,v}^d=(\operatorname{ev}_3)_*[I_{u,v}^d].
\]
The incidence scheme $I_{u,v}^d$ is reduced and locally irreducible, and all its irreducible components have codimension $\ell(u)+\ell(v)$. Indeed, the moduli space $M_d$ is irreducible by Thomsen's theorem \cite{thomsen1998irreducibility}, and Mihalcea's incidence lemma \cite[Lemma~6.3]{mihalcea2006positivity} applies after diagonal translation by $a_0=1\times w_0\in S_{2n}$, where $w_0$ is the longest element in the second block: it identifies $\tau X_u\times X_v$ with $w_0^{(2n)}X_u\times X_v$. The codimension is therefore
\[
\operatorname{codim}_X(w_0^{(2n)}X_u)+\ell(v)=\ell(u)+\ell(v).
\]
By the incidence formula of Mihalcea \cite[Section~3.1]{mihalcea2006positivity}, this gives
\begin{equation}\label{eq:xiaoyue-incidence}
\operatorname{ev}_1^*[\tau X_u]^T\,
\operatorname{ev}_2^*[X_v]^T\cap[M_d]^T=[I_{u,v}^d]_T.
\end{equation}

The cycle $\Gamma_{u,v}^d$ is effective and $B^-(\tau)$-invariant. To see invariance, write $N^-(\tau)=N^-\cap\tau N^-\tau^{-1}$ and $B^-(\tau)=TN^-(\tau)$. The varieties $X_v$ and $\tau X_u$ are stable under the corresponding actions, so their product, its inverse image under the equivariant map $F$, and its proper push-forward under $\operatorname{ev}_3$ are $B^-(\tau)$-invariant. Effectivity follows from the reducedness of $I_{u,v}^d$.

There is a unique expansion
\[
[\Gamma_{u,v}^d]_T=\sum_{r\in S_{2n}}p_{u,v}^{r,d}(\mathbf a)[X_r]_T
\quad\text{in }H_T^*(X).
\]
By Gao and Xiong's Graham--positivity theorem for $B^-(\tau)$-invariant effective cycles \cite[Theorem~2.4 and Corollary~2.5]{gao2025grahampositivitytripleschubert},
\[
p_{u,v}^{r,d}(\mathbf a)\in
\mathbb Z_{\geq0}[t_i-y_j:1\leq i,j\leq n].
\]
Finally, by equivariant Poincaré duality, \eqref{eq:xiaoyue-incidence}, and the projection formula,
\[
(\operatorname{ev}_3)^T_*
\left(\operatorname{ev}_1^*[\tau X_u]^T
\operatorname{ev}_2^*[X_v]^T\right)
=[\Gamma_{u,v}^d]_T,
\]
and hence the degree-$d$ summand of \eqref{eq:xiaoyue-eqlr} is
\[
(\pi_X)^T_*\bigl([\Gamma_{u,v}^d]_T[Y_w]^T\bigr)
=\sum_{r\in S_{2n}}p_{u,v}^{r,d}(\mathbf a)
(\pi_X)^T_*\bigl([X_r]^T[Y_w]^T\bigr)
=p_{u,v}^{w,d}(\mathbf a).
\]
Thus every coefficient $C_{u,v}^{w,d}(y,t)$ is in
$\mathbb Z_{\geq0}[t_i-y_j]$. Replacing $(u,v,y,t)$ by $(v,u,z,y)$ proves the asserted positivity of $c_{u,v}^{w,d}(y;z)$.
\end{proof}

\bibliographystyle{acm}
\bibliography{formschub}

@article{huang2015equivariant,
  title={{On equivariant quantum Schubert calculus for $G/P$}},
  author={Huang, Yongdong and Li, Changzheng},
  journal={{Journal of Algebra}},
  volume={441},
  pages={21--56},
  year={2015},
  publisher={Elsevier}
}

@article{kirillov2000quantum,
	title={{Quantum double Schubert polynomials, quantum Schubert polynomials and Vafa--Intriligator formula}},
	author={Kirillov, Anatol N and Maeno, Toshiaki},
	journal={{Discrete Mathematics}},
	volume={217},
	number={1-3},
	pages={191--223},
	year={2000},
	publisher={Elsevier}
}

@article{bjorner1997shellable,
  title={Shellable nonpure complexes and posets. {I}{I}},
  author={Bj{\"o}rner, Anders and Wachs, Michelle},
  journal={Transactions of the American Mathematical Society},
  volume={349},
  number={10},
  pages={3945--3975},
  year={1997}
}

@article{reading2006cambrian,
  title={{C}ambrian lattices},
  author={Reading, Nathan},
  journal={Advances in Mathematics},
  volume={205},
  number={2},
  pages={313--353},
  year={2006}
}

@article{reading2004congruences,
  title={Lattice congruences of the weak order},
  author={Reading, Nathan},
  journal={Order},
  volume={21},
  number={4},
  pages={315--344},
  year={2004}
}

@article{fomin1997quantum,
	title={{Quantum Schubert polynomials}},
	author={Fomin, Sergey and Gelfand, Sergei and Postnikov, Alexander},
	journal={{Journal of the American Mathematical Society}},
	pages={565--596},
	year={1997},
	publisher={JSTOR}
}

@incollection{postnikov1999quantum,
	title={{On a quantum version of Pieri’s formula}},
	author={Postnikov, Alexander},
	booktitle={Advances in geometry},
	pages={371--383},
	year={1999},
	publisher={Springer}
}

@article{benedetti2025quantum,
	title={{A quantum Murnaghan--Nakayama rule for the flag manifold}},
	author={Benedetti, Carolina and Bergeron, Nantel and Colmenarejo, Laura and Saliola, Franco and Sottile, Frank},
	journal={{Algebraic Combinatorics}},
	volume={8},
	number={3},
	pages={619--653},
	year={2025},
	doi={10.5802/alco.427}
}

@article{colmenarejo2026quantum,
	title={{The quantum $k$-Bruhat order}},
	author={Colmenarejo, Laura and Mayers, Nicholas},
	journal={{arXiv preprint arXiv:2601.03437}},
	year={2026}
}

@article{leung2012classical,
	title={Classical aspects of quantum cohomology of generalized flag varieties},
  author={Leung, Naichung Conan and Li, Changzheng},
  journal={{International Mathematics Research Notices}},
  volume={2012},
  number={16},
  pages={3706--3722},
  year={2012},
  publisher={OUP}
}

@article{lam2014quantum,
	title={{Quantum double Schubert polynomials represent Schubert classes}},
	author={Lam, Thomas and Shimozono, Mark},
	journal={{Proceedings of the American Mathematical Society}},
	volume={142},
	number={3},
	pages={835--850},
	year={2014}
}

@article{woodward2005d,
	title={{On D. Peterson's comparison formula for Gromov-Witten invariants of G/P}},
	author={Woodward, Christopher T},
	journal={Proceedings of the American Mathematical Society},
	pages={1601--1609},
	year={2005},
	publisher={JSTOR}
}

@article{fulton2004quantum,
	title={{On the quantum product of Schubert classes}},
	author={Fulton, William and Woodward, Christopher},
	journal={{Journal of Algebraic Geometry}},
	volume={13},
	number={4},
	pages={641--661},
	year={2004}
}

@misc{schubmult,
	author = {Matthew J. Samuel},
	title = {{\texttt{schubmult}: a {P}ython package for multiplying {S}chubert polynomials}},
	howpublished = {Available from the {P}ython {P}ackage {I}ndex, \url{https://pypi.org/project/schubmult/}, with sources at \url{https://github.com/matthematics/schubmult}},
	year = {2026}
}

@book{humphreys1992reflection,
  title={Reflection groups and {C}oxeter groups},
  author={Humphreys, James E},
  year={1992},
  publisher={Cambridge university press}
}

@phdthesis{samuelleibniz, 
	author = "Matthew J. Samuel",
	school = "Rutgers, The State University of New Jersey",
	year = "2014",
	title = "The {L}eibniz formula for divided difference operators associated to {K}ac-{M}oody root systems"
}

@BOOK{combcox,
  TITLE = {Combinatorics of {C}oxeter {G}roups},
  SUBTITLE = {Graduate Texts in Mathematics 231},
  AUTHOR = {Anders Bj{\"o}rner and Francesco Brenti},
  YEAR = {2005}, 
  PUBLISHER = {Springer},
	PLACE = "New York"
}

@article{samuelmolev,
    author =       "Samuel, Matthew J.",
    title =        "A {M}olev-{S}agan type formula for double {S}chubert polynomials",
    journal =      "J. Pure Appl. Algebra",
    year =         "2024",
		volume = "228",
		number = "7",
		pages = "107636"
}

@book{notes,
    author    = "Macdonald, I. G.",
    title     = "Notes on {S}chubert {P}olynomials",
    year      = "1991",
    address = "D{\'e}p. de math{\'e}matiques et d{'}informatique",
    publisher   = "Universit{\'e} Du Qu{\'e}bec {\`a} Montr{\'e}al"
}

@article{sottile,
		author = "Sottile, Frank",
		title = "Pieri{'}s rule for flag manifolds and {S}chubert polynomials",
		year = "1996",
		journal = "Ann. Inst. Fourier",
		volume = "46",
		number = "1",
		pages = "89--110"
}

@article{robinsonpieri,
  title={A {P}ieri-{T}ype {F}ormula for {$H^*T (SL_n (C)/B)$}},
  author={Robinson, Shawn},
  journal={J. Algebra},
  volume={249},
  number={1},
  pages={38--58},
  year={2002},
  publisher={Elsevier}
}

@article{buch2016puzzle,
  title={The puzzle conjecture for the cohomology of two-step flag manifolds},
  author={Buch, Anders Skovsted and Kresch, Andrew and Purbhoo, Kevin and Tamvakis, Harry},
  journal={Journal of algebraic combinatorics},
  volume={44},
  pages={973--1007},
  year={2016},
  publisher={Springer}
}

@article{buch2015mutations,
  title={Mutations of puzzles and equivariant cohomology of two-step flag varieties},
  author={Buch, Anders Skovsted},
  journal={Ann. of Math.},
  pages={173--220},
  year={2015},
  publisher={JSTOR}
}

@article{molev1999littlewood,
  title={A {L}ittlewood-{R}ichardson rule for factorial {S}chur functions},
  author={Molev, Alexander and Sagan, Bruce},
  journal={Trans. Amer. Math. Soc.},
  volume={351},
  number={11},
  pages={4429--4443},
  year={1999}
}

@article{gao2025grahampositivitytripleschubert,
	title={Positivity in triple {S}chubert calculus},
	author={Gao, Yibo and Xiong, Rui},
	journal={S{\'e}minaire Lotharingien de Combinatoire},
	volume={95B},
	note={Proceedings of FPSAC 2026, Art.~37, 11 pp.},
	year={2026},
	url={https://radon.mat.univie.ac.at/~slc/wpapers/FPSAC2026/37.html}
}

@article{kim1996equivariant,
	title={On equivariant quantum cohomology},
	author={Kim, Bumsig},
	journal={International Mathematics Research Notices},
	number={17},
	pages={841--851},
	year={1996}
}

@article{mihalcea2006positivity,
	title={Positivity in equivariant quantum {S}chubert calculus},
	author={Mihalcea, Leonardo Constantin},
	journal={American Journal of Mathematics},
	volume={128},
	number={3},
	pages={787--803},
	year={2006}
}

@article{thomsen1998irreducibility,
	title={Irreducibility of $M_{0,n}(G/P,\beta)$},
	author={Thomsen, Jesper Funch},
	journal={International Journal of Mathematics},
	volume={9},
	number={3},
	pages={367--376},
	year={1998}
}

@article{mihalcea2007equivariant,
	title={{On equivariant quantum cohomology of homogeneous spaces: Chevalley formulae and algorithms}},
	author={Mihalcea, Leonardo Constantin},
	journal={{Duke Mathematical Journal}},
	volume={140},
	number={2},
	pages={321--350},
	year={2007}
}

@article{ciocan1999partial,
	title={{On quantum cohomology rings of partial flag varieties}},
	author={Ciocan-Fontanine, Ionu\c{t}},
	journal={{Duke Mathematical Journal}},
	volume={98},
	number={3},
	pages={485--524},
	year={1999}
}

@article{buch2003qmonk,
  title={A direct proof of the quantum version of {M}onk's formula},
  author={Buch, Anders Skovsted},
  journal={Proceedings of the American Mathematical Society},
  volume={131},
  number={7},
  pages={2037--2042},
  year={2003}
}

@article{buch2005partial,
  title={Quantum cohomology of partial flag manifolds},
  author={Buch, Anders Skovsted},
  journal={Transactions of the American Mathematical Society},
  volume={357},
  pages={443--458},
  year={2005}
}

@article{buch2003gwgrass,
  title={{G}romov--{W}itten invariants on {G}rassmannians},
  author={Buch, Anders Skovsted and Kresch, Andrew and Tamvakis, Harry},
  journal={Journal of the American Mathematical Society},
  volume={16},
  pages={901--915},
  year={2003}
}

@article{buch2015curve,
  title={Curve neighborhoods of {S}chubert varieties},
  author={Buch, Anders Skovsted and Mihalcea, Leonardo C.},
  journal={Journal of Differential Geometry},
  volume={99},
  pages={255--283},
  year={2015}
}

@article{buch2011qkgrass,
  title={Quantum {K}-theory of {G}rassmannians},
  author={Buch, Anders Skovsted and Mihalcea, Leonardo Constantin},
  journal={Duke Mathematical Journal},
  volume={156},
  number={3},
  pages={501--538},
  year={2011}
}

@article{mihalcea2006equivariant,
  title={Equivariant quantum {S}chubert calculus},
  author={Mihalcea, Leonardo Constantin},
  journal={Advances in Mathematics},
  volume={203},
  number={1},
  pages={1--33},
  year={2006}
}

@article{gu2020rimhook,
	title={{A rim-hook rule for quiver flag varieties}},
	author={Gu, Wei and Kalashnikov, Elana},
	journal={{arXiv preprint arXiv:2009.02810}},
	year={2020}
}
\end{document}